\documentclass[12pt]{amsart}

\usepackage[leqno]{amsmath}
\usepackage{amssymb}
\usepackage{amsxtra}
\usepackage{amscd}
\usepackage{hyperref}
\usepackage{graphicx, color}
\usepackage{dcpic}

\newtheorem{thm}{Theorem}
\newtheorem{prop}[thm]{Proposition}

\newtheorem{defn}[thm]{Definition}

\newtheorem{lemma}[thm]{Lemma}

\newcommand{\R}{\mathbb{R}}

\newcommand{\F}{\mathbb{F}}
\newcommand{\I}{\ensuremath{\mathbb{I}}}
\newcommand{\N}{\mathbb{N}}
\newcommand{\half}{\mathbb{H}}

\newcommand{\Z}{\mathbb{Z}}
\newcommand{\Zmod}[1]{\Z/{#1}\Z}

\newcommand{\disp}[1]{\displaystyle{#1}}

\newcommand{\graph}[1]{\Gamma_{L}}
\newcommand{\circles}[1]{\ensuremath{\mathrm{\textsc{cir}}(#1)}}

\newcommand{\lra}{\longrightarrow}

\newcommand{\inlinediag}[2][0.33]{\includegraphics[scale=#1]{./images/#2}}

\newcommand{\tangle}[1]{\mathcal{#1}}
\newcommand{\numarcs}[1]{n_{A}(#1)}
\newcommand{\numcircs}[1]{n_{C}(#1)}

\newcommand{\cross}[1]{\ensuremath{\mathrm{\textsc{cr}}(#1)}}
\newcommand{\merge}[1]{\ensuremath{\mathrm{\textsc{merge}}(#1)}}
\newcommand{\leftMerge}[1]{\ensuremath{\overleftarrow{\mathrm{\textsc{merge}}}(#1)}}
\newcommand{\rightMerge}[1]{\ensuremath{\overrightarrow{\mathrm{\textsc{merge}}}(#1)}}
\newcommand{\fission}[1]{\ensuremath{\mathrm{\textsc{divide}}(#1)}}
\newcommand{\leftFission}[1]{\ensuremath{\overleftarrow{\mathrm{\textsc{divide}}}(#1)}}
\newcommand{\rightFission}[1]{\ensuremath{\overrightarrow{\mathrm{\textsc{divide}}}(#1)}}

\newcommand{\dec}[1]{\ensuremath{\mathrm{\textsc{dec}}(#1)}}
\newcommand{\interior}[1]{\ensuremath{\mathrm{\textsc{interior}}(#1)}}
\newcommand{\free}[1]{\ensuremath{\mathrm{\textsc{free}}(#1)}}
\newcommand{\match}[1]{\ensuremath{\mathrm{\textsc{Match}}(#1)}}
\newcommand{\state}[1]{\ensuremath{\mathrm{\textsc{State}}(#1)}}
\newcommand{\actor}[1]{\ensuremath{\mathrm{\textsc{active}}(#1)}}

\newcommand{\resolution}[1]{\ensuremath{\mathrm{\textsc{res}}(#1)}}
\newcommand{\complex}[1]{\ensuremath{[\![#1]\!]}}

\newcommand{\lefty}[1]{\overleftarrow{#1}}
\newcommand{\righty}[1]{\overrightarrow{#1}}

\newcommand{\leftHalf}{\overleftarrow{\half}}
\newcommand{\rightHalf}{\overrightarrow{\half}}

\newcommand{\leftGraph}[1]{\overleftarrow{\Gamma}_{#1}}
\newcommand{\rightGraph}[1]{\overrightarrow{\Gamma}_{#1}}
\newcommand{\bridgeGraph}[1]{\Gamma_{#1}}

\newcommand{\cleaved}[1]{\mathcal{C\!L}_{#1}}
\newcommand{\cleave}[1]{\widehat{\mathcal{C\!L}}_{#1}}
\newcommand{\leftBridges}[1]{\ensuremath{\overleftarrow{\mathrm{\textsc{Br}}}(#1)}}
\newcommand{\rightBridges}[1]{\ensuremath{\overrightarrow{\mathrm{\textsc{Br}}}(#1)}}
\newcommand{\bridges}[1]{\ensuremath{\mathrm{\textsc{Bridge}}(#1)}}

\newcommand{\startCircle}[1]{C_{a}(#1)}
\newcommand{\terminalCircle}[1]{C_{b}(#1)}

\newcommand{\leftnorm}[1]{\lefty{l}(#1)}
 
\newcommand{\rightComplex}[1]{[\![\, #1 \,\rangle\!\!\!\rangle}

\title{A type $D$ structure in Khovanov Homology}
\author{Lawrence P. Roberts}
\email{lproberts@as.ua.edu}
\address{Department of Mathematics\\ 345 Gordon Palmer Hall \\ The University of Alabama \\ Tuscaloosa, Al 35487}
\thanks{This research was supported in part by a grant from the Research Grants Committee of the University of Alabama, Tuscaloosa}

\begin{document}
\begin{abstract}    
We describe the first part of a gluing theory for the bigraded Khovanov homology with $\Z$-coefficients. This part associates a type $D$ structure to a tangle properly embedded in a half-space and proves that the homotopy class of the type $D$ structure is an invariant of the isotopy class of the tangle. The construction is modeled off bordered Heegaard-Floer homology, but uses only combinatorial/diagrammatic methods
\end{abstract}
\maketitle

\section{Introduction}

\noindent Let $\righty{\mathcal{T}} \subset [0,\infty) \times \R^{2}$ be a tangle, properly embedded in a half-space of $\R^{3}$. Suppose that $\partial \righty{\mathcal{T}} \subset \{0\} \times \R \times\{0\}$ is a set of $2k$ points in the $y$-axis, ordered by the linear ordering on the $y$-axis. For example, $\righty{\mathcal{T}}$ might have a diagram $\righty{T}$ such as
$$
\inlinediag[0.7]{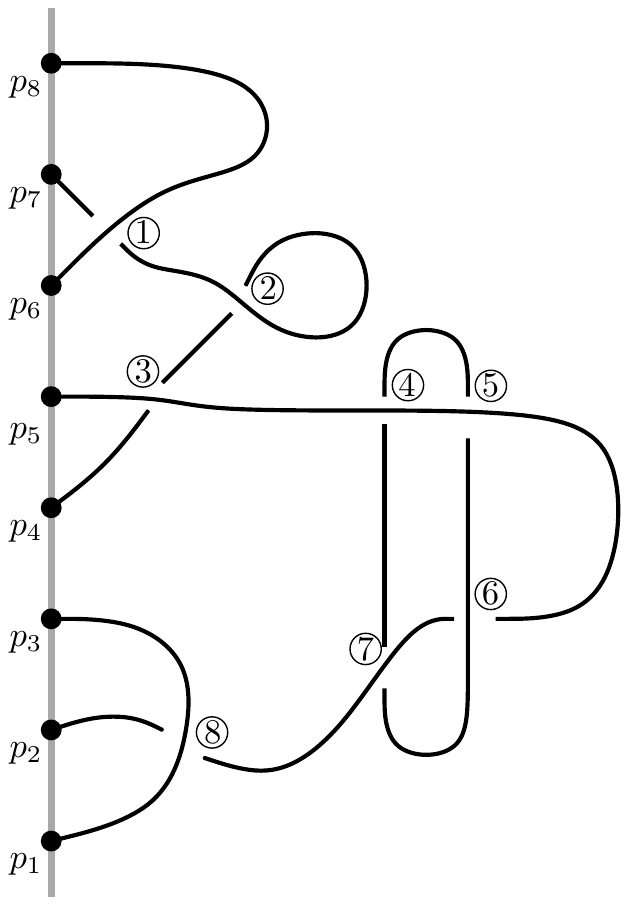}
$$
where we have enumerate the eight crossings, as shown by the circled numbers. We will call these {\em outside} tangles, since the enumeration of $\partial{\righty{\tangle{T}}}$ is compatible with the orientation on the $y$-axis coming from the ``{\em outward} pointing normal first'' convention for orienting $(-\infty,0] \times \R$. Thus $\righty{\mathcal{T}}$ is ``outside'' the fence that is specified by the ordering of its boundary points. \\
\ \\
\noindent There are several different ways to associate a Khovanov type homology to such a tangle. M. Asaeda, J. Przytycki, and A. Sikora take the simplest approach, generalizing the reduced Khovanov homology, in \cite{APS}. In each resolution, they consider those states which assign a $-$ to each arc component, and otherwise follow Khovanov's original specification (although using Viro's conventions). In \cite{PLau}, A. D. Lauda \& H. Pfeiffer increase the sophistication by using an open-closed TQFT to assign decorations to the arc components. This is a natural generalization of the Khovanov Frobenius algebra approach, which corresponds to using a $(1+1)$-TQFT to assign decorations to the circles in a resolution diagram. While both provide invariants of the tangle, computed from a diagram, neither provides a means to compute the original, bigraded Khovanov homology from component tangles (taken with coefficients in $\Z$).  Lauda and Pfeiffer's approach comes close, but will only work with characteristic 2 coefficients, and then only by working through a filtered complexes for each tangle, which arise from a Lee type theory. \\
\ \\
\noindent D. Bar-Natan takes a different approach in \cite{Bar2} by introducing an abstract categorical version of the Khovanov homology. This provides a straightforward way to glue tangles, but only if one is willing to work in the abstract categories of diagrams underlying his approach. Khovanov homology is the result of a functor applied to the complexes in this abstract category, so one could ask if we can find a functor for tangles. This is what Lauda and Pfeiffer did, but as we have noted this approach does not allow one to easily glue the resulting images in a manner corresponding to the gluing of tangles. At the most sophisticated level, Khovanov himself, \cite {Khta}, defines a tangle invariant. This invariant has been much studied and provides a way to construct the Khovanov chain complexes by gluing. However, simplifications of these chain complexes must wait until all the gluing is performed, as it is unclear how exactly states will change prior to the gluing being performed.\\
\ \\
\noindent The current paper begins an effort to resolve these problems. We aim for an invariant of tangles which will allow us to glue tangles to recover the Khovanov homology, while allowing us to simplify as we go, rather than waiting until the whole complex is exposed. We will, in effect, provide a different means for doing the bookkeeping in Khovanov's approach to tangle invariants. At the end we will obtain a gluing theory capable of recovering the module based, bigraded Khovanov homology over $\Z$ from gluing tangles. \\
\ \\
\noindent The key to this approach is to modify the algebraic machinery of bordered Heegaard-Floer homology in \cite{Bor1} to the world of Khovanov homology. We will thus provide a fully combinatorial example of a ``bordered package.'' In addition, this package will be defined for characteristic 0 coefficients, providing a convention for signs which supersedes the characteristic 2 approach in \cite{Bor1}. The package will replicate the algebraic form of \cite{Bor1} and thus allow us to provide a pairing theorem analogous to that in bordered Floer homology.\\
\ \\
\noindent More specifically, in this paper we will construct an free Abelian group $\rightComplex{\righty{T}}$ for each outside tangle. The generators of this group will consist of all diagrams obtained by first
\begin{enumerate}
\item Choosing a planar matching, $\lefty{m}$, of the points in $\partial \righty{T} \subset \{0\} \times \R$, then
\item Choosing for each crossing in $\righty{T}$ a smoothing, encoded either with a $0$ or $1$ (see section \ref{sec:APS} for the convention)
\end{enumerate}
We consider $\lefty{m}$ to lie in $(-\infty,0]\times \R$, and glue it to the diagram found by smoothing $\righty{T}$ at each crossing. We retain the $y$-axis in this diagram, and then assign a decoration in $\{+,-\}$ to each circle in the diagram. For the tangle $\righty{T}$ depicted above, the $10000001$ smoothing of the crossings, and a choice of matching $\lefty{m}$ yields
$$
\inlinediag[0.75]{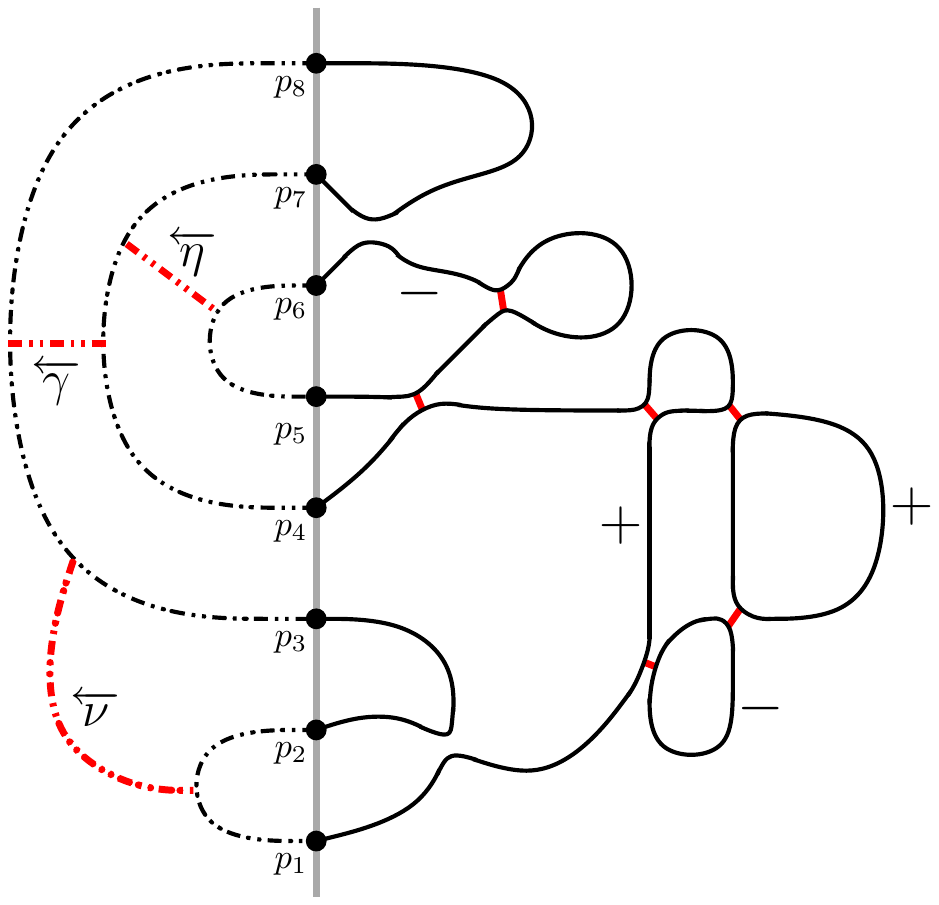}
$$
while the decorations of $\pm$ appear next to each circle. We will denote such a state by $(r,s)$ where $r$ includes both the matching and the smoothing data, and $s$ is the ``sign'' assignment. These states are bigraded, as described in section \ref{sec:APS}. The circles are divided into two classes: those intersecting the $y$-axis, which we will call {\em cleaved} circles, and those which are {\em free}. \\
\ \\
\noindent So far we a replicating the straightforward generalization of Khovanov's tangle invariant to the setting described. We now divert from that approach. To each linearly ordered collection of $2n$-points in the $y$-axis we associate a bigraded algebra $\mathcal{B}\Gamma_{n}$ equipped with a $(1,0)$-differential $d_{\Gamma_{n}}$ which satisfies an appropriate Leibniz identity. The (lengthy) specification of this algebra in terms of generators and relations is given in section \ref{sec:algebra}. For now, it suffices to note that there is an idempotent for each decorated cleaved planar link, without free circles. For example, the following diagram corresponds to an idempotent
$$
\inlinediag[0.6]{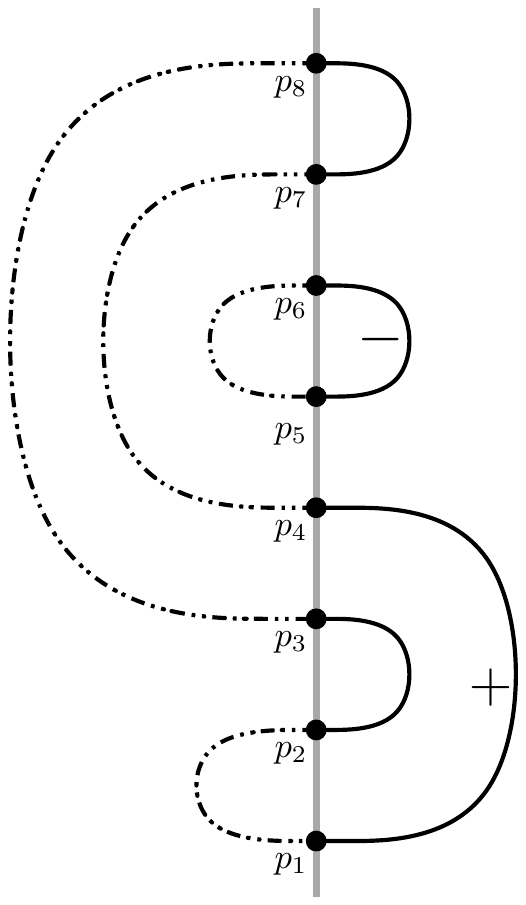}
$$ 
Each state $(r,s)$ is associated to such an idempotent by forgetting the free circles, and considering the remaining diagram up to isotopies which fix the $y$-axis. The depicted idempotent is that corresponding to the $10000001$-resolution diagram depicted above. \\
\ \\
\noindent With $\mathcal{B}\Gamma_{n}$ in hand, we can use the Khovanov procedure to equip $\rightComplex{\righty{T}}$ with a $(1,0)$-module map
$$
\delta_{\righty{T}} : \rightComplex{\righty{T}} \rightarrow \mathcal{B}\Gamma_{n} \otimes_{\mathcal{I}} \rightComplex{\righty{T}} 
$$
The definition is given in section \ref{sec:typeDcom}. We will content ourselves by computing
the image of the $10000001$-resolution. Each red arc in that diagram corresponds to a $0$-resolved crossing. Changing it to a $1$ resolved crossing gives a new diagram. Some of these do not change the corresponding idempotent, and each gives rise to a term in the image of $\delta_{\righty{T}}$. These terms are $\pm 1$ times the following four states
$$
\inlinediag[0.5]{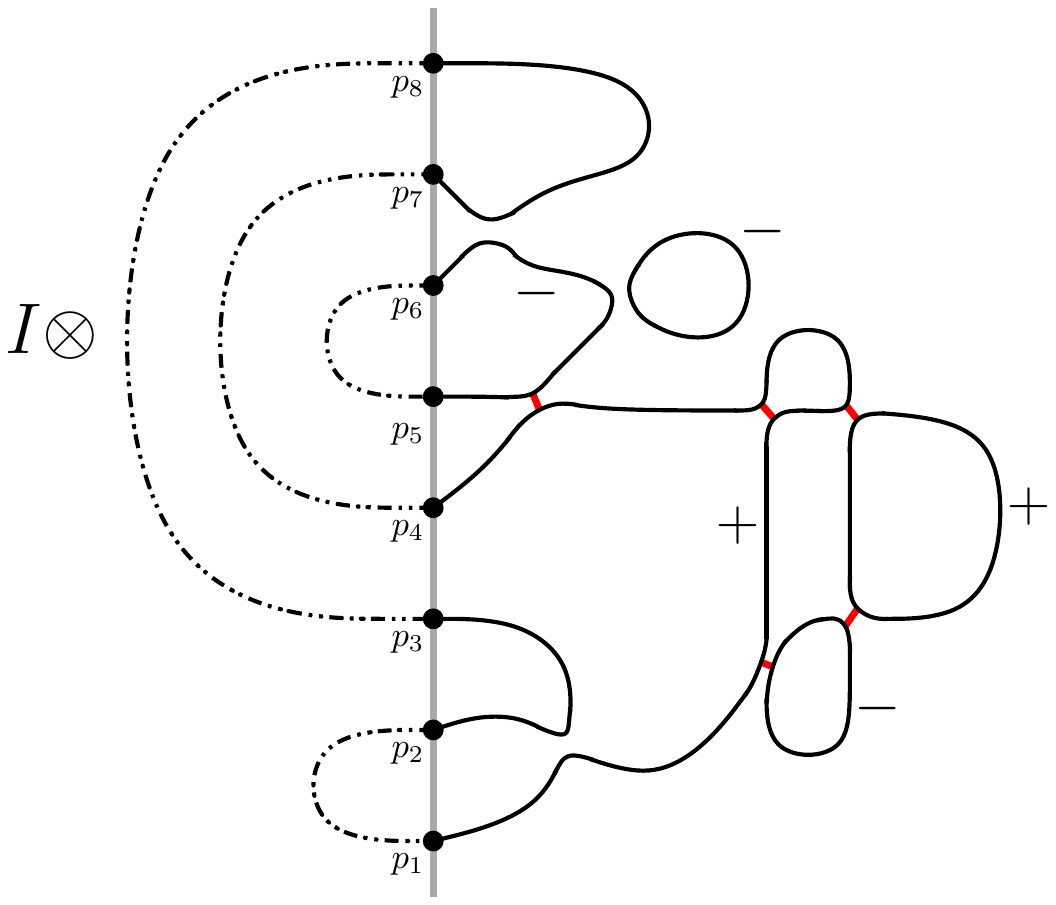} \hspace{1in} \inlinediag[0.5]{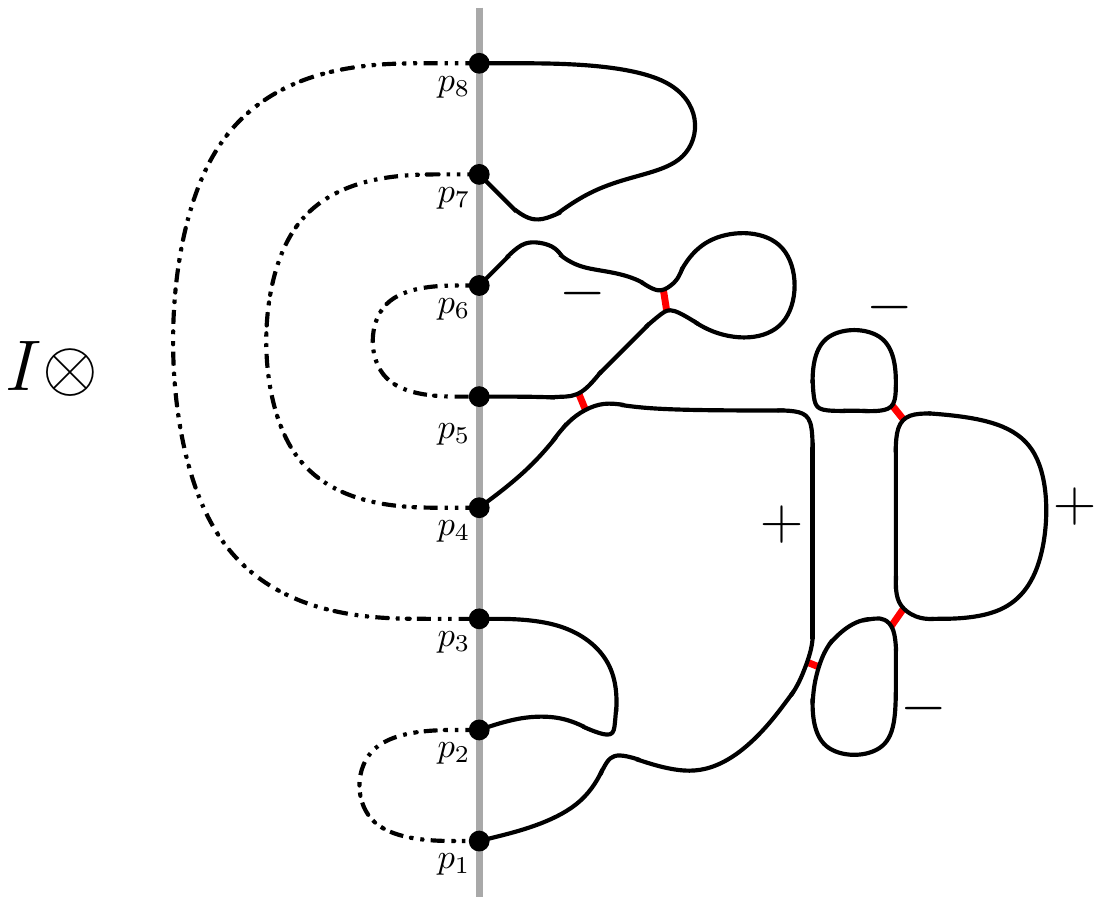}
$$
$$
\inlinediag[0.5]{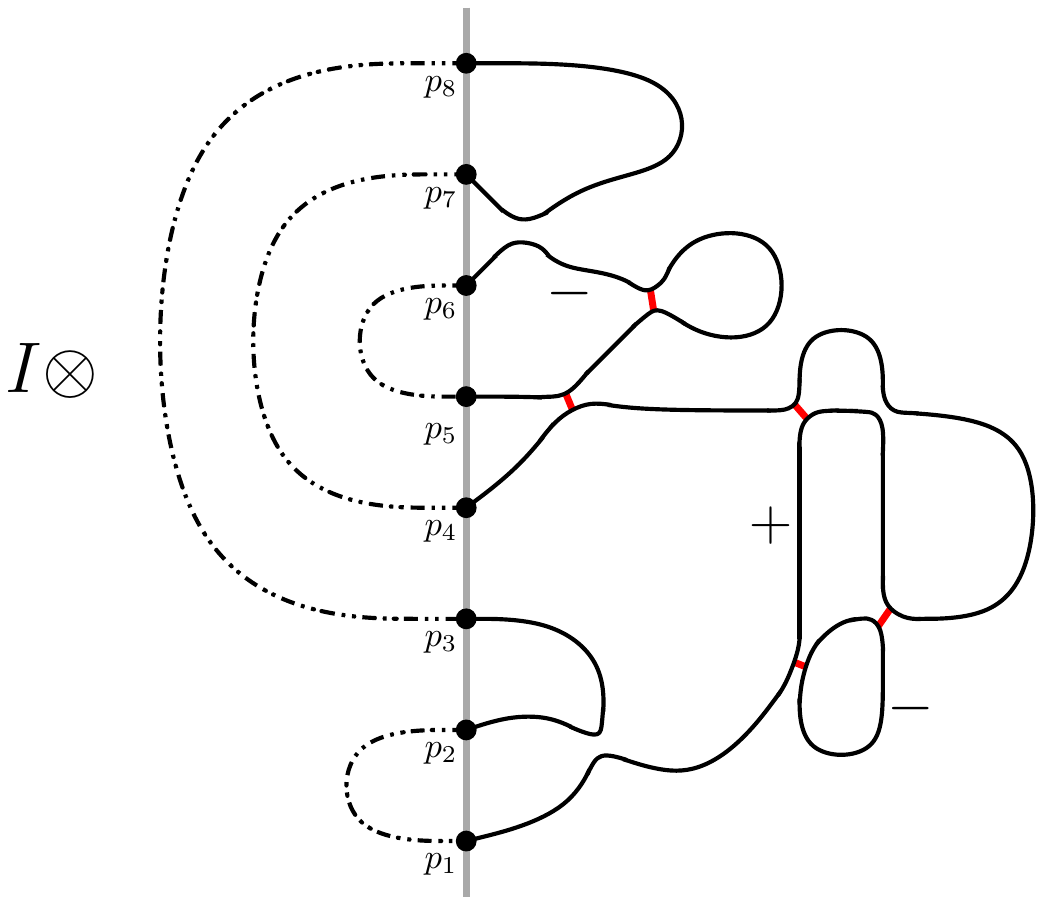} \hspace{1in} \inlinediag[0.5]{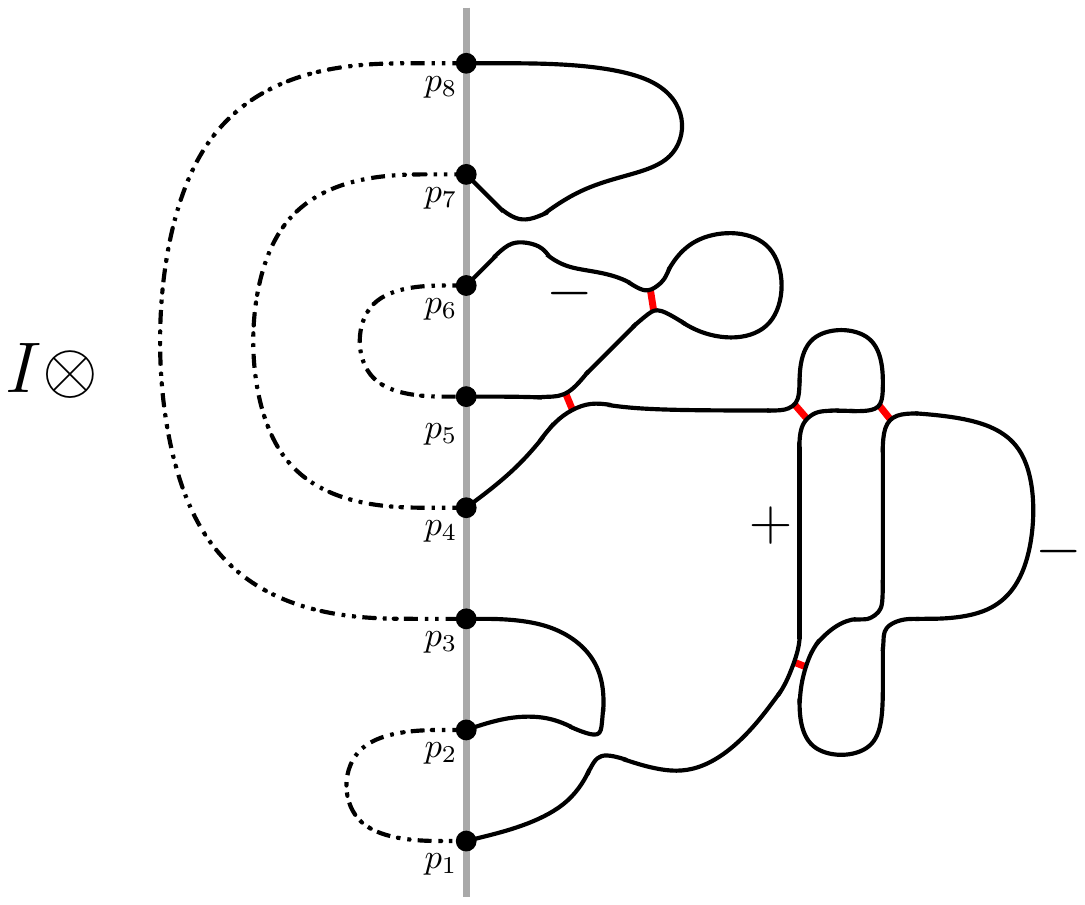} 
$$
The sign convention is defined in section \ref{sec:APS}. The four states correspond, from the left in the first row to the rightmost in the second row, to the smoothings $11000001$, $10010001$, $10001001$, and $10000101$. In the top row, the smoothing changes split a $-$ decorated circle from each of the cleaved circles. In the lower left, a $+$ decorated free circle is merged into a cleaved circle. While in the bottom right two free circles are merged. Each of these terms comes from the original Khovanov differential.  In fact, they are precisely those used in the Asaeda, Przytycki, and Sikora's construction, tensored with the idempotent for the diagram.\\
\ \\
\noindent However, when we split a $+$ cleaved circle in Khovanov's construction we should also get a term where the cleaved circle is decorated with $-$. This corresponds to a different idempotent. To account for the change we use an element of $\mathcal{B}\Gamma_{n}$. If the circle is $D$, this element will be called $\righty{e_{D}}$. There are two such terms in the image of $\delta_{\righty{T}}(r,s)$ for our example -- one for merging the $+$ cleaved circle to a $-$ free circle (smoothing $10000011$) and one from splitting the $+$ cleaved circle (smoothing $10010001$):
$$
\inlinediag[0.5]{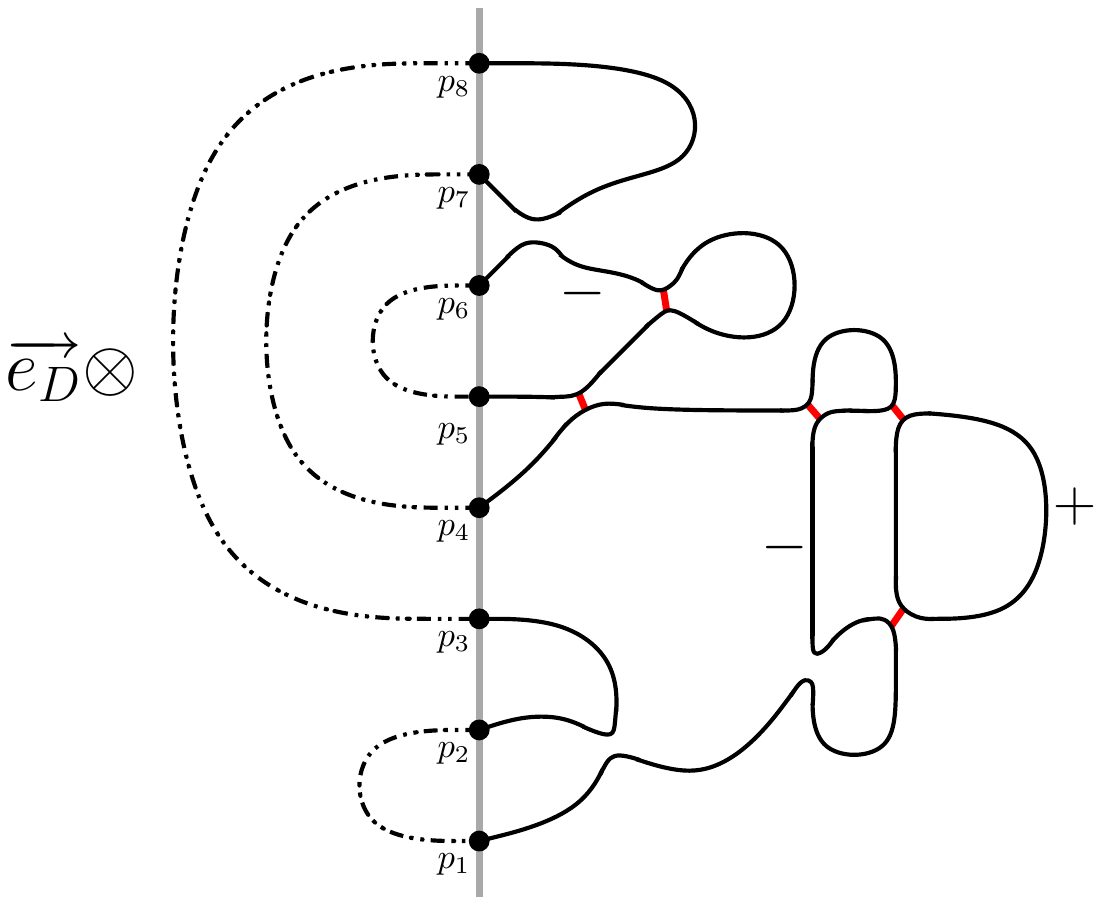} \hspace{1in} \inlinediag[0.5]{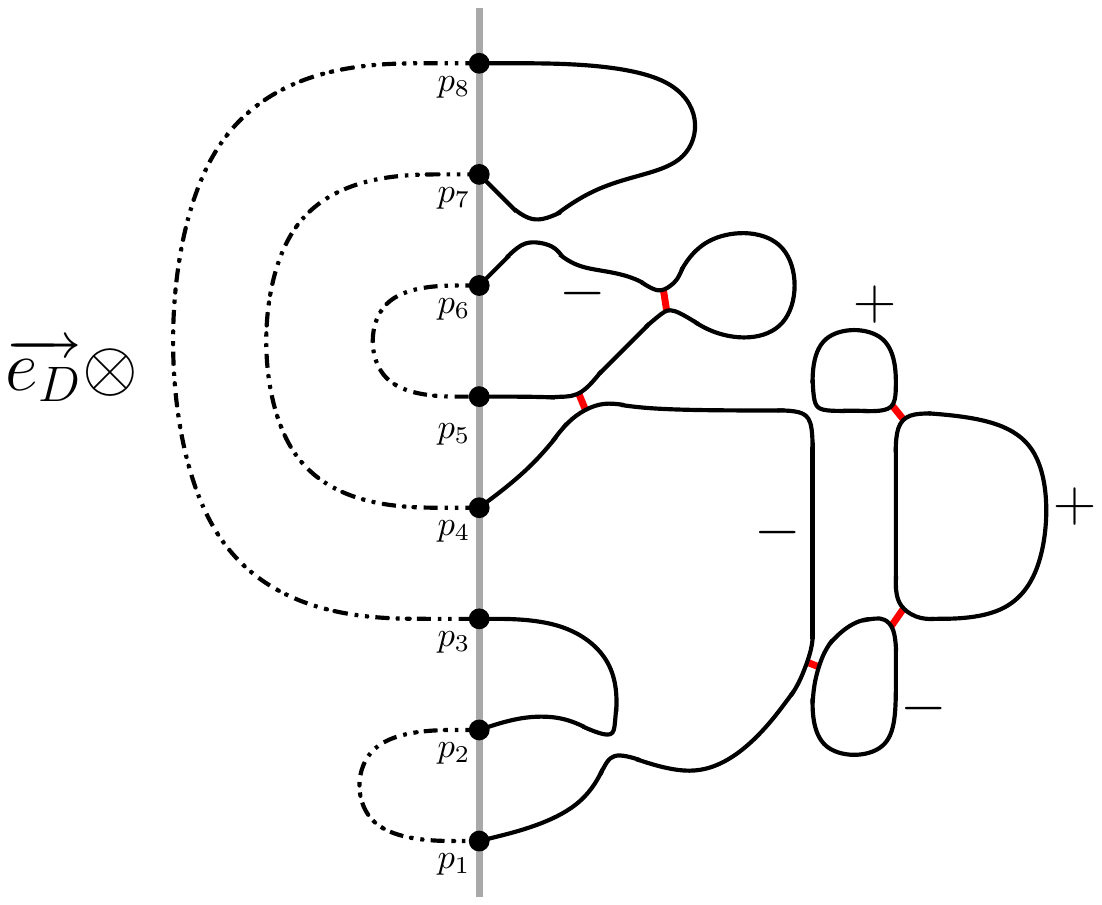}
$$
Again these are shown without the corresponding sign. Note that the smoothing $10010001$ has contributed to both types of terms -- a major source of complexity in the proofs.\\
\ \\
\noindent There is another possible way for the smoothing changes to alter the idempotent corresponding to the diagram. For example, going from $10000001$ to $10100001$ will merge the two cleaved circles. This changes the diagram for the idempotent as well as the decorations. We account for this with another type of algebra element, called $\righty{e_{\delta}}$ here:
$$
\inlinediag[0.5]{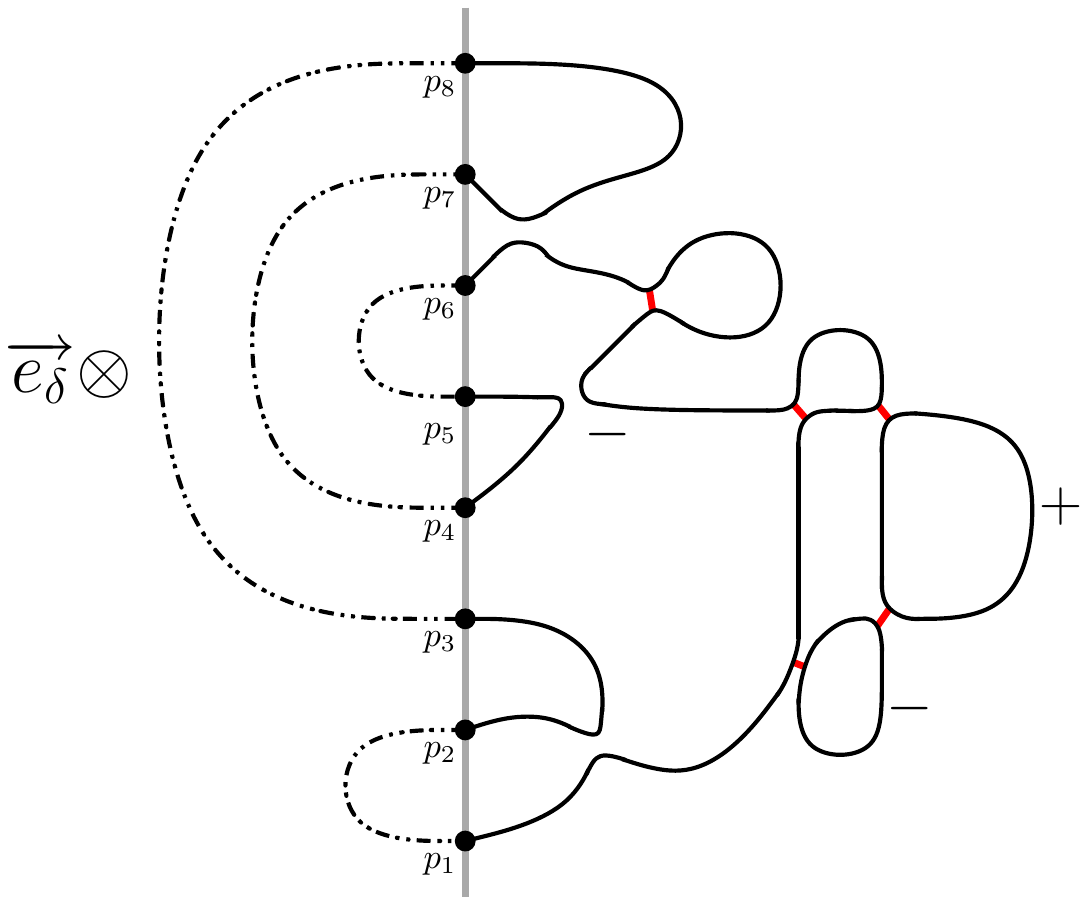}
$$
\ \\
\noindent Each of these changes occurred to the right of the $y$-axis. To the left we only have the matching $\lefty{m}$. In the gluing theory described in the sequel, $\lefty{m}$ will be the sticky end. To glue to any possible diagram, it will need to incorporate all possible changes to the cleaved circles. In particular, we can have surgery along the red arcs $\lefty{\gamma}$, $\lefty{\eta}$ and $\lefty{\nu}$. Surgery on $\lefty{\gamma}$ and $\lefty{\nu}$ each divide the cleaved circle $D$ and thus give rise to two terms each:
$$
\inlinediag[0.5]{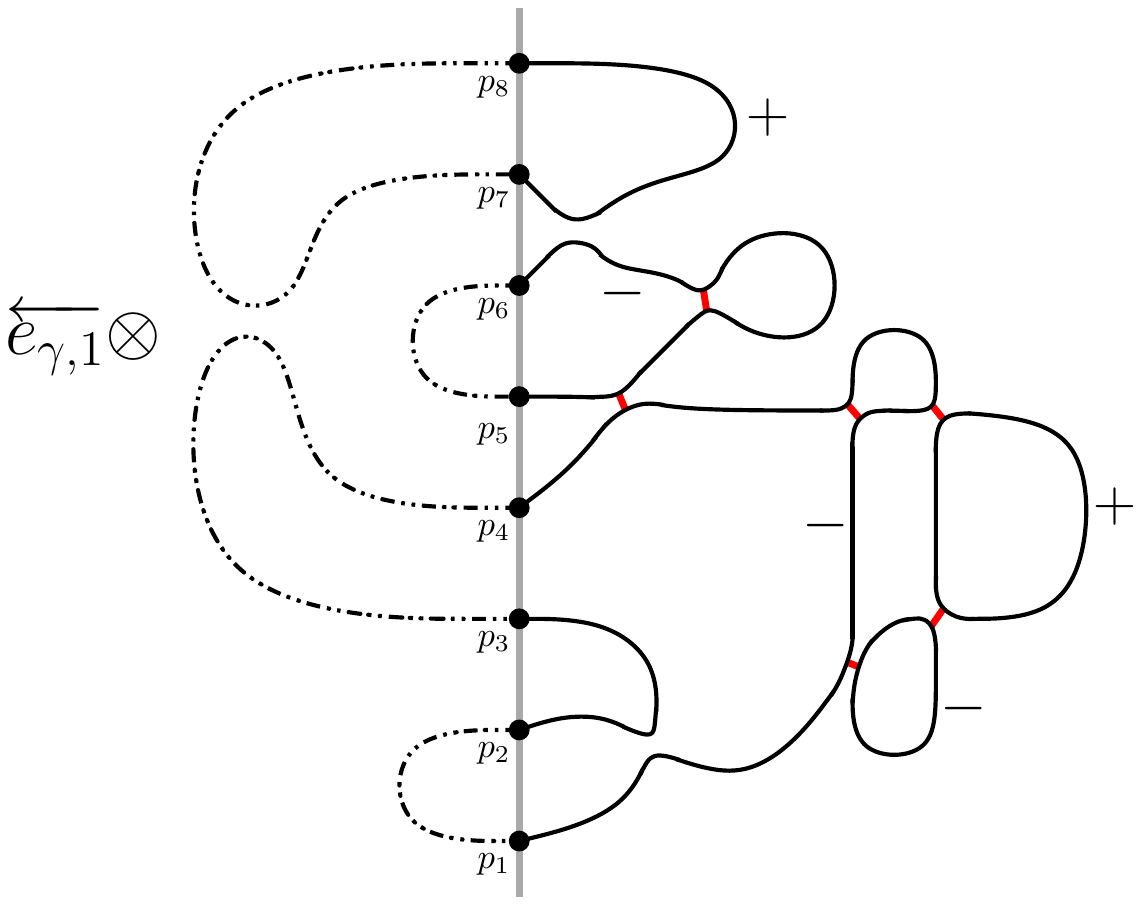} \hspace{1in} \inlinediag[0.5]{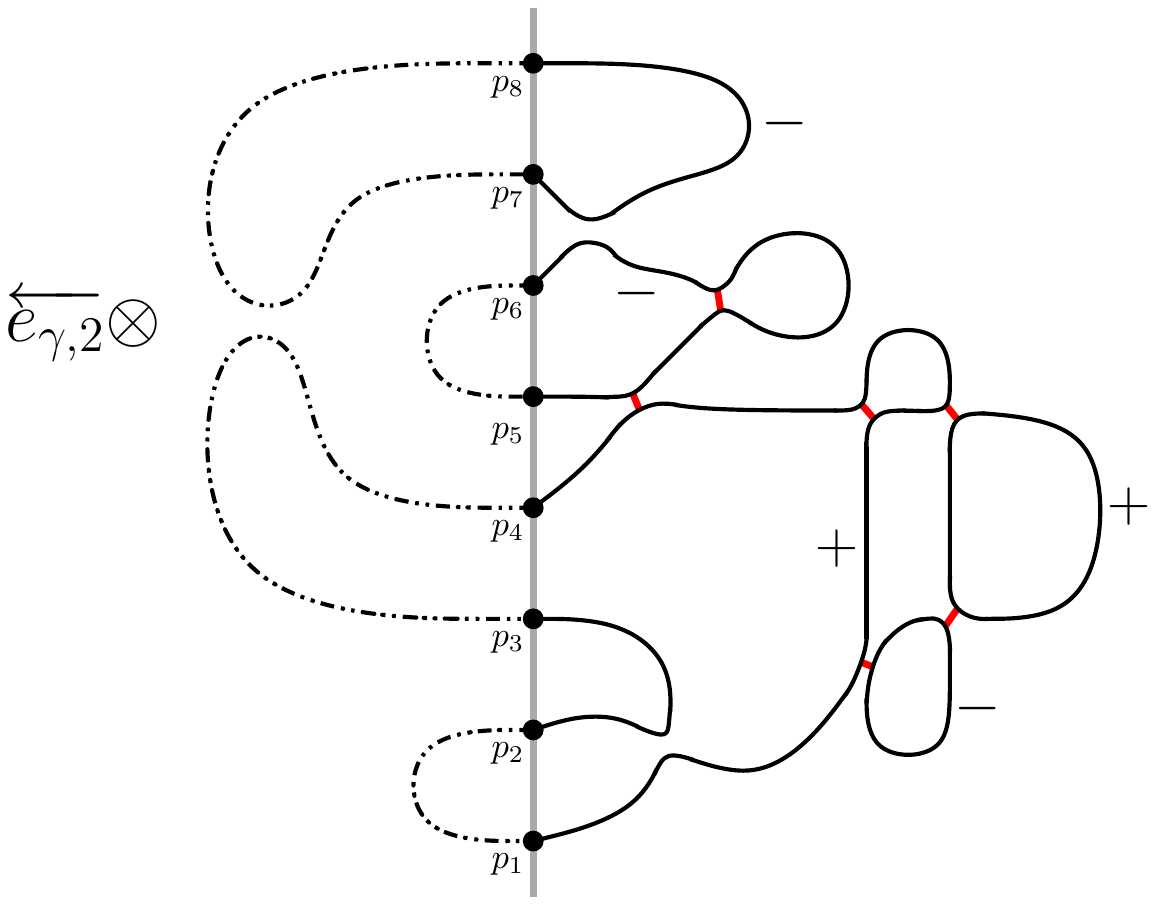}
$$
$$
\inlinediag[0.5]{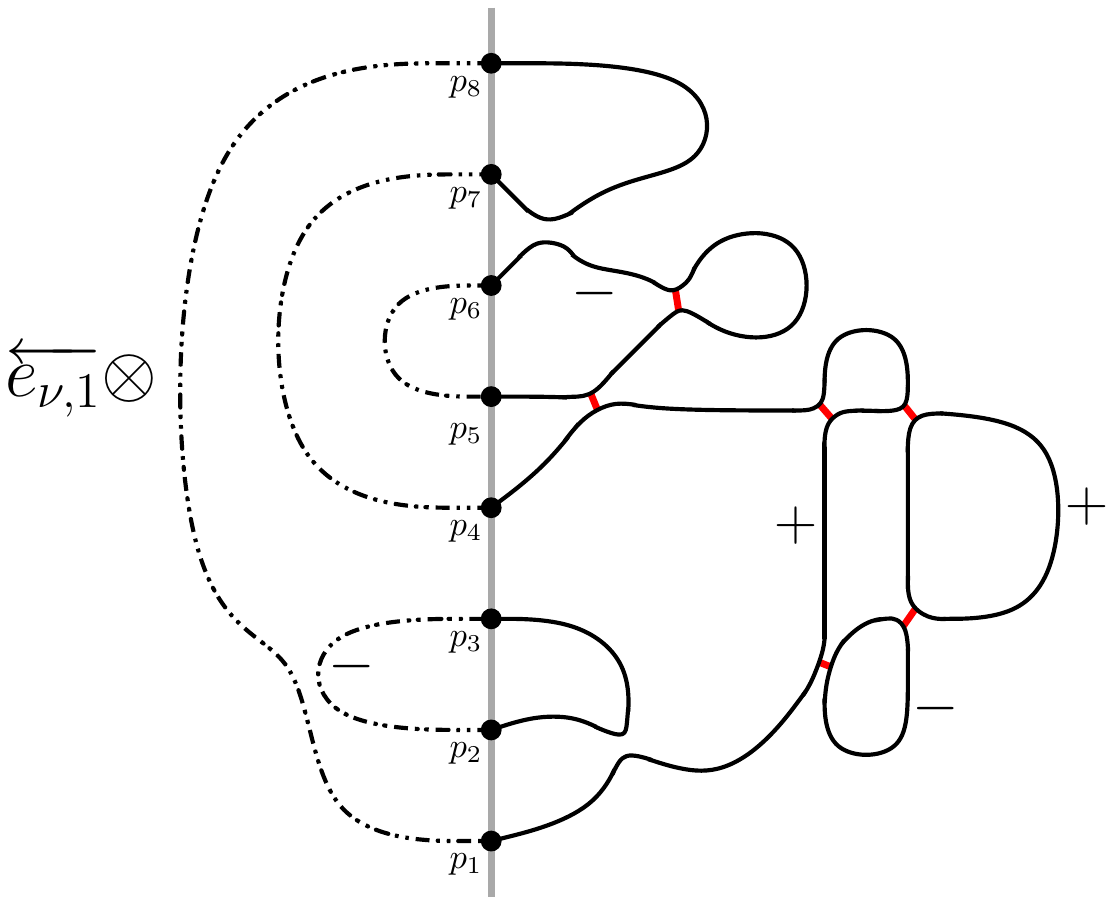} \hspace{1in} \inlinediag[0.5]{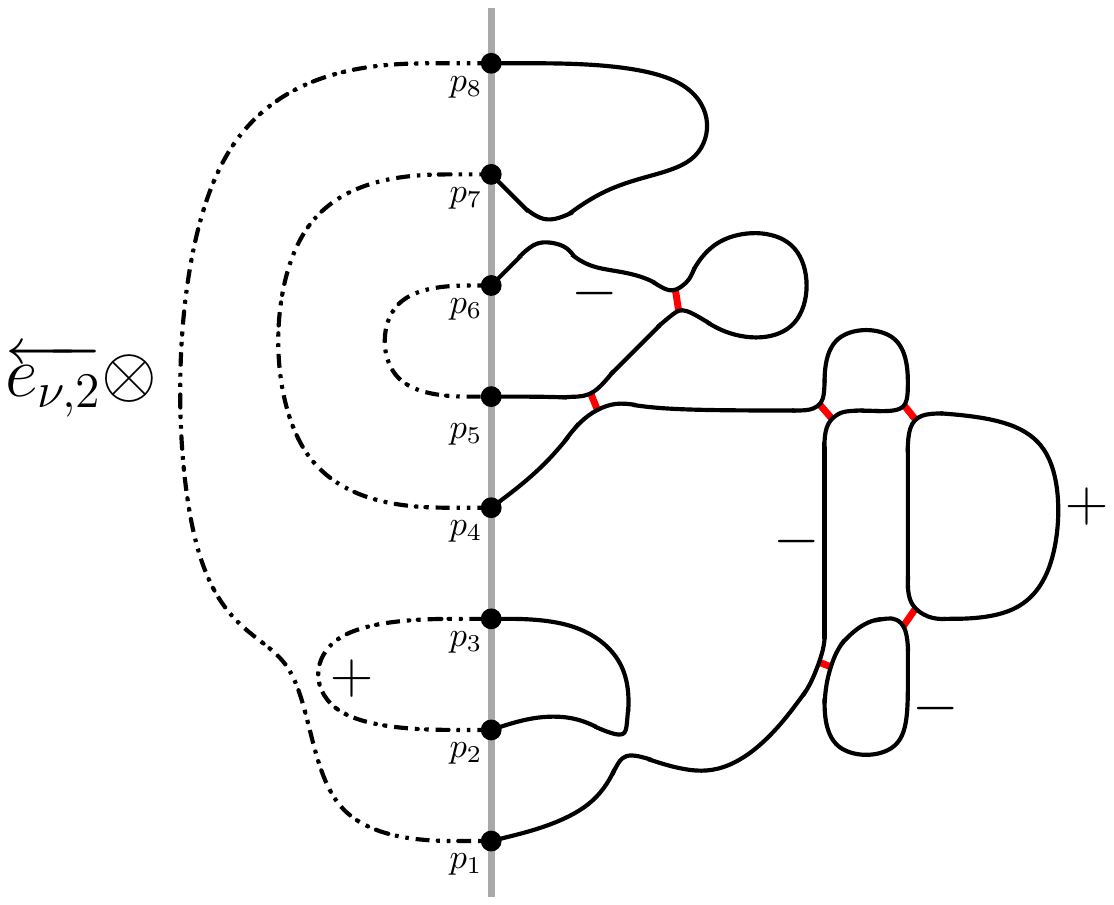}
$$
Note that the choices of smoothings do not change under these surgeries, but that in each pair of terms we have diagrams corresponding to distinct idempotents (so $\lefty{e_{\gamma,1}} \neq \lefty{e_{\gamma,2}}$). Meanwhile $\lefty{\eta}$ merges the two cleaved circles, giving rise to another term
$$
\inlinediag[0.5]{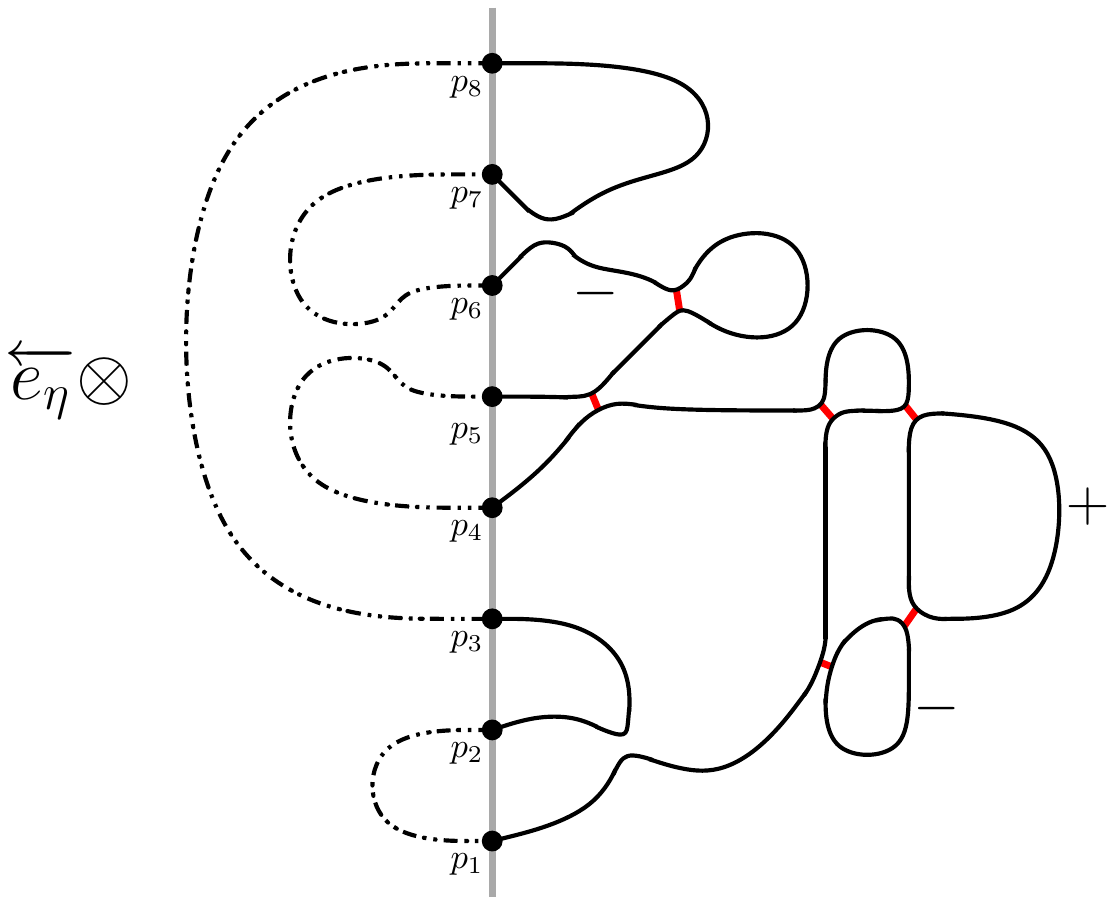}
$$
Finally, operations on the left can change the sign of a cleaved circle. So for each $+$ cleaved circle we obtain one more term, corresponding to flipping the sign to a $-$:
$$
\inlinediag[0.5]{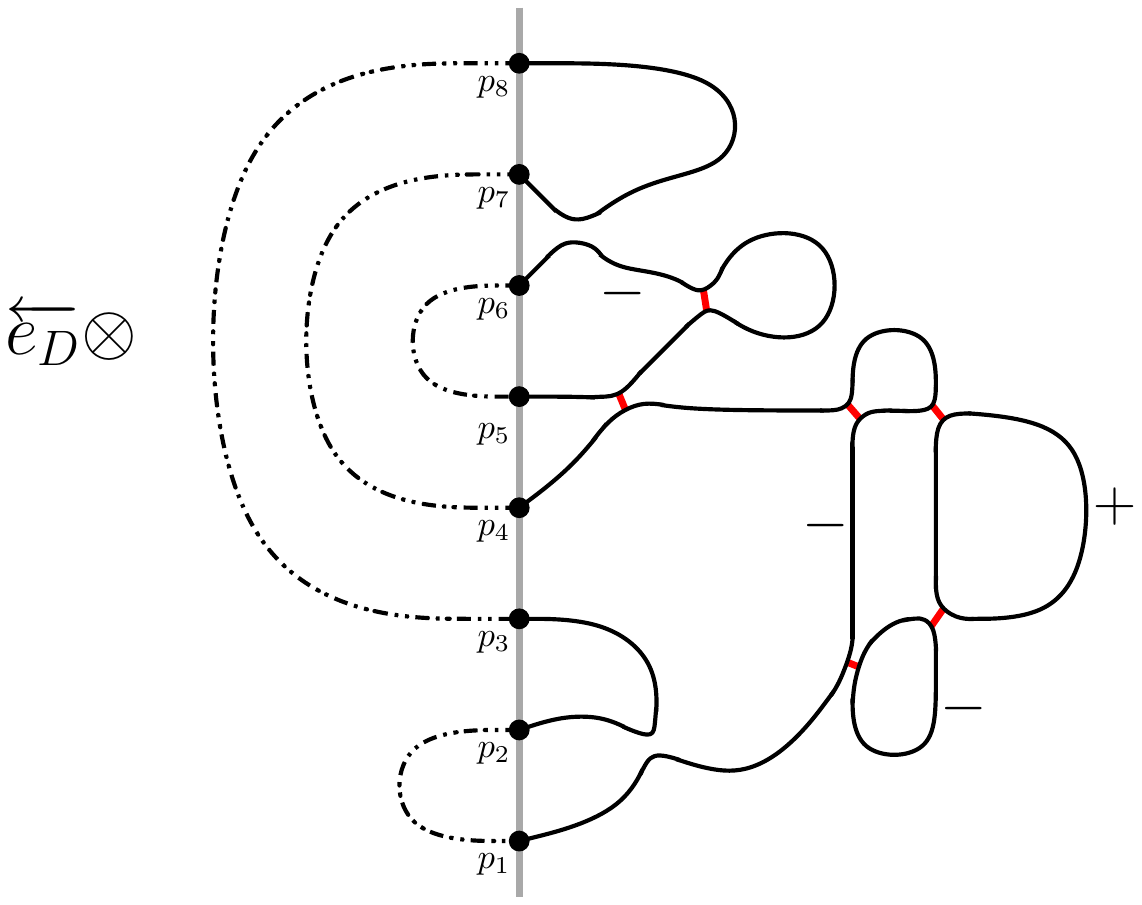}
$$
\ \\
\noindent Despite the evident complexity, the relations in $\mathcal{B}\Gamma_{n}$ allow us to prove that $\delta_{\righty{T}}$ is a type $D$ structure:
$$
(\mu_{\mathcal{B}\bridgeGraph{n}} \otimes \I)\,(\I \otimes \righty{\delta_{T}}) \, \righty{\delta_{T}}  +  (d_{\Gamma_{n}} \otimes |\I|)\,\righty{\delta_{T}} = 0
$$
where $\mu_{\mathcal{B}\bridgeGraph{n}}$ is the product in $\mathcal{B}\Gamma_{n}$ and $|\I|$ is a signed identity on $\rightComplex{\righty{T}}$. The proof is an ugly case-by-case analysis that can be found in section \ref{sec:typeDcom}.\\
\ \\
\noindent Type $D$ structures in \cite{Bor1} admit a notion of homotopy equivalence. We will provide a characteristic 0 version which allows us to prove our main theorem in section \ref{sec:invariance},

\begin{thm}
Let $\righty{\tangle{T}}$ be an outside tangle with diagram $\righty{T}$. The homotopy class of the $D$-structure $\righty{\delta_{T}}: \rightComplex{\righty{T}} \rightarrow \mathcal{B}\bridgeGraph{n} \otimes \rightComplex{\righty{T}}$ is an invariant of the tangle $\righty{\tangle{T}}$.
\end{thm}

\noindent Along the way we will be more precise about the signs, bigrading, and shifts required to obtain this result.\\
\ \\
\noindent Unfortunately, like most Khovanov type theories, the complexity of the computations renders explicit examples prohibitively difficult to include. Nevertheless, one specific case is not much more complicated than calculating the Khovanov homology: tangles with two boundary points. In section \ref{sec:examples} we sketch the calculation for such a tangle coming from the left-handed trefoil, a case which can be done by hand. The Khovanov homology of the left handed trefoil has $2$-torsion. It is interesting to see the interplay between the formalism described here and the occurrence of this torsion. \\
\ \\
\noindent {\bf Note:} After posting this paper to the arXiv, the author was informed by Cotton Seed that he had independently discovered a similar construction of a type D structure in Khovanov homology. The author would like to thank Andy Manion for pointing out an error in the first version of this article.
\section{The algebra from cleaved links}\label{sec:algebra}

\noindent We will define the the algebra $\mathcal{B}\Gamma_{n}$ as the quotient of a quiver algebra on an appropriate digraph. A directed graph $(G,s,t)$, with source and target maps $s, t: E(G) \rightarrow V(G)$, freely generates a category defined by 

\begin{enumerate}
\item the objects are the vertices $v$ of $G$
\item the morphisms from $v_{1}$ to $v_{2}$ consist of all directed paths from $v_{1}$ to $v_{2}$, i.e. all $\rho = e_{1}\ldots e_{m}$ where $e_{i}$ is an edge of $G$, and $s(e_{k+1}) = t(e_k)$ for each $k$.  
\item the identity morphism $\mathrm{I}_{v}$  at a vertex $v$ is the {\em empty} path from $v$ to itself.
\item The composition of morphism $\rho_{1} = e_{1}\ldots e_{m}$ from $v_{1}$ to $v_{2}$ and $\rho_{2}=f_{1}\ldots f_{n}$ from $v_{2}$ to $v_{3}$ is defined to be the concatenation $\rho_{1} |\!| \rho_{2}$. 
\end{enumerate}

\noindent Given a path $\rho = e_{1}\ldots e_{m}$ we define the source and target of the path to be $s(\rho) = s(e_{1})$, $t(\rho) = t(e_{m})$. The length of $\rho$ is $m$, while the identity morphisms have length $0$.\\
\ \\
\noindent This category can further be enhanced to be pre-additive over a field $R$ by taking $R$-linear combinations of the morphisms from $v_{1}$ to $v_{2}$ and extending the composition linearly. This pre-additive category is the {\em quiver category} associated to $G$ and will be denoted $\mathcal{Q}(G)$. When $G$ is finite and has no directed cycles, $\mathcal{Q}G$ is finite-dimensional in the sense that $\mathrm{Mor}(v_{1},v_{2})$ is finite dimensional for each of the finitely many pairs of vertices $v_{1}$ and $v_{2}$. \\ 
\ \\
\noindent If we further extend the composition to a product, by defining the product of paths to be zero unless specified above, we obtain the quiver algebra for $G$. In this case, $I_{v}$ is an idempotent for each $v \in V(G)$.\\
\ \\
\noindent For each $n \in \N$ we will now define a directed graph $\Gamma_{n}$, to be used in this construction, and follow that description by imposing relations on $\mathcal{Q}\Gamma_{n}$ (which the reader can take to be either the pre-additive category or the corresponding quiver algebra).

\subsection{The vertices of $\Gamma_{n}$}

\noindent Let $P_{n}$ be the set of points $p_{1}=(0,1),\ldots, p_{2n}=(0,2n)$ on the $y$-axis of $\R^{2}$, ordered by the second coordinate. We denote the closed half-plane $(-\infty,0] \times \R \subset \R^{2}$ by $\leftHalf$ while  $\rightHalf = [0,\infty) \times \R$. We will take $\R^{2}$ to have its standard orientation. Then the ordering of $P_{n}$ is compatible with that from the boundary orientation of $\leftHalf$ (using the ''outward pointing normal first'' convention) while it is opposite the ordering inherited from the orientation on $\partial \rightHalf$.

\begin{defn}
A {\em $n$-cleaved link} $L$ is any embedding of circles in $\R^{2}$ such that 
\begin{enumerate}
\item the circles of $L$ are disjoint and transverse to the $y$-axis,
\item each point in $P_{n}$ is on a circle in $L$,
\item each circle in $L$ contains at least two points in $P_{n}$
\end{enumerate}
We will denote the circle components of an $n$-cleaved link $L$ by $\circles{L}$.
\end{defn}

\noindent We can think of $L$ as the result of gluing a planar matching of $P_{n}$ in $\leftHalf$ to a planar matching in $\rightHalf$, where

\begin{defn}
A {\em planar matching} $M$ of $P_{n}$ in a half plane $\mathbb{H}$ is an proper embedding of $n$ arcs: $\alpha_{i}: [0,1] \hookrightarrow \mathbb{H}$, with $\alpha_{i}(0), \alpha_{i}(1) \in P_{n}$ for $i=1,\ldots, n$.
\end{defn}
 
\noindent Two planar matchings $M_{1}$ and $M_{2}$ will be considered equivalent if the there is an isotopy of $\mathbb{H}$ which takes $M_{1}$ to $M_{2}$, while pointwise fixing $\partial \mathbb{H}$. The equivalence classes of matchings on $P_{n}$ will be denoted $\match{n}$.

\begin{defn}
The {\em constituents} of an $n$-cleaved link $L$ are the planar matchings which are glued to obtain $L$:
\begin{equation}
\lefty{L} = \leftHalf \cap L  \hspace{1in}
\righty{L} = \rightHalf \cap L 
\end{equation}
\end{defn}
 
\noindent We take two $n$-cleaved links to be equivalent if they are related by isotopy of $\R^{2}$ which pointwise fixes the $y$-axis.

\begin{defn}    
The set of equivalence classes of $n$-cleaved links will be denoted by $\cleave{n}$.  
\end{defn}

\noindent  Thus, equivalence of $n$-cleaved links is the result of gluing the equivalence relation for the constituent planar matchings.\\
\ \\
\noindent A {\em decoration} for an $n$-cleaved link $L$ is a map $\sigma\!: \circles{L} \lra \{+,-\}$. 

\begin{defn}
$\cleaved{n}$ is the set of decorated, $n$-cleaved links:
\begin{equation}
\cleaved{n}= \big\{\,(L,\sigma)\,\big|\,L\in \cleave{n}, \sigma\mathrm{\ is\ a\ decoration\ for\ }L\big\}
\end{equation} 
\end{defn}

\noindent The restriction of a decoration $\sigma$ to $\lefty{L}$ assigns decorations to each arc. We will denote the pair by $(\lefty{L},\sigma)$, and likewise for $\righty{L}$. The sets of these restrictions will be denoted $\lefty{\cleaved{n}}$ and $\righty{\cleaved{n}}$.\\
\ \\

\begin{defn}[Vertices of $\Gamma_{n}$] 
The vertex set of $\Gamma_{n}$ is the set $\cleaved{n}$
\end{defn}

\subsection{The edges of $\Gamma_{n}$}

\noindent The edges of $\Gamma_{n}$ will decompose into a set corresponding to $\leftHalf$ and a set corresponding to $\rightHalf$. We will thus start by giving some definitions based on the planar matchings we might obtain as $\lefty{L}$ and $\righty{L}$. \\

\begin{defn}
A {\em bridge} for a planar matching $M$ is an embedding $\gamma : [0,1] \rightarrow \mathrm{int(}\mathbb{H}\rm{)}$ such that 
\begin{enumerate}
\item $\gamma(0)$ and  $\gamma(1)$ are on distinct {\em arcs} of $M$
\item the image under $\gamma$ of $(0,1)$ is disjoint from $M$ 
\end{enumerate} 
\end{defn}

\noindent The isotopies used to define equivalence for planar matchings $M_{1}$ and $M_{2}$ will take a bridge for $M_{1}$ to a bridge for $M_{2}$. We will thus consider a bridge $\gamma_{1}$ for $M_{1}$ to be equivalent to a bridge $\gamma_{2}$ for $M_{2}$ if there is a planar isotopy, fixing $\partial \mathbb{H}$, which takes $\gamma_{1}$ into $\gamma_{2}$ and $M_{1}$ into $M_{2}$. When $M = M_{1} = M_{2}$ this will include sliding the feet of the bridge along the arcs they abut in $M$. 

\begin{defn}
The equivalence classes of bridges for a planar tangle $T$, without closed components, will be denoted $\bridges{T}$.
\end{defn}

\noindent Given a bridge $\gamma$ for a planar matching $M$, we can construct a new planar matching by surgery along $\gamma$.



\begin{defn}
Let $T$ be an equivalence class of planar matchings with a bridge $\gamma$. By $T_{\gamma}$ we will mean the equivalence class of planar matching found from surgery along $\gamma$. 
\end{defn}

\noindent $T_{\gamma}$ has a special bridge $\gamma^{\dagger}$ which is the image of the co-core of the surgery.\\
\ \\
\noindent We now lift these concepts to the links $L \in \cleave{n}$. However, by ``arc'' we will mean an arc in one of the two planar matchings, while ``circle'' will mean a component of $L$ considered as a planar link.  

\begin{prop}
Given any bridge $\gamma$ for $L$, $\bridges{L}\backslash\{\gamma\}$ can be decomposed as a disjoint union $B_{\pitchfork}(L,\gamma) \cup B_{|\,|}(L,\gamma)$, where
\begin{enumerate}
\item $B_{\pitchfork}(T,\gamma)$ consists of the classes of bridges all of whose representatives intersect $\gamma$
\item $B_{|\,|}(L,\gamma)$ consists of the classes of bridges containing a representative which does not intersect $\gamma$,
\end{enumerate}  
Furthermore, we can divide $B_{|\,|}(L,\gamma)$ into the disjoint union $B_{d}(L,\gamma) \cup B_{s}(L,\gamma) \cup B_{o}(L,\gamma)$ where
\begin{enumerate}
\item $B_{d}(L,\gamma)$ consists of those bridges neither of whose ends is on an arc with $\gamma$,
\item $B_{s}(L,\gamma)$ consists of those bridges with a single end on the same arc as $\gamma$ and lying on the {\em same} side of the arc as $\gamma$
\item $B_{o}(L,\gamma)$ consists of those bridges with a single end on the same arc as $\gamma$ and lying on the {\em opposite} side of the arc as $\gamma$
\end{enumerate}
\end{prop}

\noindent If $\eta \in B_{\circ}(L,\gamma)$ where $\circ$ represents a specific choice of one of the above then $\gamma \in B_{\circ}(L,\eta)$. Furthermore, if $\delta$ has a different location than $\gamma$ then $\delta \in B_{d}(L,\gamma)$. We consider how these sets change under surgery on $\gamma$.

\begin{prop} 
Surgery on $\gamma$ induces an identification $B_{d}(L,\gamma)$ with $B_{d}(L_{\gamma}, \gamma^{\dagger})$ and a $2:1$ map $B_{s}(L,\gamma) \lra B_{o}(L_{\gamma},\gamma^{\dagger})$. 
\end{prop}

\noindent Dually there is a $2:1$ map $B_{s}(L_{\gamma}, \gamma^{\dagger}) \lra B_{o}(L,\gamma)$.\\
\ \\
\noindent{\bf Proof:} Let $\eta \in B_{|\,|}(L,\gamma)$. Pick a representative arc for $\eta$ which does not intersect the representative arc for $\gamma$. Then $\eta$ also represents a bridge in $B_{|\,|}(L_{\gamma}, \gamma^{\dagger})$, and vice-versa. If $\eta \in B_{d}(L,\gamma)$, any isotopy of the representaive arc occurs in a region disjoint from $\gamma$ and its endpoints, since the isotopy will occur along arcs disjoint from those intersecting $\gamma$. This isotopy also survives into $(L_{\gamma}, \gamma^{\dagger})$. Reversing this construction for $B_{d}(L_{\gamma}, \gamma^{\dagger})$ proves the identification. Note that an isotopy of $\eta \in B_{s}(L,\gamma)$ missing $\gamma$ can likewise be pushed forward. However, for each $\eta$ we can slide $\eta$ over $\gamma$ to get another bridge $\eta'\in B_{s}(L,\gamma)$. In $L_{\gamma}$ $\eta \simeq \eta'$ and both are on the opposite side of $\gamma^{\dagger}$. By looking at a local model, this is the only type of collision, so the map is $2:1$ on $B_{s}(L,\gamma)$. We can apply the same argument to $(L_{\gamma},\gamma^{\dagger})$ to obtain the $2:1$ map in the other direction. $\Diamond$. \\
\ \\
\noindent Each of these notions: bridges, surgery, etc, can be lifted from the planar matchings $\lefty{L}$ and $\righty{L}$ to  the $n$-cleaved link $L$, as long as we do not alter the $y$-axis. {\em We will do this without comment, but retain the ``arrow'' notation to distinguish restriction to each of the two sides of the cutting line.} Thus a bridge for $\lefty{L}$ will be interpreted as a bridge for $L$, and the equivalence classes of bridges for $\lefty{L}$ \rm{(}$\righty{L}$\rm{)} will be denoted by $\leftBridges{L}$ \rm{(}$\rightBridges{L}$\rm{)}. The bridges for $L$ will be one or other of these types: $\bridges{L} = \leftBridges{L} \cup \rightBridges{L}$.\\
\ \\
\noindent We now introduce some notation for handling the bridges for a cleaved link $L$, relative to the circles in $L$:

\begin{defn}
Let $\gamma \in \bridges{L}$ for an $n$-cleaved link $L$, then
\begin{enumerate}
\item $L_{\gamma}$ is the $n$-cleaved link found by surgering the appropriate constituent of $L$,
\item The support of $\gamma$ is the set of three circles in $L$ and $L_{\gamma}$ which contain the feet of $\gamma$ and $\gamma^{\dagger}$
\item $\gamma$ is in $\merge{L}$ if there are two distinct circles $\{\startCircle{\gamma}, \terminalCircle{\gamma}\}$ containing
$\gamma(0)$ and $\gamma(1)$. In this case, $C_{\gamma}$ is the circle in $\circles{L_{\gamma}}$ which contains both feet of $\gamma^{\dagger}$.
\item $\gamma$ is in $\fission{L}$ if both feet of $\gamma$ are on the same circle $C$ of $L$. In this case, $C^{a}_{\gamma}$ and $C^{b}_{\gamma}$ are the circles in $\circles{L_{\gamma}}$ which contain the feet of $\gamma^{\dagger}$.   
\end{enumerate}
\end{defn}

\noindent We are now in a position to define the edges of  $\bridgeGraph{n}$. The edges come in two types; we start by defining the {\em bridge} edges in $\bridgeGraph{n}$:
   
\begin{defn}[Bridge edges in $\bridgeGraph{n}$]\label{def:bridgeGraph}
Let $(L,\sigma) \in \cleaved{n}$. There is an edge $(L,\sigma) \longrightarrow (L_{\gamma}, \sigma_{\gamma})$ of $\bridgeGraph{n}$ for each $\gamma \in \bridges{L}$ whenever $\sigma(C)=\sigma_{\gamma}(C)$ for $C \in \circles{L}$ not in the support of $\gamma$, and for circles in the the support of $\gamma$ the decorations match one of the following cases:
\begin{enumerate}
\item\label{bridge1} when  $\gamma \in \merge{L}$ and $\sigma$ and $\sigma_{\gamma}$ restrict to the support of $\gamma$ as one of 
\begin{equation}
\begin{array}{llcl}
\sigma(\startCircle{\gamma})=+ & \sigma(\terminalCircle{\gamma})= + &\hspace{1in} &\sigma_{\gamma}(C_{\gamma})= + \\
\sigma(\startCircle{\gamma})=- & \sigma(\terminalCircle{\gamma})= + &\  &\sigma_{\gamma}(C_{\gamma})= - \\
\sigma(\startCircle{\gamma})=+ & \sigma(\terminalCircle{\gamma})= - &\ &\sigma_{\gamma}(C_{\gamma})= - \\
\end{array}
\end{equation} 
\item\label{bridge2} when $\gamma \in \fission{L}$, $C \in \circles{L}$ is the circle containing both feet of $\gamma$, and
\begin{enumerate}
\item  $\sigma(C) = +$,  if $\sigma$ and $\sigma_{\gamma}$ restrict to the support of $\gamma$ as either of
\begin{equation}
\begin{array}{llcl}
\sigma(C)=+ &\hspace{1in} & \sigma_{\gamma}(C^{a}_{\gamma})= + & \sigma_{\gamma}(C^{b}_{\gamma})= - \\
\sigma(C)=+ &\ & \sigma_{\gamma}(C^{a}_{\gamma})= - & \sigma_{\gamma}(C^{b}_{\gamma})= + \\
\end{array}
\end{equation}
\item $\sigma(C) = -$, if $\sigma$ and $\sigma_{\gamma}$ restrict to the support of $\gamma$ as 
\begin{equation}
\begin{array}{llcl}
\sigma(C)=- & \hspace{1in} & \sigma^{-}_{\gamma}(C^{a}_{\gamma})= - & \sigma^{-}_{\gamma}(C^{b}_{\gamma})= - \\
\end{array}
\end{equation} 
\end{enumerate}
\end{enumerate}
\end{defn}

\noindent The second type of edge in $\bridgeGraph{n}$ correspond to changing the decoration on a single circle, and will be called {\em decoration} edges:

\begin{defn}[Decoration edges in $\bridgeGraph{n}$]
 For each circle $C \in \circles{L}$ with $\sigma(C)=+$ there are {\em two} distinct edges $(L,\sigma) \longrightarrow (L,\sigma_{C})$, where $\sigma_{C}$ is the decoration on $L$ with $\sigma_{C}(C) = -$ and  $\sigma_{C}(C')=\sigma(C')$ for $C' \in \circles{L} \backslash \{C\}$. 
\end{defn} 

\noindent The circle whose decoration changes along a decoration edge will be called the {\em support} of that edge.\\
\ \\
\noindent  The bridge edges naturally partition into those coming from $\lefty{L}$ and those from $\righty{L}$. We reflect this partition in the decoration edges by (arbitrarily) assigning one of the decoration edges for each $C \in \circles{L}, \sigma(C) = +$ to $\lefty{L}$ and one to $\righty{L}$. $\leftGraph{n}$ is the subgraph of $\bridgeGraph{n}$ with the same vertices but whose edges correspond to $\lefty{L}$. Likewise, $\rightGraph{n}$ is the subgraph whose edges correspond to bridges in $\righty{L}$. The {\em location} of an edge is the subgraph $\leftGraph{n}$ or $\rightGraph{n}$ which contains the edge.   \\
\ \\
\noindent We have completed the construction of the directed graph needed to specify the quiver algebra. 
 
\begin{defn}
$\mathcal{Q}\bridgeGraph{n}$ is the quiver category induced from the graph $\bridgeGraph{n}$ 
\end{defn}

\begin{prop}\label{prop:noCycle}
$\mathcal{Q}\bridgeGraph{n}$ is finite dimensional
\end{prop}

\noindent{\bf Proof:} We show that $\bridgeGraph{n}$ contains no directed cycle. For each $(L,\sigma) \in \cleaved{n}$ define 
\begin{equation}
\iota(L,\sigma) = \#\big\{C \in \circles{L} \big| \sigma(C) = +\big\} - \#\big\{C \in \circles{L} \big| \sigma(C) = -\big\}
\end{equation}
\noindent For any bridge edge $(L,\sigma) \rightarrow (L_{\gamma},s^{i}_{\gamma})$,  $\iota(L_{\gamma},s^{i}_{\gamma}) = \iota(L,s) - 1$. Furthermore, for a decoration edge $(L,s) \rightarrow (L,s_{C})$ we have $\iota(L,s_{C}) = \iota(L,s) - 2$. Consequently, a directed path in $\bridgeGraph{n}$ strictly decreases the value of $\iota$, so there can be no directed cycle. $\Diamond$\\
\ \\
\noindent It will be convenient to have a shortened notation for the algebra element corresponding to each edge in $\bridgeGraph{n}$:
 
\begin{enumerate}
\item {\bf For a bridge edge:} An edge
$$
(L,\sigma) \longrightarrow (L_{\gamma}, \sigma')
$$
corresponding to a bridge $\gamma \in \bridges{L}$ will be labeled $e_{(\gamma; \sigma , \sigma')}$ (or $e_{\gamma}$ if the context is clear). When we wish to emphasize that $\gamma \in \leftBridges{L}$ we will write $\lefty{e}_{(\gamma; \sigma, \sigma')}$ and likewise for $\rightBridges{L}$. When we wish to emphasize that this is a merge we will use $m_{(\gamma; \sigma, \sigma')}$ (and $f_{(\gamma; \sigma, \sigma')}$ for a division). 
\item {\bf For a decoration edge:} When $\sigma(C) = +$,  $\lefty{e_{C}}$ will denote the decoration edge
$$
(L,\sigma) \longrightarrow (L, \sigma_{C})
$$
that been assigned to $\lefty{L}$. For a decoration edge assigned to $\righty{L}$ we will use $\righty{e_{C}}$. 
\end{enumerate} 

\noindent {\bf Cautionary note:} The label drops the link $L$ from the notation since it is the domain of $\sigma$. However, even though the same bridge $\gamma$ may appear to be in different $L$'s when we draw our diagrams, they have not been identified, and the corresponding edges in $\bridgeGraph{n}$ should be considered distinct, as should the algebra elements. \\
\ \\
\noindent We  equip $\mathcal{Q}\bridgeGraph{n}$ with some additional structure, which will become important shortly

\begin{defn}
The bigrading on $\mathcal{Q}\bridgeGraph{n}$ is the bigrading of paths in $\bridgeGraph{n}$ induced by homomorphically extending the following bigrading of the edges: 
$$
\begin{array}{lcl}
\mathrm{I}_{(L,\sigma)} & \longrightarrow & (0,0)\\
\righty{e_{C}} & \longrightarrow & (0,-1)\\
\lefty{e_{C}} & \longrightarrow & (1,1)\\
\righty{e_{\gamma}} &\longrightarrow & (0,-1/2)\\
\lefty{e_{\gamma}} &\longrightarrow & (1,1/2)\\
\end{array}
$$
The first entry of the bigrading on the path $\alpha$ will be denoted by  $\leftnorm{\alpha}$, while the second element will be denoted $q(\alpha)$.
\end{defn}

\noindent With this choice, elements corresponding to $\rightHalf$ will act as even elements for the $\Z/2\Z$-grading from $\leftnorm{\alpha}$, while elements in $\leftHalf$ will be {\em odd} elements. This distinction will play an important role as we impose relations on $\mathcal{Q}\bridgeGraph{n}$ since the relations will supply a form of graded commutativity based on this $\Z/2\Z$-grading.

\subsection{The cleaver category}

We now describe a collection of relations we will impose on $\mathcal{Q}\bridgeGraph{n}$. The pre-additive category we obtain by quotienting by these relation we will denote $\mathcal{B}\bridgeGraph{n}$ and call the {\em cleaver category}. As each of these relation will be homogeneous for the bigrading, $\mathcal{B}\bridgeGraph{n}$ will also be bigraded.\\
\ \\
\noindent First, we will identify a number of different types of squares in $\bridgeGraph{n}$ of the form
\begin{equation}
\begin{CD}
(L,\sigma) @>e_{\alpha}>> (L_{\alpha},\sigma_{\alpha}) \\
@Ve_{\beta}VV 				@VVe_{\beta'}V\\
(L_{\beta}, \sigma_{\beta}) @>e_{\alpha'}>> (L_{2}, \sigma_{2})\\
\end{CD}
\end{equation}
For each we will impose the  relations
\begin{equation}
e_{\alpha}e_{\beta'} = (-1)^{\leftnorm{e_{\alpha}}\leftnorm{e_{\beta}}}e_{\beta}e_{\alpha'} 
\end{equation}

\noindent{\bf Relations from disjoint support:} Suppose we have two edges $(L, \sigma) \stackrel{e_{\alpha}}{\longrightarrow} (L_{\alpha}, \sigma_{\alpha})$ and $(L, \sigma) \stackrel{e_{\beta}}{\longrightarrow} (L_{\beta}, \sigma_{\beta})$ with the same source and disjoint supports. We consider three possibilities:
\begin{enumerate}

\item Both edges $e_{\alpha}$ and $e_{\beta}$ are bridge edges for the different equivalence casses of bridges $\gamma$, $\eta$ respectively. Since the
$\gamma$ and $\eta$ are distinct $\eta \in \bridges{L_{\gamma}}$ and $\gamma \in \bridges{L_{\eta}}$. Furthermore, $L_{\gamma,\eta} = L_{\eta,\gamma}$. In addition, since the support of $\gamma$ and $\eta$ are disjoint $\sigma_{\gamma}$ assigns the same decorations to the support of $\eta$ as does $\sigma$, so there is an edge $e_{\beta'} : (L_{\gamma},\sigma_{\gamma}) \longrightarrow (L_{\gamma,\eta}, \sigma_{\gamma,\eta})$
where $\sigma_{\gamma,\eta}$ assigns the same decorations to circles in the support of $\gamma$ as $\sigma_{\gamma}$ and to circles in the support of $\eta$ as $\sigma_{\eta}$ (and $\sigma$ for every other circle). The same argument applies with $\eta$ and $\sigma$ in the opposite order, and it is easy to see that $\sigma_{\gamma,\eta} = \sigma_{\eta,\gamma}$. Consequently, there is a square
\begin{equation}
\begin{CD}
(L,\sigma) @>e_{\gamma}>> (L_{\gamma},\sigma_{\gamma}) \\
@Ve_{\eta}VV 				@VVe_{\eta}V\\
(L_{\eta}, \sigma_{\eta}) @>e_{\gamma}>> (L_{\gamma,\eta}, \sigma_{\gamma,\eta})\\
\end{CD}
\end{equation}
If $\gamma$ or $\eta$ is in $\rightBridges{L}$ then $e_{\gamma}e_{\eta} = e_{\eta}e_{\gamma}$, whereas if both $\gamma, \eta \in \leftBridges{L}$ then $e_{\gamma}e_{\eta} = -e_{\eta}e_{\gamma}$. 

\item One edge, $e_{\alpha}$, is a bridge edge for $\gamma$ while $e_{\beta}$ is a decoration edge for $C \in \circles{L}$. In this case $\sigma(C) = +$. By the disjoint support assumption $\sigma_{\gamma}(C) = +$ as well, so there is a decoration edge $e_{C}: (L_{\gamma}, \sigma_{\gamma}) \longrightarrow (L_{\gamma}, \sigma_{\gamma,C})$ {\em with the same location}  as $e_{\beta}$. Call this edge $e_{\beta'}$. On the other hand $\sigma_{C}$ will assign the same decorations to the circle in $\gamma$'s support, thus there is an edge $e_{\alpha'}=e_{(\gamma,\sigma_{C}, \sigma_{C,\gamma})}: (L,\sigma_{C}) \longrightarrow (L, \sigma_{C,\gamma})$. The edge is uniquely identified by requiring $\sigma_{C,\gamma} = \sigma_{\gamma,C}$ (since $\sigma_{\gamma}$ is given in the original data). Consequently, we obtain a square
\begin{equation}
\begin{CD}
(L,\sigma) @>e_{\gamma}>> (L_{\gamma},\sigma_{\gamma}) \\
@Ve_{C}VV 				@VVe_{C}V\\
(L, \sigma_{C}) @>e_{\gamma}>> (L_{\gamma}, \sigma_{\gamma,C})\\
\end{CD}
\end{equation}
If $\gamma$  is in $\rightBridges{L}$ or $e_{C} = \righty{e_{C}}$ then $e_{\gamma}e_{C} = e_{C}e_{\gamma}$, whereas if both $\gamma \in \leftBridges{L}$ and $e_{C} = \lefty{e_{C}}$ then $e_{\gamma}\lefty{e_{C}} = - \lefty{e_{C}}e_{\gamma}$.

\item Both edges are decoration edges for circles $C_{\alpha}$ and $C_{\beta}$. To have disjoint supports we need $C_{\alpha} \neq C_{\beta}$. Furthermore, $\sigma(C_{\alpha}) = \sigma(C_{\beta}) = +$ for the edges to exist. Since $\sigma_{\alpha} = \sigma_{C_{\alpha}}$ will assign a $+$ to $C_{\beta}$, there is a decoration edge $(L,\sigma_{\alpha}) \longrightarrow (L,\sigma_{\alpha, C_{\beta}})$ with the same location as $e_{\beta}$. We denote this edge by $e_{\beta'}$. Likewise there is a decoration edge $e_{\alpha'} : (L,\sigma_{\beta}) \longrightarrow (L,\sigma_{\beta, C_{\alpha}})$ with the same location as $e_{\alpha}$. $\sigma_{\beta, C_{\alpha}} = \sigma_{\alpha, C_{\beta}}$ since both are $\sigma$ away from $C_{\alpha}$ and $C_{\beta}$, but assign these circles $-$. We call this decoration $\sigma_{\alpha,\beta}$. Thus there is a square in $\bridgeGraph{n}$ of the following form:
\begin{equation}
\begin{CD}
(L,\sigma) @>e_{C_{\alpha}}>> (L,\sigma_{C_{\alpha}}) \\
@Ve_{C_{\beta}}VV 				@VVe_{C_{\beta}}V\\
(L, \sigma_{C_{\beta}}) @>e_{C_{\alpha}}>> (L, \sigma_{C_{\alpha},C_{\beta}})\\
\end{CD}
\end{equation}

This situation gives rise to four relations 
$$
\begin{array}{ccc}
\righty{e_{C_{\alpha}}}\righty{e_{C_{\beta}}} =  \righty{e_{C_{\beta}}}\righty{e_{C_{\alpha}}} & \hspace{1in} & \righty{e_{C_{\alpha}}}\lefty{e_{C_{\beta}}} =  \lefty{e_{C_{\beta}}}\righty{e_{C_{\alpha}}} \\
\lefty{e_{C_{\alpha}}}\righty{e_{C_{\beta}}} =  \righty{e_{C_{\beta}}}\lefty{e_{C_{\alpha}}} & \hspace{1in} & \lefty{e_{C_{\alpha}}}\lefty{e_{C_{\beta}}} =  -\lefty{e_{C_{\beta}}}\lefty{e_{C_{\alpha}}} \\
\end{array}
$$
\end{enumerate}
\noindent We note that the type (decoration vs. bridge) and the location are the same for $e_{\alpha}$ and $e_{\alpha'}$ as well as for the pair $e_{\beta}$ and $e_{\beta'}$.\\
\ \\

\noindent{\bf Relations for decoration edges:} When the support of $e_{C}$ is not disjoint from that of $e_{\gamma} : (L, \sigma) \rightarrow (L_{\gamma}, \sigma_{\gamma})$ we must distinguish $e_{C}$ based on its location. 

\begin{enumerate}
\item {\bf The relations for $\righty{e_{C}}$:} Suppose first that $\gamma \in \merge{L}$ merges $C_{1}$ and $C_{2}$ to get $C \in \circles{L_{\gamma}}$, and $\sigma(C_{1}) = \sigma(C_{2}) = +$. The only decoration edges without disjoint support are those from $C_{1}$, $C_{2}$ and $C$ (which has $\sigma_{\gamma}(C) = +$ as well). It is straightforward to verify that the following squares
exist in $\bridgeGraph{n}$:
\begin{equation}
\begin{CD}
(L,\sigma) @>e_{\gamma}>> (L_{\gamma},\sigma_{\gamma}) \\
@V\righty{e_{C_{1}}}VV 				@VV\righty{e_{C}}V\\
(L, \sigma_{C_{1}}) @>e_{\gamma}>> (L_{\gamma}, \sigma_{C})\\
\end{CD}
\hspace{0.75in}
\begin{CD}
(L,\sigma) @>e_{\gamma}>> (L_{\gamma},\sigma_{\gamma}) \\
@V\righty{e_{C_{2}}}VV 				@VV\righty{e_{C}}V\\
(L, \sigma_{C_{2}}) @>e_{\gamma}>> (L_{\gamma}, \sigma_{C})\\
\end{CD}
\end{equation}
Both of these squares should provide relations, stated specifically as
\begin{equation}
\righty{e_{C_{1}}}m_{(\gamma,\sigma_{C_{1}},\sigma_{C})} = \righty{e_{C_{2}}}m_{(\gamma,\sigma_{C_{2}},\sigma_{C})} = m_{(\gamma,\sigma,\sigma_{\gamma})} \righty{e_{C}}
\end{equation}
Note that if $\sigma(C_{i}) = -$ for either $i=1$ or $2$, then no such squares exist in $\bridgeGraph{n}$. \\
\ \\
Dually, if surgery on $\gamma \in \fission{L}$ divides circle $C$ into $C_{1}$ and $C_{2}$ in $\circles{L_{\gamma}}$, and $\sigma$ assigns $+$ to $C$, then for paths starting at $(L,\sigma)$ 
\begin{equation}
\righty{e_{C}}f_{(\gamma,\sigma_{C},\sigma_{C,\gamma})} =  f_{(\gamma,\sigma,\sigma^{1}_{\gamma})}\righty{e_{C_{1}}} = f_{(\gamma,\sigma,\sigma^{2}_{\gamma})} \righty{e_{C_{2}}}
\end{equation}
where $\sigma^{i}_{\gamma}$ assigns $+$ to $C_{i}$ and $-$ to $C_{3-i}$. Thus, after following the edge $f_{(\gamma,\sigma,\sigma^{1}_{\gamma})}$ only $\righty{e_{C_{1}}}$ is available for the next step. There are no such relations if $\sigma(C) = -$.

\item {\bf The relations for $\lefty{e_{C}}$:} The squares found for $\righty{e_{C}}$ still exist for $\lefty{e_{C}}$. However, {\em we do not use them to define relations}. Instead, we note that the two paths around the edges of the squares have the same source and target, and thus contribute to the same vector space of morphisms in $\mathcal{Q}\bridgeGraph{n}$. Consequently, we can add them to obtain a new morphism in the same set of morphisms. In particular, we impose the following relation when $\sigma(C_{1}) = \sigma(C_{2}) = +$ and $\gamma \in \rightMerge{L}$:
\begin{equation}\label{rel:lefty1}
\lefty{e_{C_{1}}}m_{(\gamma,\sigma_{C_{1}},\sigma_{C})} + \lefty{e_{C_{2}}}m_{(\gamma,\sigma_{C_{2}},\sigma_{C})} - m_{(\gamma,\sigma,\sigma_{\gamma})} \lefty{e_{C}} = 0
\end{equation}
and when $\sigma(C) = +$ and $\gamma \in \rightFission{L}$ divides $C$ into $C_{1}$ and $C_{2}$
\begin{equation}\label{rel:lefty2}
\lefty{e_{C}}f_{(\gamma,\sigma_{C},\sigma_{C,\gamma})} +  f_{(\gamma,\sigma,\sigma^{1}_{\gamma})}\lefty{e_{C_{1}}} - f_{(\gamma,\sigma,\sigma^{2}_{\gamma})} \lefty{e_{C_{2}}} = 0
\end{equation}
whereas for $\gamma \in \leftMerge{L}$
\begin{equation}\label{rel:lefty3}
\lefty{e_{C_{1}}}m_{(\gamma,\sigma_{C_{1}},\sigma_{C})} + \lefty{e_{C_{2}}}m_{(\gamma,\sigma_{C_{2}},\sigma_{C})} + m_{(\gamma,\sigma,\sigma_{\gamma})} \lefty{e_{C}} = 0
\end{equation}
and when $\sigma(C) = +$ and $\gamma \in \leftFission{L}$ divides $C$ into $C_{1}$ and $C_{2}$
\begin{equation}\label{rel:lefty4}
\lefty{e_{C}}f_{(\gamma,\sigma_{C},\sigma_{C,\gamma})} +  f_{(\gamma,\sigma,\sigma^{1}_{\gamma})}\lefty{e_{C_{1}}} + f_{(\gamma,\sigma,\sigma^{2}_{\gamma})} \lefty{e_{C_{2}}} = 0
\end{equation}
\end{enumerate}

\noindent{\bf Bridge relations:} Within $\bridgeGraph{n}$ there are also squares of the following form:
\begin{equation} \label{cd:commute}
\begin{CD}
(L,\sigma) @>e_{(\gamma,\sigma_{00},\sigma_{01})}>> (L_{\gamma},\sigma_{01}) \\
@Ve_{(\gamma',\sigma_{00},\sigma_{10})}VV 				@VVe_{(\gamma',\sigma_{01},\sigma_{11})}V\\
(L_{\gamma'}, \sigma_{10}) @>e_{(\gamma,\sigma_{01},\sigma_{11})}>> (L_{\gamma,\gamma'}, \sigma_{11})\\
\end{CD}
\end{equation}
where
\begin{enumerate}
\item $\gamma \in \bridges{L}$, $\gamma' \in B_{|\,|}(L,\gamma)$,
\item $\sigma_{00}$, $\sigma_{01}$, $\sigma_{10}$, $\sigma_{11}$ are any decorations for which
such a square exists.
\end{enumerate}
The crucial element distinguishing these squares is the pattern of $\gamma$ for both horizontal maps, and $\gamma'$ for both vertical maps. These should be taken as the actual arcs representing the bridges in order for the identification to proceed.\\
\ \\
\noindent {\bf Note:} Since $\gamma'$ is parallel and not equal to $\gamma$, there are classes $\gamma' \in \bridges{L_{\gamma}}$ and $\gamma \in \bridges{L_{\gamma'}}$. In addition, $L_{\gamma,\gamma'} = L_{\gamma',\gamma}$. Thus the only difficulty in filling in a square given the top and right (or left and bottom) sides will be in finding an appropriate decoration for the fourth corner.  Such squares exist frequently:

\begin{prop}
If $\gamma \in \bridges{L}$, $\gamma' \in B_{|\,|}(L,\gamma)$, and there are edges
$$
\begin{CD}
(L,s) @>e_{(\gamma,\sigma,\sigma')}>> (L_{\gamma},\sigma') @>e_{(\gamma',\sigma',\sigma'')}>> (L_{\gamma,\gamma'}, \sigma'')\\
\end{CD}
$$
then there is at least one decoration $\overline{\sigma}$ on $L_{\gamma'}$ for which we can find edges
$$
\begin{CD}
(L,\sigma) @>e_{(\gamma',\sigma,\overline{\sigma})}>> (L_{\gamma'},\overline{\sigma}) @>{e_{(\gamma,\overline{\sigma},\sigma'')}}>> (L_{\gamma,\gamma'}, \sigma'')\\
\end{CD}
$$
\end{prop}

\noindent For some of the squares in equation \ref{cd:commute} we impose the following relation:

\begin{equation}
e_{(\gamma,\sigma_{00},\sigma_{01})}e_{(\gamma',\sigma_{01},\sigma_{11})} = (-1)^{\leftnorm{e_{\gamma}}\leftnorm{e_{\gamma}'}} e_{(\gamma',\sigma_{00},\sigma_{10})}e_{(\gamma,\sigma_{10},\sigma_{11})}
\end{equation}

\noindent We do so when
\begin{enumerate}
\item either of the arcs $\gamma$ and $\gamma'$ is in $\rightHalf$, or
\item the arcs $\gamma$ and $\gamma'$ represent {\em left} bridges in $(L,\sigma)$, with $\lefty{\gamma} \in B_{o}(L,\lefty{\gamma}')$
\end{enumerate}

\noindent Such squares also exist in the case where $\lefty{\gamma} \in B_{s}(L,\lefty{\gamma}')$.  However, in this case, they occur in triples: for pairs chosen from $\lefty{\gamma}$, $\lefty{\gamma}'$ and $\lefty{\gamma}''$, the result of sliding the end of $\lefty{\gamma}$ over $\lefty{\gamma}'$. We impose a different relation in this setting. \\
\ \\
\noindent In $L_{\lefty{\gamma}}$, $\lefty{\gamma}'$ and $\lefty{\gamma}''$ represent the same bridge $\lefty{\delta}$. Likewise, in $L_{\lefty{\gamma}'}$, both $\lefty{\gamma}$ and $\lefty{\gamma}''$ represent the same bridge $\lefty{\zeta}$. Finally, let $\lefty{\eta}$ be the image of $\lefty{\gamma}$ and $\lefty{\gamma}'$ in $L_{\lefty{\gamma}''}$. Then
$$
\lefty{e_{\gamma}}\lefty{e_{\delta}} + \lefty{e_{\gamma}}'\lefty{e_{\zeta}} + \lefty{e_{\gamma}}''\lefty{e_{\eta}} = 0
$$
whenever there are comptaible decorations on $L_{\lefty{\gamma}}$, $L_{\lefty{\gamma}'}$, $L_{\lefty{\gamma}''}$. Note that if one of the terms does not exist, we obtain ``anti-commutativity'' for the other two terms.\\

\noindent {\bf Other bridge relations:} Suppose $\gamma \in \leftBridges{L}$ and $\eta \in B_{\pitchfork}(L_{\gamma}, \gamma^{\dagger})$. Then $\eta \in \leftBridges{L}$ as well. For a path
$$
\begin{CD}
(L,s) @>e_{(\gamma,\sigma,\sigma')}>> (L_{\gamma},\sigma') @>e_{(\eta,\sigma',\sigma'')}>> (L_{\gamma,\eta}, \sigma'')\\
\end{CD}
$$
we require that
\begin{equation}
e_{(\gamma,\sigma,\sigma')}e_{(\eta,\sigma',\sigma'')} = 0
\end{equation}
\ \\
\noindent{\bf There and back again..} There is one final set of relations. When composing the morphisms for two bridges $\gamma$ and $\gamma'$, we need $\gamma'$ to be a bridge for $L_{\gamma}$. There is a new bridge $\gamma^{\dagger}$ in $L_{\gamma}$ not coming from one for $L$. Thus we can always try to form path(s) by surgering first along $\gamma$ and then along $\gamma^{\dagger}$. The result of these surgeries is $L$ again with a different decoration. 

\begin{prop}
If there is a path $\rho$
$$
\begin{CD}
(L,\sigma) @>e_{(\gamma,\sigma,\sigma')}>> (L_{\gamma},\sigma') @>e_{(\gamma^{\dagger},\sigma',\sigma'')}>> (L, \sigma'')
\end{CD}
$$
then there is a circle $C \in \circles{L}$, in the support of $\gamma$, where $\sigma(C)=+$ and $\sigma''=\sigma_{C}$. 
\end{prop}

\noindent{\bf Proof:} Each edge is a bridge edge and thus lowers $\iota$ by $1$. The composition lowers $\iota$ by $2$ but returns to the same cleaved link, equipped with a different decoration. By a case-by-case analysis we see that $\sigma''$ cannot change a $-$ to a $+$ sign, so it must change a single circle from $+$ to $-$. Since the edges only affect decorations on circles in the support of $\gamma$ and $\gamma^{\dagger}$ the result follows. $\Diamond$\\
\ \\
\begin{defn}
For each path $\rho$ as above, the circle $C$ will be called the {\em active circle} for $\rho$. 
\end{defn}

\noindent For each path
$$
\begin{CD}
(L,\sigma) @>e_{(\gamma,\sigma,\sigma')}>> (L_{\gamma},\sigma') @>e_{(\gamma^{\dagger},\sigma',\sigma_{C})}>> (L, \sigma_{C})
\end{CD}
$$
with active circle $C$ and $\gamma \in \rightBridges{L}$ we impose the relation
\begin{equation}\label{rel:special}
\righty{e}_{(\gamma,\sigma,\sigma')}\righty{e}_{(\gamma^{\dagger},\sigma',\sigma_{C})} = \righty{e_{C}}
\end{equation}
where $\sigma'$ is any choice of decoration compatible with surgery on $\gamma$.\\
\ \\
\noindent Quotienting by all these relations gives a pre-additive category $\mathcal{B}\bridgeGraph{n}$. 

\begin{prop}
The relations are homogeneous for the bigrading on $\mathcal{Q}\bridgeGraph{n}$, and thus induce a bigrading on the morphism sets in $\mathcal{B}\bridgeGraph{n}$.
\end{prop}
\noindent {\bf Proof:} The bigrading of each generator depends only on its location (left or right) and its type (bridge or decoration), and not on the vertices the corresponding edge joins. In each of the relations above, each term has the same set of locations and types appearing, except for the relation in equation \ref{rel:special}. However, in that case the left side has bigrading $(0,-1/2) + (0,-1/2) = (0,-1)$ which is the same as the bigrading for the element on the right side, $\righty{e_{C}}$. $\Diamond$\\ 
\ \\
\begin{defn}
$\mathcal{I}_{n}$ is the sub-algebra generated by the idempotents $I_{(L,\sigma)}$ corresponding to length $0$ paths in $\bridgeGraph{n}$. 
\end{defn}

\subsection{A differential on $\mathcal{B}\bridgeGraph{n}$}

\noindent Surgery along a bridge $\gamma$ followed by that on $\gamma^{\dagger}$ gives a relation when  $\gamma$ is in $\rightHalf$ (as in equation \ref{rel:special}); however, the same path exists for bridges in $\leftHalf$. This {\em does not} correspond to a relation. Rather, it gives a differential on $\mathcal{B}\bridgeGraph{n}$. The sign conventions for this differential are somewhat unusual:

\begin{defn}
Let $\mathcal{C}$ be a pre-additive category over a ring $R$. A (right) differential on $\mathcal{C}$ is a collection of $R$-linear maps
$$
d_{v_{1},v_{2}} : \mathrm{Mor}(v_{1},v_{2}) \rightarrow \mathrm{Mor}(v_{1},v_{2})
$$
such that
\begin{enumerate}
\item for each identity morphism $I_{v} \in \mathrm{Mor}(v,v)$, $d_{v,v}(I_{v}) = 0$
\item for $\alpha \in \mathrm{Mor}(u,v)$ and $\beta \in \mathrm{Mor}(v,w)$ the Leibniz identity holds:
\begin{equation}\label{eqn:Leibniz}
d_{u,w}(\alpha\beta) = (-1)^{|\beta|}(d_{u,v}(\alpha)\big)\beta + \alpha\big(d_{v,w}(\beta)\big)
\end{equation}
\end{enumerate}
\end{defn}

\begin{prop}
For each $(L,\sigma)$ and circle $C \in \circles{L}$ with $\sigma(C) = +$. Let $d_{(L,\sigma), (L,\sigma_{C})}$ be the $k$-linear map on $\mathrm{Mor}\big((L,\sigma),(L,\sigma_{C})\big)$ defined by homorphically extending to all paths the
specification 
\begin{equation}\label{eq:dga}
d_{(L,\sigma), (L,\sigma_{C})}(\lefty{e_{C}}) = -\sum e_{(\gamma,\sigma,\sigma')}e_{(\gamma^{\dagger},\sigma',\sigma_{C})}
\end{equation}
where the sum is over all length two paths 
$$ 
\begin{CD} 
(L,\sigma) @>e_{(\gamma,\sigma,\sigma')}>> (L_{\gamma},\sigma') @>e_{(\gamma^{\dagger},\sigma',\sigma_{C})}>> (L,\sigma_{C}) 
\end{CD}
$$ 
with $\gamma \in \leftBridges{L}$, and $d(e) = 0$ for every other edge $e$ in $\bridgeGraph{n}$, and for the identity morphisms $I_{(L,\sigma)}$. Then $d$ can be extended using the (right) Leibniz rule in equation \ref{eqn:Leibniz} to a $(1,0)$ differential on the bigraded cateogy $\mathcal{B}\bridgeGraph{n}$ with $|\alpha| = \leftnorm{\alpha}$.
\end{prop}

\noindent We will call the resulting differential on $\mathcal{B} \bridgeGraph{n}$, $d_{\Gamma_{n}}$. Note that the Leibniz rule extends $d$ to $\mathcal{Q}\bridgeGraph{n}$ without difficulty. The import of the proposition is that $d$ is compatible with the relations defining $\mathcal{B}\bridgeGraph{n}$. \\
\ \\

\noindent{\bf Proof:} The maps $d_{(L,s),(L',s')}$ specified in the proposition define a collection of $k$-linear maps on the morphism sets corresponding to length $0$ and length $1$ paths. We first verify that this is a $(1,0)$ differential on the length $0$ and $1$ paths. $d$ is trivial on length $0$ paths, and thus is a $(1,0)$ differential on them. On length $1$ paths, $d$ is only non-trivial on $\lefty{e_{C}}$, which has bigrading $(1,1)$. It image consists of a linear combination of terms of the form $\lefty{e}_{(\gamma,\sigma,\sigma')}\lefty{e}_{(\gamma^{\dagger},\sigma',\sigma_{C})}$ in bigrading $(1,1/2) + (1,1/2) = (2,1)$. Thus, the image is in bigrading $(1,0)$ greater than that of $\lefty{e_{C}}$, as required.\\
\ \\
\noindent The Leibniz rule uniquely extends $d$ to a $(1,0)$ differential on all the other generating paths in $\mathcal{Q}\bridgeGraph{n}$ since it will only expand one $\lefty{e_{C}}$ and otherwise be trivial. If this extension is compatible with the relations defining the quotient $\mathcal{B}\bridgeGraph{n}$ then there is an induced differential on $\mathcal{B}\bridgeGraph{n}$. The reaminder of the proof consists of verifying this compatibility. \\
\ \\
\noindent We can start by by noting that the $d$ maps will be trivial on any path which does not include at least one $\lefty{e_{C}}$ edge for some circle $C$. Thus, any relation not involving the paths $\lefty{e_{C}}$ will immediately be compatible with the Leibniz identity. The remaining relations consist of two types: commuting relations for disjoint supports, and the relations \ref{rel:lefty1} and \ref{rel:lefty2}. \\
\ \\

\noindent {\bf Verifying the relations for disjoint support:} Suppose $(L,\sigma)$ has a circle $C$ with $\sigma(C) = +$ and there is a square in $\bridgeGraph{n}$
$$ 
\begin{CD} 
(L,\sigma) @>e_{(\gamma,\sigma,\sigma_{\gamma})}>> (L_{\gamma},\sigma_{\gamma}) \\
@V\righty{e_{C}}VV           @VV\righty{e_{C}}V\\
(L,\sigma_{C}) @>e_{(\gamma,\sigma_{C},\sigma_{C,\gamma})}>> (L_{\gamma},\sigma_{C,\gamma})
\end{CD}
$$
for a bridge $\gamma$ where $C$ is not in the support of $\gamma$ and $\sigma_{\gamma}$. Suppose that $\gamma \in \rightBridges{L}$. Then we need to show that $d(e_{(\gamma,\sigma,\sigma_{\gamma})}\lefty{e_{C}}) = d(\lefty{e_{C}}e_{(\gamma,\sigma_{C},\sigma_{C,\gamma})})$. Using the Leibniz rule we need only  verify that $e_{(\gamma,\sigma,\sigma_{\gamma})}d(\lefty{e_{C}}) = (-1)^{\leftnorm{e_{\gamma}}}d(\lefty{e_{C}})e_{(\gamma,\sigma_{C},\sigma_{C,\gamma})}$. But $\leftnorm{e_{\gamma}} = 0$ so we need only verify that $e_{(\gamma,\sigma,\sigma_{\gamma})}d(\lefty{e_{C}}) = d(\lefty{e_{C}})e_{(\gamma,\sigma_{C},\sigma_{C,\gamma})}$. Using equation \ref{eq:dga} it suffices to show for any bridge $\eta \in \leftBridges{L}$, which includes $C$ in its support, that there is a $1-1$ correspondence $\sigma' \leftrightarrow \sigma''$ such that
$$
e_{(\gamma,\sigma,\sigma_{\gamma})}e_{(\eta,\sigma_{\gamma},\sigma')}e_{(\eta^{\dagger},\sigma',\sigma_{\gamma,C})} = e_{(\eta,\sigma,\sigma'')}e_{(\eta^{\dagger},\sigma'',\overline{\sigma})}e_{(\gamma,\overline{\sigma},\sigma_{\gamma,C})} 
$$
under the requirement that $\sigma$ and $\sigma'$ ensure that $C$ is the active circle for the $\eta\eta^{\dagger}$ paths (and all required paths exist). \\
\ \\
\noindent If the support of $\eta$ does not intersect the support of $\gamma$ then the support of $\eta$ is $C$ and $C'$ and the support of $\gamma$ is $D$ and $D'$ where $C=C'$ is possible as is $D=D'$. If both are merges, then $C'$ can also equal $D'$. 
For all these cases interchanging the $\gamma$ path and the $\eta$ path can be done uniquely, when the two-step path exists. More specifically, examining the bridge relations, we see that there is a unique $\sigma''$ which makes $e_{(\gamma,\sigma,\sigma_{\gamma})}e_{(\eta,\sigma_{\gamma},\sigma')} = e_{(\eta,\sigma,\sigma'')}e_{(\gamma,\sigma'',\sigma')}$
and a unique $\overline{\sigma}$ which makes $e_{(\gamma,\sigma'',\sigma')}e_{(\eta^{\dagger},\sigma',\sigma_{\gamma,C})} = 
e_{(\eta^{\dagger},\sigma'',\overline{\sigma})}e_{(\gamma,\overline{\sigma},\sigma_{\gamma,C})}$. Putting these together gives the desired identity and the $1-1$ correspondence.\\
\ \\
\noindent When $\gamma \in \leftBridges{L}$ we show that $d(e_{(\gamma,\sigma,\sigma_{\gamma})}\lefty{e_{C}}) = d(-\lefty{e_{C}}e_{(\gamma,\sigma_{C},\sigma_{C,\gamma})})$. Again we only need to verify that $e_{(\gamma,\sigma,\sigma_{\gamma})}d(\lefty{e_{C}}) = (-1)^{\leftnorm{e_{\gamma}}}d(-\lefty{e_{C}})e_{(\gamma,\sigma_{C},\sigma_{C,\gamma})}$. Since $\leftnorm{e_{\gamma}} = 1$ this reduces to  $e_{(\gamma,\sigma,\sigma_{\gamma})}d(\lefty{e_{C}}) = d(\lefty{e_{C}})e_{(\gamma,\sigma_{C},\sigma_{C,\gamma})}$. The argument is the same as before, except that $e_{(\gamma,\sigma,\sigma_{\gamma})}e_{(\eta,\sigma_{\gamma},\sigma')} = - e_{(\eta,\sigma,\sigma'')}e_{(\gamma,\sigma'',\sigma')}$ and $e_{(\gamma,\sigma'',\sigma')}e_{(\eta^{\dagger},\sigma',\sigma_{\gamma,C})} = -
e_{(\eta^{\dagger},\sigma'',\overline{\sigma})}e_{(\gamma,\overline{\sigma},\sigma_{\gamma,C})}$. These signs cancel in the three term product, so the argument will still apply.\\
\ \\
\noindent The second type of disjoint support occurs if we look at $d(\lefty{e_{C}}\righty{e_{C'}})$ or $d(\lefty{e_{C}}\lefty{e_{C'}})$ with $C \neq C'$. Such paths only occur for $(L,\sigma)$ with $\sigma(C) = \sigma(C') = +$. In the first case we need to verify that $d(\lefty{e_{C}}\righty{e_{C'}}) = d(\righty{e_{C'}}\lefty{e_{C}})$, or equivalently $(-1)^{0}d(\lefty{e_{C}})\righty{e_{C'}} = \righty{e_{C'}}d(\lefty{e_{C}})$. Again, we can reduce this to verifying equalities 
$$
\righty{e_{C'}}e_{(\eta,\sigma_{C'},\sigma')}e_{(\eta^{\dagger},\sigma',\sigma_{C,C'})} = e_{(\eta,\sigma,\sigma'')}e_{(\eta^{\dagger},\sigma'',s_{C})}\righty{e_{C'}} 
$$ 
where $\eta\eta^{\dagger}$ has $C$ as its active circle. If the support of $\eta$ does not include $C'$, then the elements commute since they have disjoint supports, and $\righty{e_{C'}}$ is even. Furthermore, the only difference between $\sigma'$ and $\sigma''$ will be the change in decoration on $C'$. If $C'$ is in the support of $\eta$ then we can mark the edges by whether they are merges or divisions and see that
$$
\righty{e_{C'}}m_{(\eta,\sigma_{C'},\sigma')}f_{(\eta^{\dagger},\sigma'',\sigma_{C,C'})} = m_{(\eta,\sigma,\sigma'')}\righty{e_{C\#C'}}f_{(\eta^{\dagger},\sigma''_{C\#C'},\sigma_{C,C'})}= m_{(\eta,\sigma,\sigma'')}f_{(\eta^{\dagger},\sigma'',\sigma_{C})}\righty{e_{C'}} 
$$ 
where $\sigma''$ is uniquely determined by the merge of the $+$ decorated circles, and the last step is forced by the
requirement that $C$ be the active circle for $\eta$.\\
\ \\
\noindent There is more to verifying $d(\lefty{e_{C}}\lefty{e_{C'}}) = - d(\lefty{e_{C'}}\lefty{e_{C}})$. Our Leibniz rule requires $d(\lefty{e_{C}}\lefty{e_{C'}}) = -d(\lefty{e_{C}})\lefty{e_{C'}} + \lefty{e_{C}}d(\lefty{e_{C'}})$ and
$-d(\lefty{e_{C'}}\lefty{e_{C}}) = d(\lefty{e_{C'}})\lefty{e_{C}} - \lefty{e_{C'}}d(\lefty{e_{C}})$. 
The terms in $d(\lefty{e_{C}})\lefty{e_{C'}}$ consist of $e_{\gamma}e_{\gamma^{\dagger}}\lefty{e_{C'}}$, summed over bridges $\gamma$. $e_{\gamma}e_{\gamma^{\dagger}}$ is an even element which commutes with $\lefty{e_{C'}}$ unless $\gamma$ has one foot on $C$ and one foot on $C'$. Excepting that case, $e_{\gamma}e_{\gamma^{\dagger}}\lefty{e_{C'}} = \lefty{e_{C'}}e_{\gamma}e_{\gamma^{\dagger}}$, a term contributing to $\lefty{e_{C'}}d(\lefty{e_{C}})$. Furthermore, both occur with a minus sign in the Leibniz identity. A similar argument shows an identification between the terms in $\lefty{e_{C}}d(\lefty{e_{C'}})$ and $d(\lefty{e_{C'}})\lefty{e_{C}}$ for those bridges with only one foot on $C$ and $C'$. \\
\ \\
\noindent Returning to the exceptional case: let $\gamma \in \leftBridges{L}$ join $C$ and $C'$. This contibutes terms to both $-d(\lefty{e_{C}})\lefty{e_{C'}}$ and $\lefty{e_{C}}d(\lefty{e_{C'}})$ (since $\sigma(C)=\sigma(C') = +$)whose sum equals
$$
-m_{(\gamma,\sigma,\sigma_{\gamma})}f_{(\gamma^{\dagger},\sigma_{\gamma},\sigma_{C})}\lefty{e_{C'}} + \lefty{e_{C}}m_{(\gamma,\sigma_{C},\sigma_{C,\gamma})}f_{(\gamma^{\dagger},\sigma_{C,\gamma},\sigma_{C,C'})}
$$
We use equation \ref{rel:lefty3} to transform
$$
\lefty{e_{C}}m_{(\gamma,\sigma_{C},\sigma_{C,\gamma})} = -\lefty{e_{C'}}m_{(\gamma,\sigma_{C'}, \sigma_{C',\gamma}} - m_{(\gamma,\sigma,\sigma_{\gamma})}\lefty{e_{C\#C'}}
$$
Multiplying on the right by $f_{(\gamma^{\dagger},s'',s_{C,C'})}$ gives
$$
\lefty{e_{C}}m_{(\gamma,\sigma_{C},\sigma_{C,\gamma})}f_{(\gamma^{\dagger},\sigma_{C,\gamma},\sigma_{C,C'})} = - \lefty{e_{C'}}m_{(\gamma,\sigma_{C'}, \sigma_{C',\gamma}}f_{(\gamma^{\dagger},\sigma_{C,\gamma},\sigma_{C,C'})} - m_{(\gamma,\sigma,\sigma_{\gamma})}\lefty{e_{C\#C'}}f_{(\gamma^{\dagger},\sigma_{C,\gamma},\sigma_{C,C'})}
$$
On the other hand, using equation \ref{rel:lefty4} we obtain
$$
-m_{(\gamma,\sigma,\sigma_{\gamma})}f_{(\gamma^{\dagger},\sigma_{\gamma}, \sigma_{C,\gamma})}\lefty{e_{C'}} = + m_{(\gamma,\sigma,\sigma_{\gamma})}f_{(\gamma^{\dagger},\sigma_{\gamma},\sigma_{C',\gamma})}\lefty{e_{C}} + 
m_{(\gamma,\sigma,\sigma_{\gamma})}\lefty{e_{C\#C'}}f_{(\gamma^{\dagger},\sigma_{C\#C',\gamma},\sigma_{C,C'})}
$$
Adding the results yields
$$
-m_{(\gamma,\sigma,\sigma_{\gamma})}f_{(\gamma^{\dagger},\sigma_{\gamma}, \sigma_{C,\gamma})}\lefty{e_{C'}} + \lefty{e_{C}}m_{(\gamma,\sigma_{C},\sigma_{C,\gamma})}f_{(\gamma^{\dagger},\sigma_{C,\gamma},\sigma_{C,C'})} =
$$
$$
\hspace{1in}  -\lefty{e_{C'}}m_{(\gamma,\sigma_{C'}, \sigma_{C',\gamma}}f_{(\gamma^{\dagger},\sigma_{C,\gamma},\sigma_{C,C'})} + m_{(\gamma,\sigma,\sigma_{\gamma})}f_{(\gamma^{\dagger},\sigma_{\gamma},\sigma_{C',\gamma})}\lefty{e_{C}}
$$
The right hand side is precisely the terms from $\gamma$ in $+d(\lefty{e_{C'}})\lefty{e_{C}} - \lefty{e_{C'}}d(\lefty{e_{C}})$, thereby verifying the Leibniz identity. \\
\ \\
\noindent{\bf Verifying the relations from bridges:} Finally, we need to see that the definition of $d$ is compatible with the
relations
\begin{equation}
\lefty{e_{C_{1}}}m_{(\gamma,\sigma_{C_{1}},\sigma_{\gamma,C})} + \lefty{e_{C_{1}}}m_{(\gamma,\sigma_{C_{1}},\sigma_{\gamma,C})} \pm m_{(\gamma,\sigma,\sigma_{\gamma})} \lefty{e_{C}} = 0
\end{equation}
and
\begin{equation}
\lefty{e_{C}}f_{(\gamma,\sigma_{C},\sigma_{C,\gamma})} +  f_{(\gamma,\sigma,\sigma^{1}_{\gamma})}\lefty{e_{C_{1}}} \pm f_{(\gamma,\sigma,\sigma^{2}_{\gamma})} \lefty{e_{C_{2}}} = 0
\end{equation}  
for any bridge $\gamma$. \\
\ \\
\noindent We will only verify the first as the second is similar. Note that all the terms vanish when $\sigma(C_1) = -$ or $\sigma(C_{2}) = -$ so we will assume that $\sigma(C_{2}) = \sigma(C_{1}) = +$. A $+$ occurs in the relation for $\gamma \in \leftBridges{L}$, while a minus occurs for $\gamma \in \rightBridges{L}$. If we apply $d$ to all the terms, and simplify using the Leibniz identity, we will need for $\gamma \in \leftBridges{L}$ that
$$
-d(\lefty{e_{C_{1}}})m_{(\gamma,\sigma_{C_{1}},\sigma_{\gamma,C})} - d(\lefty{e_{C_{2}}})m_{(\gamma,\sigma_{C_{2}},\sigma_{\gamma,C})} + m_{(\gamma,\sigma,\sigma_{\gamma})} d(\lefty{e_{C}}) = 0
$$
while for $\gamma$ in $\rightBridges{L}$ we will need 
$$
+d(\lefty{e_{C_{1}}})m_{(\gamma,\sigma_{C_{1}},\sigma_{\gamma,C})} + d(\lefty{e_{C_{2}}})m_{(\gamma,\sigma_{C_{2}},\sigma_{\gamma,C})} - m_{(\gamma,\sigma,\sigma_{\gamma})} d(\lefty{e_{C}}) = 0
$$
Both cases therefore require that we verify the relation
$$
m_{(\gamma,\sigma,\sigma_{\gamma})} d(\lefty{e_{C}}) = d(\lefty{e_{C_{1}}})m_{(\gamma,\sigma_{C_{1}},\sigma_{\gamma,C})} + d(\lefty{e_{C_{2}}})m_{(\gamma,\sigma_{C_{2}},\sigma_{\gamma,C})} 
$$
where $C= C_{1}\#_{\gamma} C_{2}$ is the result of merging $C_{1}$ and $C_{2}$ along $\gamma$. \\
\ \\
\noindent To compute $d(\lefty{e_{C}})$ we divide the bridges $\eta$ in $\leftBridges{L_{\gamma}}$ abutting $C$ into three groups: 1) those abutting $C_{1}$ and not $C_{2}$, 2) those abutting $C_{2}$ and not $C_{1}$, and 3) those between arcs in $C_{1}$ and $C_{2}$ (including $\gamma^{\dagger}$ when $\gamma \in \leftBridges{L}$). \\
\ \\
\noindent For a bridge in 1) commutativity of the bridge maps yields a unique choice of decoration $\sigma_{\eta}$ so that
$
m_{(\gamma,\sigma,\sigma_{\gamma})}e_{(\eta, \sigma_{\gamma}, \sigma_{\eta,\gamma})} = \pm e_{(\eta, \sigma, \sigma_{\eta})}m_{(\gamma,\sigma_{\eta},\sigma_{\eta,\gamma})}
$
and
$
m_{(\gamma,\sigma_{\eta},\sigma_{\eta,\gamma})}e_{(\eta^{\dagger}, \sigma_{\eta,\gamma}, \sigma_{C,\gamma})} = \pm
e_{(\eta^{\dagger}, \sigma_{\eta}, \sigma_{C_{1}})}m_{(\gamma,\sigma_{C_{1}},\sigma_{\gamma,C})}
$
where the sign depend on whether $\gamma \in \leftBridges{L}$ or not. In either case, together these imply
$$
m_{(\gamma,\sigma,\sigma_{\gamma})}e_{(\eta, \sigma_{\gamma}, \sigma_{\eta,\gamma})}e_{(\eta^{\dagger}, \sigma_{\eta,\gamma}, \sigma_{C,\gamma})} = 
e_{(\eta, \sigma, \sigma_{\eta})}e_{(\eta^{\dagger}, \sigma_{\eta}, \sigma_{C_{1}})}m_{(\gamma,\sigma_{C_{1}},\sigma_{\gamma,C})}
$$
where the left side is minus a term in $m_{\gamma}d(\lefty{e_{C}})$ and the right side is minus a term in $d(\lefty{e_{C_{1}}})m_{\gamma}$. 
Since the roles of $C_{1}$ and $C_{2}$ are symmetric, the same argument applies to bridges from 2) which abut $C_{2}$ and not $C_{1}$ to give terms in $m_{\gamma}d(\lefty{e_{C}})$ and $d(\lefty{e_{C_{2}}})m_{\gamma}$\\
\ \\
\noindent Now suppose $\eta$ abuts arcs from $C_{1}$ and $C_{2}$, but does not equal $\gamma^{\dagger}$. Then surgery on $C$ along $\eta$ must be a division. For each division there will be two choices of decoration, so we study the terms
\begin{equation}\label{eqn:special}
m_{(\gamma,\sigma,\sigma_{\gamma})}f_{(\eta, \sigma_{\gamma}, \sigma^{a}_{\eta,\gamma})}m_{(\eta^{\dagger}, \sigma^{a}_{\eta,\gamma}, \sigma_{C,\gamma})}
+ m_{(\gamma,\sigma,\sigma_{\gamma})}f_{(\eta, \sigma_{\gamma}, \sigma^{b}_{\eta,\gamma})}m_{(\eta^{\dagger}, \sigma^{b}_{\eta,\gamma}, \sigma_{C,\gamma})}
\end{equation}
where $C_{a}$ and $C_{b}$ are the result of dividing $C$ along $\eta$. These terms are the $\eta$ contribution to $m_{(\gamma,\sigma,\sigma_{\gamma})} d(\lefty{e_{C}})$. Now
$$
m_{(\gamma,\sigma,\sigma_{\gamma})}f_{(\eta, \sigma_{\gamma}, \sigma^{a}_{\eta,\gamma})}m_{(\eta^{\dagger}, \sigma^{a}_{\eta,\gamma}, \sigma_{C,\gamma})}
= \pm m_{(\eta,\sigma,\sigma_{\eta})}f_{(\gamma, \sigma_{\eta}, \sigma^{a}_{\gamma,\eta})}m_{(\eta^{\dagger}, \sigma^{a}_{\eta,\gamma}, s_{C,\gamma})}
$$
but 
$$
f_{(\gamma, \sigma_{\eta}, \sigma^{a}_{\gamma,\eta})}m_{(\eta^{\dagger}, \sigma^{a}_{\eta,\gamma}, \sigma_{C,\gamma})} = \pm f_{(\eta^{\dagger}, \sigma_{\eta}, s_{C_{1}})}m_{(\gamma, \sigma_{C_{1}}, \sigma_{C,\gamma})} = \pm f_{(\eta^{\dagger}, \sigma_{\eta}, \sigma_{C_{2}})}m_{(\gamma, \sigma_{C_{2}}, \sigma_{C,\gamma})}
$$
This is true also if we replace $a$ by $b$ in the superscripts. The sum in \ref{eqn:special} will equal the sum
$$
m_{(\eta,\sigma,\sigma_{\eta})}f_{(\eta^{\dagger}, \sigma_{\eta}, \sigma_{C_{1}})}m_{(\gamma, \sigma_{C_{1}}, \sigma_{C,\gamma})} + m_{(\eta,\sigma,\sigma_{\eta})}f_{(\eta^{\dagger}, \sigma_{\eta}, \sigma_{C_{2}})}m_{(\gamma, \sigma_{C_{2}}, \sigma_{C,\gamma})}
$$
where the first term is the $\eta$ contribution to $d(\lefty{e_{C_{1}}})m_{(\gamma,\sigma_{C_{1}},\sigma_{\gamma,C})}$ and the second term is the $\eta$ contribution to $d(\lefty{e_{C_{2}}})m_{(\gamma,\sigma_{C_{2}},\sigma_{\gamma,C})}$. \\
\ \\
\noindent This leaves us with the remaining case where $\eta = \gamma^{\dagger}$ in $m_{(\gamma,\sigma,\sigma_{\gamma})} d(\lefty{e_{C}})$. In this case we obtain the terms
$$
m_{(\gamma,\sigma,\sigma_{\gamma})}\big(f_{(\gamma^{\dagger}, \sigma_{\gamma}, \sigma_{C_{1}})}m_{(\gamma, \sigma_{C_{1}}, \sigma_{C, \gamma})}\big) + m_{(\gamma,\sigma,\sigma_{\gamma})}\big(f_{(\gamma^{\dagger}, \sigma_{\gamma}, \sigma_{C_{2}})}m_{(\gamma, \sigma_{C_{2}}, \sigma_{C,\gamma})}\big)
$$
If use associativity on each of the terms
$$
m_{(\gamma,\sigma,\sigma_{\gamma})}\big(f_{(\gamma^{\dagger}, \sigma_{\gamma}, \sigma_{C_{i}})}m_{(\gamma, \sigma_{C_{i}}, \sigma_{C,\gamma})}\big) = 
\big(m_{(\gamma,\sigma,\sigma_{\gamma})}f_{(\gamma^{\dagger}, \sigma_{\gamma}, \sigma_{C_{i}})}\big)m_{(\gamma, \sigma_{C_{i}}, \sigma_{C})}
$$
This is the $\gamma$ contribution to $d(\lefty{e_{C_{1}}})m_{(\gamma,\sigma_{C_{1}},\sigma_{C,\gamma})}$. The term for $i=2$ can likewise be refactored to correspond to a term in $d(\lefty{e_{C_{2}}})m_{(\gamma,\sigma_{C_{2}},\sigma_{C,\gamma})}$. $\Diamond$\\
\ \\

\noindent This completes the construction. From now on we will consider $ (\mathcal{B}\bridgeGraph{n}, d_{\Gamma_{n}})$ to be a differential, bigraded $\Z$-algebra.  

\subsection{Examples:} We give a complete description of $(\mathcal{B}\bridgeGraph{1}, d_{\Gamma_{1}})$ and a partial description of $(\mathcal{B}\bridgeGraph{2}, d_{\Gamma_{2}})$.\\
\ \\
\noindent $(\mathcal{B}\bridgeGraph{1}, d_{\Gamma_{1}})$: $P_{1}$ consists of two points, and there is only one planar matching in $\leftHalf$ and $\rightHalf$. Consequently, the only $1$-cleaved link is a circle intersecting the $y$-axis in two points. Thus, there are two vertices in $\bridgeGraph{1}$: when this circle is decorated with a $+$ and when it is decorated with a $-$. We will call these $C^{\pm}$. There are no bridges in either $\leftHalf$ or $\rightHalf$, so the only edges are $\lefty{e}_{C}: C^{+} \longrightarrow C^{-}$ and  $\righty{e}_{C}: C^{+} \longrightarrow C^{-}$. Thus $\bridgeGraph{1}$ looks like  

\begin{center}
\includegraphics[scale=0.5]{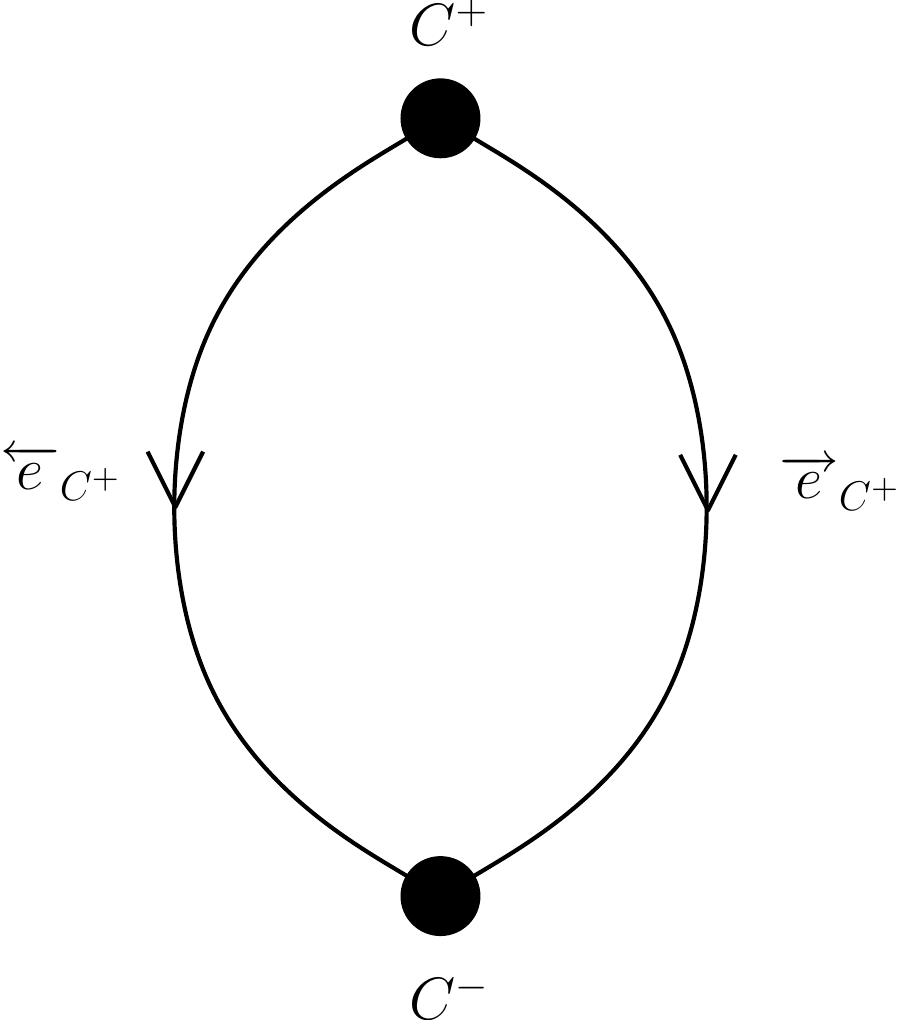}
\end{center}

\begin{figure}
\begin{center}
\includegraphics[scale=0.75]{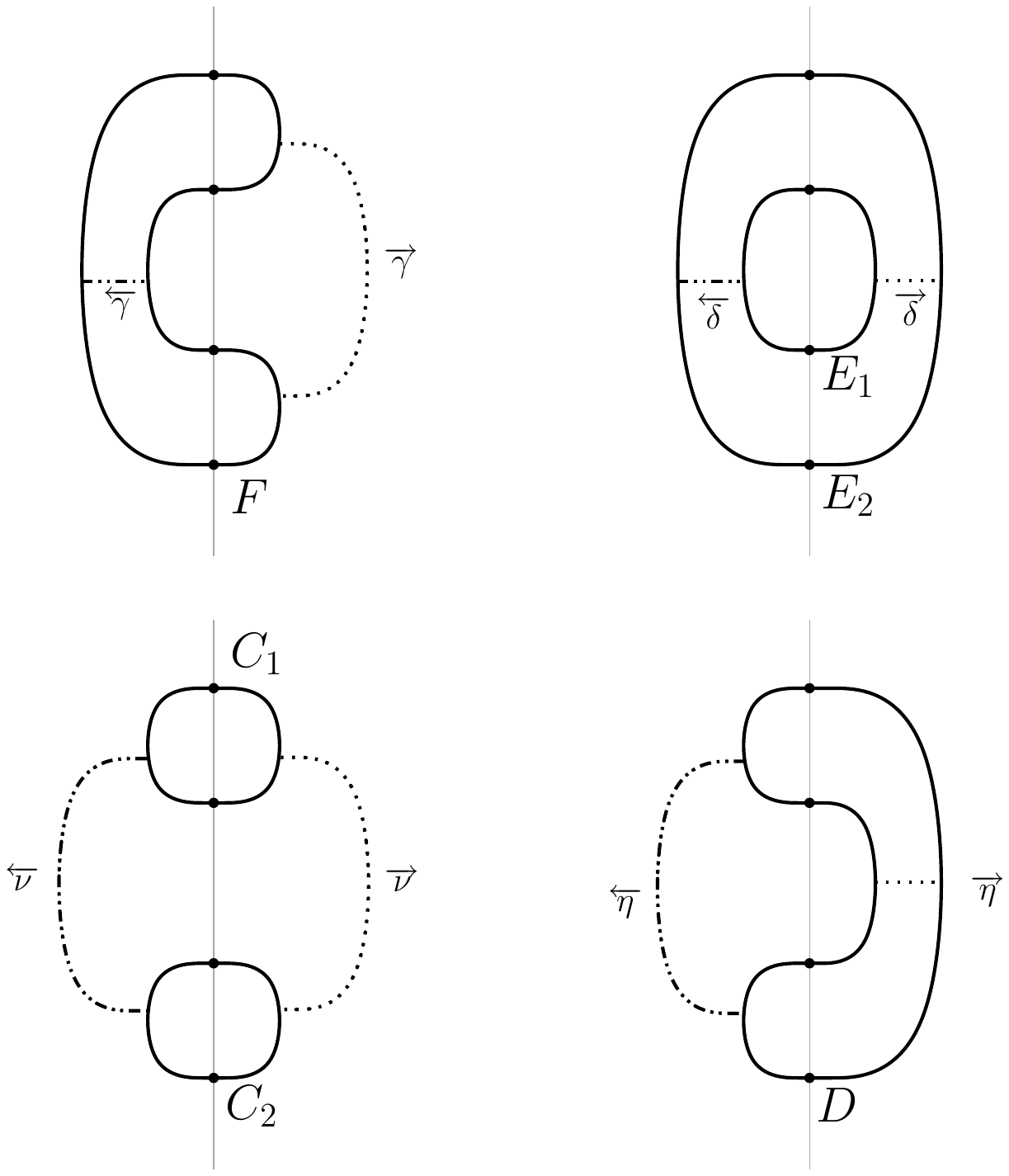}
\end{center}
\caption{The paired matchings on $P_{2}$.} 
\label{fig:bg2Verts}
\end{figure}

\noindent Thus, $\mathcal{B}\bridgeGraph{1}$ consists of four elements $I_{C^{+}}, I_{C^{-}}$ in grading $(0,0)$, $\lefty{e}_{C}$ in grading $(1,1)$, and $\righty{e}_{C}$ in grading $(0,-\frac{1}{2})$. The product of any two of these is trivial except for the actions of the idempotents: $I_{C^{+}}\lefty{e}_{C} = \lefty{e}_{C} = \lefty{e}_{C}I_{C^{-}}$, and similarly for $\righty{e}_{C}$. The differential $d_{\Gamma_{1}} \equiv 0$ since its image is in the set generated by paths of bridge edges.\\
\ \\
\noindent $(\mathcal{B}\bridgeGraph{2}, d_{\Gamma_{2}})$: $P_{2}$ consists of four points. There are two planar matchings $m_{1}$ and $m_{2}$ in this case. Thus there are four cleaved links $\lefty{m}_{1}\#\righty{m}_{1}$, $\lefty{m}_{1}\#\righty{m}_{2}$, $\lefty{m}_{2}\#\righty{m}_{1}$, $\lefty{m}_{2}\#\righty{m}_{2}$. $\lefty{m}_{1}\#\righty{m}_{2}$, $\lefty{m}_{2}\#\righty{m}_{1}$ will have only one circle, while $\lefty{m}_{1}\#\righty{m}_{1}$ and $\lefty{m}_{2}\#\righty{m}_{2}$ have two. These are depicted in Figure \ref{fig:bg2Verts}\\
\ \\

\begin{figure}
\begin{center}
\includegraphics[scale=0.8]{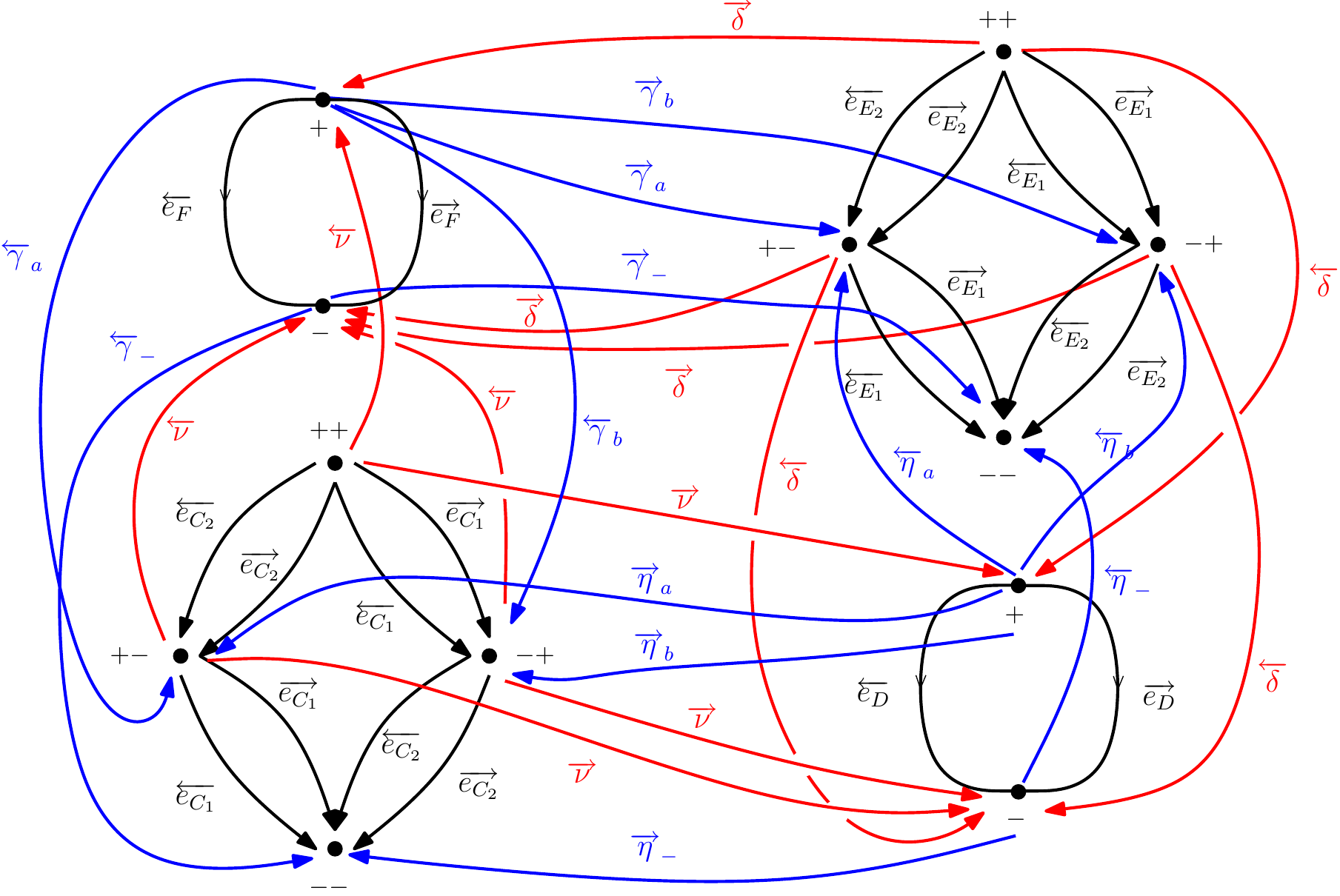}
\end{center}
\caption{A depiction of $\bridgeGraph{2}$. The vertices correspond to the circles in Figure \ref{fig:bg2Verts}. The decoration edges are shown, and labeled, in black. The bridge edges are shown in blue (for divides) and red (for merges), which correspondingly colored labels. To ease the clutter, the edges are labeled by the bridge and an arrow, using the notation from Figure \ref{fig:bg2Verts}. }
\label{fig:bg2Edges}
\end{figure}

\noindent Once we include decorations we will have $2 + 2 + 4 + 4 = 10$ different vertices in the graph, as shown by the black dots in \ref{fig:bg2Edges}. There are then 44 edges in the graph. As an example of some of the relations in $\mathcal{B}\Gamma_{2}$:
\begin{enumerate}
\item There are no relations corresponding to disjoint supports that involve bridges. For the decoration edges, we have some for the pairs of circles $\{E_{1},E_{2}\}$ and $\{C_{1},C_{2}\}$. For example, $\righty{e_{E_{1}}}\lefty{e_{E_{2}}} = (-1)^{0\cdot 1} \lefty{e_{E_{2}}}\righty{e_{E_{1}}}$, while $\lefty{e_{E_{1}}}\lefty{e_{E_{2}}} = (-1)^{1\cdot 1} \lefty{e_{E_{2}}}\lefty{e_{E_{1}}}$.
\item For the edges $\righty{e_{E_{1}}}$ and $\righty{e_{E_{2}}}$ we have, from the merge relations, 
$$
\righty{e_{E_{1}}}\righty{\delta} = \righty{e_{E_{2}}}\righty{\delta} = \righty{\delta}\righty{e_{F}}
$$
and
$$
\righty{e_{E_{1}}}\lefty{\delta} = \righty{e_{E_{2}}}\lefty{\delta} = \lefty{\delta}\righty{e_{D}}
$$
On the other hand, since $\righty{\gamma}$ is a divide, we obtain
$$
\righty{e_{F}}\righty{\gamma_{-}} = \righty{\gamma_{a}}\righty{e_{E_{1}}} = \righty{\gamma_{b}}\righty{e_{E_{2}}}
$$
\item Since $\righty{\eta} = \righty{\delta}^{\dagger}$ we get $\righty{\delta}\righty{\eta_{a}} = \righty{e_{E_{2}}}$ with $E_{2}$ being the active circle, while $\righty{\delta}\righty{\eta_{b}} = \righty{e_{E_{1}}}$.
\item We also have relations such as $\righty{\gamma_{a}}\lefty{\delta} = \lefty{\gamma_{a}}\righty{\nu}$ and $\righty{\delta}\righty{\eta_{a}} = \lefty{\gamma_{b}}\righty{\nu}$ that come from pairs of bridges in the same diagram. Due to our labels in Figure \ref{fig:bg2Edges}, these look slightly different than above.
\item Finally, we give an example of $d_{\Gamma_{2}}$:
$$
d_{\Gamma_{2}}\lefty{e_{F}} = -\lefty{\gamma_{a}}\lefty{\nu} - \lefty{\gamma_{b}}\lefty{\nu}
$$
whereas
$$
d_{\Gamma_{2}}\lefty{e_{C_{2}}} = - \lefty{\nu}\lefty{\gamma_{a}}
$$
\end{enumerate}

\section{Tangles and Resolutions}\label{sec:APS}

\noindent  Let $\lefty{\mathbb{R}^{3}} = \leftHalf \times \R$ be the half space corresponding to $\leftHalf \subset \R^2$ under the standard projection $\pi$ to the $xy$-plane. Let $\righty{\mathbb{R}^{3}}$ be defined similarly. 

\begin{defn}
An (outside) tangle $\righty{\tangle{T}}$ is a smooth, proper embedding of 
\begin{enumerate}
\item[] i) $n = \numarcs{\righty{\tangle{T}}}$ copies of the interval $[0,1]$, and 
\item[] ii) $\numcircs{\righty{\tangle{T}}}$ copies of $S^{1}$
\end{enumerate}
in $\righty{\mathbb{R}^{3}}$, whose boundary is the set of $2n$ points $P_{n}$ in $\partial \rightHalf$. $\righty{\tangle{T}_{1}}$ and $\righty{\tangle{T}_{2}}$ are equivalent if there is an isotopy of $\righty{\mathbb{R}^{3}}$ taking $\righty{\tangle{T}_{1}}$ to $\righty{\tangle{T}_{2}}$ and pointwise fixing the boundary $\partial\righty{\mathbb{R}^{3}}$.  
\end{defn}

\noindent{\bf Convention:} We will generally work in $\righty{\mathbb{R}^{3}}$ to describe our constructions. However, similar definitions apply for $\lefty{\mathbb{R}^{3}}$. Tangles in $\lefty{\mathbb{R}^{3}}$ will be called {\em inside} tangles. The names come from using the outward-first normal convention for defining the orientation on the boundary of an oriented manifold. If we take the orientation of the $y$-axis to come from the usual coordinates, then the $y$-axis is the boundary of $\leftHalf$, so $\rightHalf$ is the ``outside'' and $\leftHalf$ is the ``inside''.\\ 
\ \\
\noindent As usual, each inside tangle $\righty{\mathbb{R}^{3}}$ can be isotoped to have a generic projection to the $xy$-plane. The image will be a four valent graph in the interior of $\rightHalf$ with one valent vertices at the points $P_{n} \subset \{0\} \times \R$. We record the over/under-crossing of the strands of $\righty{\tangle{T}}$ to obtain a tangle diagram for $\righty{\tangle{T}}$. Different diagrams for $\righty{\tangle{T}}$ are related by sequences of Reidemeister moves, and planar isotopies, in the interior of $\rightHalf$. We will denote a tangle diagram for a tangle by the corresponding roman letter:  $\righty{T}$ will be a diagram for $\righty{\tangle{T}}$. \\
\ \\
\noindent The crossings of $\righty{T}$ form a set $\cross{\righty{T}}$. If we orient $\righty{\tangle{T}}$ then each crossing in a diagram $\righty{T}$ is either a positive crossing, or a negative crossing according to the usual rule: \\
$$
\inlinediag[0.5]{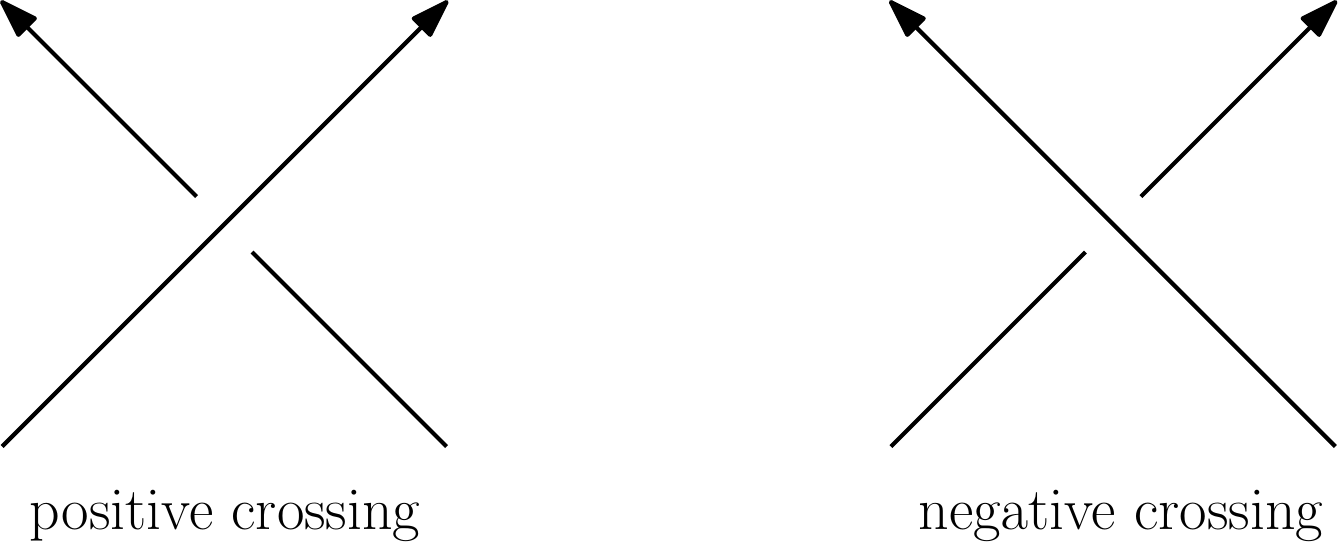} 
$$
\noindent We will denote by $n_{+}(\righty{T})$ the number of positive crossings, and $n_{-}(\righty{T})$ the number of negative crossings. 

\subsection{The tangle homology of Asaeda, Przytycki, and Sikora}
\noindent In \cite{APS} M. Asaeda, J. Przytycki, and A. Sikora define a reduced Khovanov homology for tangles, which we will describe using our conventions. Suppose we have an oriented tangle diagram $\righty{T}$ in $\rightHalf$, and we have {\em ordered its crossings}.

\begin{defn}
An APS-resolution $\rho$ of $\righty{T}$ is a map $\rho: \cross{\righty{T}} \lra \{0,1\}$. For each resolution, $\rho$, there is a resolution diagram, $\rho(\righty{T})$. $\rho(\righty{T})$ is the crossingless, planar link in $\rightHalf$ obtained by locally replacing (disjoint) neighborhoods of the crossings using the following rule for each crossing $c \in \cross{\righty{T}}$:\\
$$
\inlinediag[0.5]{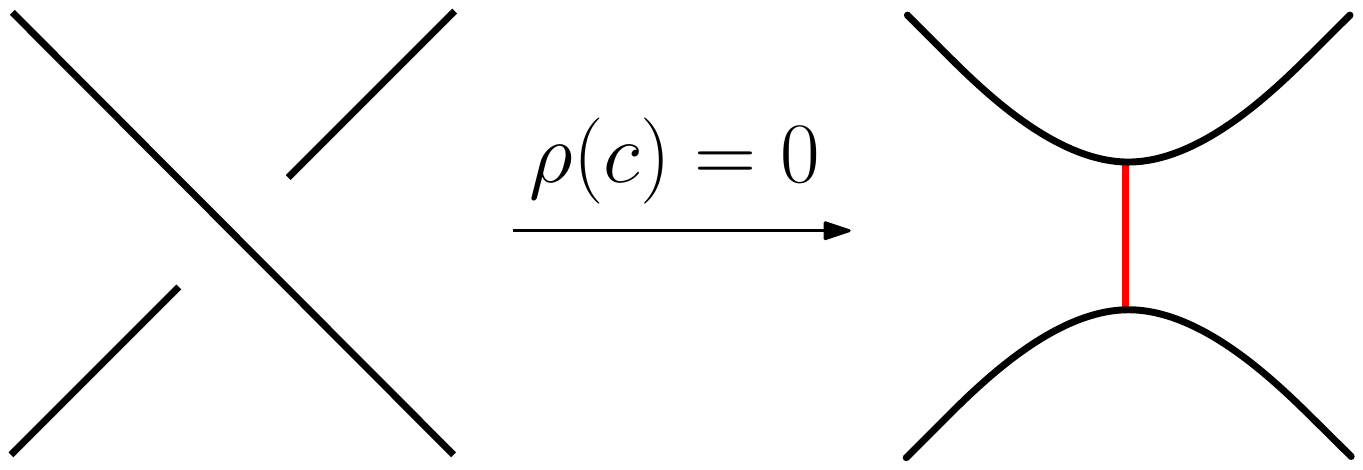} \hspace{.5in} \inlinediag[0.5]{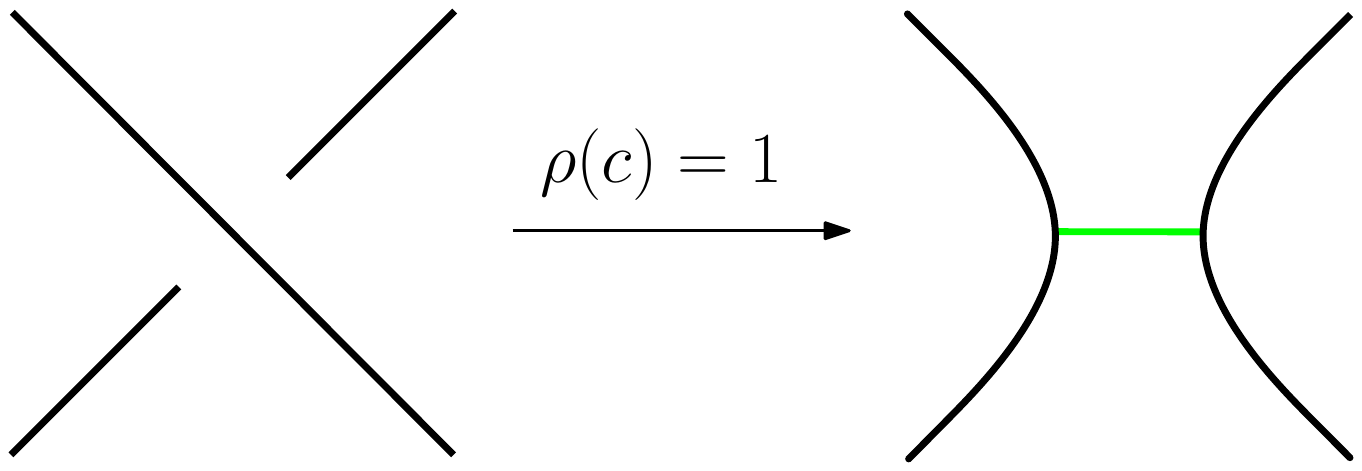} 
$$ 
\noindent The set of resolutions will be denoted $\resolution{\righty{T}}$. 
\end{defn}

\noindent The local arcs introduced by $\righty{T} \rightarrow \rho(\righty{T})$ are called {\em resolution bridges}. When a diagram $\righty{T}$ is understood we will not distinguish between $\rho$ as a map and $\rho(\righty{T})$ as a diagram equipped with resolution bridges. \\
\ \\
\noindent  The resolution bridges for a crossing $c \in \cross{\righty{T}}$ will be denoted $\gamma_{\rho, c}$ (or just $\gamma_{c}$ when the resolution is understood). If $\rho(c) = 0$ we will call $\gamma_{c}$ an {\em active} bridges for $\rho$, while if $r(c)=1$ it will be called {\em inactive}. The active bridges for $\rho$ are the elements of the set $\actor{\rho}$. We will denote by $\gamma_{c}(\rho)$ the resolution obtained by surgering $\rho(T)$ along $\gamma_{c}$. This resolution has $\gamma_{c}(\rho)(c) = 1 - \rho(c)$ and $\gamma_{c}(\rho)(\tilde{c})=\rho(\tilde{c})$ for $\tilde{c} \in \big(\cross{T}\backslash\{c\}\big)$. The resolution bridges at $c$ for $\gamma_{c}(\rho)$ will be denoted by $\gamma_{c}^{\dagger}$ when considered from $\rho$. \\
\ \\
\begin{defn}
For an APS resolution $\rho$, let
$$
h(\rho) = \sum_{c \in \cross{T}} \rho(c)\\
$$
\end{defn}

\noindent A resolution diagram $\rho(\righty{T})$ consists of a planar diagram of circles and bridges. The construction in \cite{APS} starts by assigning to this diagram a bi-graded free module over $\F$. Let $V = \F\,v_{+} \oplus \F\,v_{-}$ be the rank $2$ module generated by $v_{\pm}$ with bigrading $\mathrm{gr}(v_{\pm}) = (0, \pm 1)$. The {\em homological} grading will be the first entry in this order pair, while the second entry is the {\em quantum grading}.
We also let $W= \F\,v_{0}$ be the rank $1$ module in bigrading $(0,0)$. As we will need to shift the bigrading, we introduce our notation for these shifts:

\begin{defn}
If $V = \bigoplus_{i,j} V_{i,j}$ is a $\Z$-bigraded module, then $V[(\delta ,\eta)]$ is the bigraded module with $\big(V[(\delta ,\eta)]\big)_{i,j} = V_{i - \delta, j - \eta}$. 
\end{defn}

\noindent To each circle $C$ in $\rho(\righty{T})$ we assign $V_{C}$, a copy of $V$. To each arc $A$ we assign $W_{A}$, a copy of $W$. To the resolution $\rho$ we assign $V_{\rho}$ the bigraded tensor product of these modules, shifted in bigrading by $[(h(\rho),h(\rho))]$.\\
\ \\
\noindent For each active bridge $\gamma_{\rho,c}$ there is a $(1,0)$ $\F$-linear map $D_{\gamma,\rho}: V_{\rho} \rightarrow V_{\gamma_{c}(\rho)}$ defined by the following recipe.
\begin{enumerate}
\item If surgery on $\gamma_{c}$ merges the circles $C_{1}$ and $C_{2}$ in $\rho$ to get a circle $C$ in $\gamma_{c}(\rho)$, then we use
$$
\mu = \left\{\begin{array}{l}
v_{+} \otimes v_{+} \longrightarrow v_{+}[(1,1)]\\
v_{+} \otimes v_{-} \longrightarrow v_{-}[(1,1)]\\
v_{-} \otimes v_{+} \longrightarrow v_{-}[(1,1)]\\
v_{-} \otimes v_{-} \longrightarrow 0     \\
\end{array}\right. 
$$
to define a map $V_{C_{1}} \otimes V_{C_{2}} \rightarrow V_{C}[(1,1)]$, then tensoring this with the identity applied to each of the other factors, and shifting the bigrading by $[(h(\rho),h(\rho)]$.
\item If $\gamma_{c}$ has both ends on the same circle $C$, and surgering $C$ along $\gamma_{c}$ produces two new circles $C_{1}$ and $C_{2}$, then on $V_{C}$ we use the coproduct $\Delta : V_{C} \longrightarrow V_{C_{1}}\otimes V_{C_{2}}$ given by
$$
\Delta = \left\{\begin{array}{l}
v_{+} \longrightarrow v_{+} \otimes v_{-}[(1,1)] + v_{-} \otimes v_{+}[(1,1)] \\
v_{-} \longrightarrow v_{-} \otimes v_{-}[(1,1)]
\end{array}\right. 
$$ 
We again tensor with the modules for the remaining circles and shift accordingly.
\item
\noindent If $\gamma_{c}$  has both feet on an arc $A$, $\gamma_{c}(\rho)$ will have a new circle component $C$. We then use the map $\delta: W_{A} \rightarrow W_{A'}\otimes V_{C}[(1,1)]$ given by 
$$
\delta(v_{A,0}) = v_{A',0} \otimes v_{C,-}[(1,1)]
$$
and tensor it with the identity on the other factors.
\item If $\gamma_{c}$ has on foot on an arc $A$ and the other foot on a circle $C$ we use the map $\lambda$: 
$$
\begin{array}{ll}
\lambda(v_{C, +} \otimes v_{A,0}) = v_{A',0}[(1,1)]\\
\lambda(v_{C, -} \otimes v_{A,0}) = 0 
\end{array}
$$ 
on the factors corresponding to $C$ and $A$, and the identity on the other factors.
\item In all other cases, the map is taken to be the zero map. 
\end{enumerate}

\noindent We now define the chain groups to be the bigraded module
$$
C_{\mathrm{APS}}(\righty{T}) = \bigoplus_{\rho \in \resolution{\righty{T}}} V_{\rho}
$$
equipped with the $(1,0)$-map $d$ given by
$$
d\big|_{V_{\rho}} = \sum_{\gamma \in \actor{\rho}} (-1)^{I(\rho,\gamma)}\,D_{\gamma,\rho}
$$
where $I(\rho,\gamma) = \sum_{c_{\gamma} < c'} \rho(c')$ is the number of $\rho$-inactive crossings which occur after the crossing corresponding to $\gamma$.  \\
\ \\
\noindent The main result of \cite{APS} is

\begin{thm}
The pair $(C_{\mathrm{APS}}(\righty{T}), d_{\mathrm{APS}})$ is a chain complex. In addition, the isomorphism class of $H_{\ast}(\righty{CK}(C_{\mathrm{APS}}(\righty{T}), d_{\mathrm{APS}})$, as a relatively bigraded module, is an invariant of the isotopy class of $\righty{\mathcal{T}}$.  
\end{thm} 

\noindent Using the orientation on $\righty{T}$, we can shift again to get $(C_{\mathrm{APS}}(\righty{T}), d_{\mathrm{APS}})[(-n_{-},n_{+} - 2 n_{-})]$. The homology of this complex is an isotopy invariant of the oriented tangle $\righty{T}$ as a bi-graded module.

\subsection{The expanded tangle homology}

\noindent We can use $(C_{\mathrm{APS}}(\righty{T}), d_{\mathrm{APS}})$ to construct a larger chain complex, $\rightComplex{\righty{T}}$. To this end we start by picking a matching $m \in \match{n}$. From $(C_{\mathrm{APS}}(\righty{T}), d_{\mathrm{APS}})$ we can pick out all the $V_{\rho}$ where the arcs in $\rho(\righty{T})$ represent the matching $m$. A careful review of 
the construction of the maps $\lambda$ and $\delta$ show that they do not alter the matching $m$, since they merge or divide circles. Thus, $d_{\mathrm{APS}}$ is the direct sum of maps which preserve the matchings. Consequently, 
$$
(C_{\mathrm{APS}}(\righty{T}), d_{\mathrm{APS}}) = \bigoplus_{m \in \match{n}} (C_{\mathrm{APS}}(\righty{T}, m), d_{\mathrm{APS, m}})
$$
which we obtain by considering $m$ in $\rightHalf$. \\
\ \\
\begin{defn}
Given $(L,\sigma) \in \cleaved{n}$ let
$$
\righty{CK}(\righty{T}, L, \sigma) = (C_{\mathrm{APS}}(\righty{T}, \righty{L}), d_{\mathrm{APS, \righty{L}}})[(0,\frac{1}{2}\iota(L,\sigma))]
$$
and
$$
\rightComplex{\righty{T}} = \bigoplus_{(L,\sigma) \in \cleaved{\numarcs{\righty{T}}}} \righty{CK}(\righty{T}, L, \sigma)
$$
\end{defn}

\noindent The invariance for each summand $(C_{\mathrm{APS}}(\righty{T}, \righty{L}), d_{\mathrm{APS, \righty{L}}})$ immediately implies that $\rightComplex{\righty{T}}$ is an invariant as well. Furthermore, if  $\righty{T}$ is oriented, the homology of $$\righty{CK}(\righty{T}, L, \sigma)\{(-n_{-}, n_{+} - 2 n_{-})\}$$, considered as an absolutely bi-graded module, is a isotopy invariant of an oriented tangle in $\righty{\R^{3}}$. This is true whether or not we can extend the orientation to $L$. However, at this point, $\rightComplex{\righty{T}}$ is just a stack of copies of the summands of $C_{\mathrm{APS}}(\righty{T})$. We will improve on this situation in the coming sections. First, we provide some more insight into this construction. \\
\ \\
\noindent Let $m \in \match{n}$ and consider $m$ to be embedded in $\leftHalf$. In this case, we will write $\lefty{m}$ and call $\lefty{m}$ and inside matching. We can use $\lefty{m}$ to close up the tangle $\righty{T}$.

\begin{defn}\ \\
For an outside tangle $\righty{T}$ and an inside matching $\lefty{m}$, $\lefty{m}\#\righty{T}$ will denote the link in $\R^{3}$ found by gluing the ends of $\righty{T}$ to those of $\lefty{m}$. If $T$ is a projection of $\righty{T}$, then $\lefty{m} \# T$ is the projection of $\lefty{m} \# \righty{T}$.
\end{defn}
\ \\
\noindent{\bf Warning:} {\em We will consider $\lefty{m} \# \righty{T}$ up to isotopy of $\R^{3}$ supported in $\mathrm{int\,}\righty{\R^{3}}$} and not up to ambient isotopy in $\R^{3}$. In particular, $\lefty{m}$ will be taken to be a specific representative of its equivalence class as a matching, and will never again be moved. \\
\ \\
\noindent We now enlarge our set of resolutions to account for the possible inside matchings:

\begin{defn}
A resolution $r$ of $\righty{T}$ is a pair $(\rho, \lefty{m})$ where $\rho: \cross{\righty{T}} \lra \{0,1\}$ is an APS-resolution, and $\lefty{m}$ is an inside matching. The resolution diagram, $r(\righty{T})$ is the crossingless, planar link $\lefty{m}\#\rho(\righty{T})$ 
The set of resolutions will be denoted $\resolution{\righty{T}}$.
\end{defn}

\noindent Thus a resolution is simply a way to close up the arcs of an APS-resolution using a matching in $\leftHalf$. We will use the same notation as for APS resolutions, when working in $\rightHalf$.\\
\ \\
\noindent A resolution diagram $r(\righty{T})$ consists of a planar diagram of circles. These circles form a set $\circles{r}$ which is partitioned into those circles from $\rho(\righty{T})$ that are contained in $\mathrm{int\,}\rightHalf$ and those circles intersecting the $y$-axis. The circles from $\rho(\righty{T})$ will be called {\em free circles}, and form the set $\free{r}$, while the remaining circles will be called cleaved circles, and together comprise an element $\mathrm{cl}(r)$ of $\cleave{n}$ for $n=\numarcs{\righty{\tangle{T}}}$.  

\begin{defn}
A {\em state} for $\righty{T}$ is a pair $(r, s)$ where
\begin{enumerate}
\item $r$ is a resolution of $\righty{T}$,
\item $s$ is an assignment of an element of $\{+,-\}$ to each circle of $r(\righty{T})$. This assignment will be called a decoration on $r(\righty{T})$. 
\end{enumerate}
The states for $\righty{T}$ will be denoted $\state{T}$.
\end{defn}

\begin{defn}
The boundary of a state $(r, s)$ for $\righty{T}$ is the element $\partial(r,s) =(L,\sigma) \in \cleaved{n}$ consisting of the pair $L = \mathrm{cl}(r)\in \cleave{n}$ and the decoration $\sigma= s|_{L}$. 
\end{defn}

\noindent To a state $(r,s)$ we thus obtain an element of $\cleave{n}$. We also obtain a generator of $C_{\mathrm{APS}}(\righty{T})$ by assigning $v_{pm}$ to each circle $C \in \free{r}$ according to the value of $s(C)$, and assigning $v_{0}$ to the arcs. Thus we can think of 
$$
\righty{CK}(\righty{T}, L, \sigma) \cong \bigoplus_{\partial(r,s) = (L,\sigma)} \F\!(r,s)
$$
if we disregard the bigrading. It is straightforward to add back in the shifts from $C_{APS}$. The only additional information is the additional quantum grading shift $[(0,\iota(L,\sigma)]$ which corresponds to half that if we assigned $v_{\pm}$ to each cleaved circle according to the value of $\sigma$, instead of using $v_{0}$.\\
\ \\
\noindent Thus we can think of $\rightComplex{\righty{T}}$ as generated by the states $(r,s)$ where $r = (\rho, \lefty{m})$ corresponds to the planar link $\lefty{m} \# \rho$, and $s$ prescribes the generator from $V$ to be assigned to each circle. Cleaved circles each count for half the quantum grading of a free circle. However, we use the APS-boundary map, without changing the decoration on each cleaved circle.  \\
\ \\
\noindent Since $\lambda$ and $\delta$ preserve the image of $(r,s) \rightarrow (L,\sigma)$ and act identically regardless of the choice of $\sigma$, we have $\righty{CK}(\righty{T}, L, \sigma) \cong \righty{CK}(\righty{T}, L', \sigma')$ when $L' \cap \rightHalf = L \cap \rightHalf$ and we consider the complexes as relatively bigraded.\\
\ \\
\noindent We encode this decomposition by defining an action of $\mathcal{I}_{n_{A}}$ on $\rightComplex{T}$:
$$
I_{(L,\sigma)} (r,s) = \left\{\begin{array}{ll} (r,s) & \partial(r,s) = (L,\sigma) \\ 0 & \mathrm{else} \end{array} \right.
$$
Thus $I_{(L,\sigma)}$ acts non-trivially only on the summand $\righty{CK}(\righty{T}, L, \sigma)$.

\subsection{Some more on states}

\noindent For a state $(r,s)$, we can use $s$ to group the bridges in $\actor{r}$ into (overlapping) classes: 
\begin{enumerate}
\item $\interior{r,s}$ is the subset of active bridges $\gamma_{c}$ such that either
\begin{enumerate}
\item both feet of $\gamma_{c}$ are on elements of $\free{r}$, or
\item one foot of $\gamma_{c}$ is on $C \in \mathrm{cl}(r)$ and the other foot is on $C' \in \free{r}$ with $s(C') = +$, or
\item both feet are on the same arc of $C \cap \rightHalf$ for some $C \in \mathrm{cl}(r)$
\end{enumerate}
These are the active bridges that can make contributions to $d_{APS}$.
\item $\dec{r,s}$ is the subset of active bridges $\gamma_{c}$ such that either
\begin{enumerate}
\item both feet are on the same arc of $C \cap \rightHalf$ for some $C \in \mathrm{cl}(r)$ with $s(C) = +$, or
\item one foot of $\gamma_{c}$ is on $C \in \mathrm{cl}(r)$ with $s(C) = +$ and the other foot is on $C' \in \free{r}$ with $s(C') = -$
\end{enumerate}
\item $\rightBridges{r}$ is the subset of active bridges $\gamma_{c}$ such that either
\begin{enumerate}
\item $\gamma_{c}$ has one foot on $C_{1} \in \mathrm{cl}(r)$ and the other on a distinct circle $C_{2} \in \mathrm{cl}(r)$, or
\item $\gamma_{c}$ has both feet on some $C \in \mathrm{cl}(r)$, but they are on different arcs in $C \cap \rightHalf$.
\end{enumerate}
\end{enumerate}

\noindent To clarify the notation, we note that there is a natural map $\rightBridges{r} \longrightarrow \rightBridges{\mathrm{cl}(r)}$. If $r = (\rho, \lefty{m})$ we
will let $\leftBridges{r} = \leftBridges{\lefty{m}}$ and $\bridges{r} = \leftBridges{\lefty{m}} \cup \rightBridges{r}$.

\section{Type $D$-structures}
\label{sec:typeD}

\noindent We review the definition of a type $D$ structure given in \cite{Bor1}. Unlike in \cite{Bor1}, where the coefficients are always from a characteristic 2 field, we will give the definition with signs corresponding to the conventions of this paper.

\begin{defn}
Let $W = W_{0} \oplus W_{1}$ be a $\Zmod{2}$-graded module.  $|\I_{W}| : W \rightarrow W$ is the signed identity defined by
linearly extending
$$
|\I_{W}|(w) = (-1)^{\mathrm{gr}(w)}w
$$
for homogeneous $w \in A$.  
\end{defn}

\begin{defn}
Let $(A,d)$ be a differential, unital, $\Z$-graded algebra over a ring $R$, and let $N$ be a module over a ring $R$. A {\em (left) $D$-structure} on $N$ is a linear map
\begin{equation}
\delta : N \longrightarrow \big(A \otimes_{R} N\big)\,[-1]
\end{equation}
such that
$$
(\mu_{A} \otimes \I_{N})\,(\I_{A} \otimes \delta)\,\delta + (d \otimes |\I_{N}|)\,\delta \equiv 0
$$
A {\em morphism} of $D$-structures $(N,\delta) \rightarrow (N',\delta')$ is a map $\psi : N \rightarrow A \otimes N'$ satisfying
$$
(\mu_{A} \otimes \I_{N})\,(\I_{A} \otimes \delta')\,\psi - (\mu_{A} \otimes \I_{N})\,(\I_{A} \otimes \psi)\,\delta + (d \otimes |\I_{N}|)\,\psi \equiv 0
$$
The {\em identity morphism} is the map $I_{N} : N \rightarrow A \otimes N$ given by $x \rightarrow 1_{A} \otimes x$ where $1_{A}$ is the identity in $A$.
If $\psi$ is a $D$-structure morphism from $N$ to $N'$ and $\phi$ is a {\em $D$-structure morphism} from $N'$ to $N''$, then $\phi \ast \psi$ is the morphism
$$
\phi \ast \psi = (\mu_{A} \otimes \I_{N})\,(\I_{A} \otimes \phi)\,\psi : N \rightarrow A \otimes N''
$$
Two morphisms $\psi, \phi : N \rightarrow A \otimes N'$ are {\em homotopic}, $\psi \simeq \phi$, if there is a map $H: N \rightarrow \big(A \otimes N'\big)[1]$
$$
\psi - \phi = (\mu_{A} \otimes \I_{N})\,(\I_{A} \otimes H)\,\delta + (\mu_{A} \otimes \I_{N})\,(\I_{A} \otimes \delta')\,H + (d \otimes |\I_{N}|)\,H
$$
$(N,\delta)$ and $(N',\delta')$ are {\em homotopy equivalent} if there are $D$-structure morphisms $\psi: N \rightarrow A \otimes N'$ and $\phi : N' \rightarrow A \otimes N$ such that $\phi \ast \psi \simeq I_{N}$ and $\psi \ast \phi \simeq I_{N'}$.
\end{defn}

\noindent It is straightforward to verify the following properties:
\begin{enumerate}
\item $I_{N} \ast \psi = \psi \ast I_{N} = \psi$
\item $\psi \ast(\phi \ast \nu) = (\psi \ast \phi) \ast \nu$
\item If $\psi \simeq_{H} \psi'$ then $\phi \ast \psi \simeq_{\phi \ast H} \phi \ast \psi'$ and $\psi \ast \nu \simeq_{H \ast \nu} \psi' \ast \nu$ 
\end{enumerate}

\noindent From this it readily follows that

\begin{prop}
Homotopy equivalence of $D$-structures is an equivalence relation.
\end{prop}
\ \\
\noindent Let $N$ be a free graded $\Z$-module with basis $B=\{x_{1}, \ldots, x_{n}\}$ and suppose $a_{ij} \in A$ so that $\mathrm{gr}(a_{ij}) = |x_{i}|-|x_{j}|+1$ and
\begin{equation}\label{eqn:structure}
(-1)^{|x_{k}|} d(a_{ik}) + \sum_{j=1}^{n} a_{ij}\,a_{jk} = 0  \hspace{0.5in} i,k \in \{1, \ldots, n\}
\end{equation}
Then the $a_{ij}$ can be used as structure coefficients in the definition of a map $\delta : N \rightarrow \big(A \otimes_{R} N\big)[-1]$
$$
\delta(x_{i}) = \sum_{j=1}^{n} a_{ij} \otimes x_{j}
$$
The relation on the structure coefficients expressed in equation \ref{eqn:structure} ensures that $\delta$ is a type $D$-structure for $A$ and $N$.\\
\ \\
\noindent Note that the description of $(N,\delta)$ via structure coefficients will not necessarily be unique, even for a given $B$. However, for a type $D$-structure $(N,\delta)$ described by structure coefficients in an approproiate way, the following proposition will allow us to simplify $\delta$. 

\begin{prop}[Cancellation]\label{prop:cancelD}
Let $\delta$ be a $D$-structure on $N$. Suppose there is a basis $B$ for $N$ where $\delta$ can be described by structure coefficients satisfying $a_{ii} = 0$ and $a_{12} = 1_{A}$. Let $\overline{N} = \mathrm{span}_{R}\{\overline{x}_{3}, \ldots, \overline{x}_{n}\}$. Then 
$$
\overline{\delta}(\overline{x}_{i}) = \sum_{j \geq 3}(a_{ij} - a_{i2}\,a_{1j}) \otimes \overline{x}_{j}
$$
defines a $D$-structure on $\overline{N}$. Furthermore, the maps
$$
\begin{array}{lcl}
\iota : \overline{N} \rightarrow A \otimes N & \hspace{0.75in} & \iota(\overline{x}_{i}) = 1_{A} \otimes x_{i} -  a_{i2} \otimes x_{1} \\
\ &\ \\
\pi :  N \rightarrow A \otimes \overline{N} & \hspace{0.75in} & \pi(x_{i}) = \left\{
\begin{array}{ll}
0 & i = 1\\
-\sum_{j \geq 3}a_{1j} \otimes \overline{x}_{j} & i = 2\\ 
1_{A} \otimes \overline{x}_{i} & i \geq 3
\end{array}
\right.
\end{array}
$$
realize $\overline{N}$ as a deformation retraction of $N$ with $\iota \ast \pi \simeq_{H} \I_{N}$ using the {\em homotopy} $H: N \rightarrow A \otimes N[-1]$
$$
H(x_{i}) = \left\{\begin{array}{ll} -1_{A} \otimes x_{1} & i = 2\\ 0 & i \neq 2 \end{array}\right.
$$ 
\end{prop}

\noindent Although this result is likely known to others, we need to verify it using the sign conventions of this paper. However, as this verification is unrelated to the remainder of the paper, the proof of proposition \ref{prop:cancelD} is given in  appendix \ref{sec:cancelDProof}.

\begin{defn}
$\overline{N}$ will be called the deformation of $N$ found by canceling $x_{1}$ and $x_{2}$.
\end{defn}

\noindent We can also cancel when $a_{12}$ is a unit in the base ring $R$:

\begin{prop}
Let $\delta$ be a $D$-structure on $N$ over a ring $R$. Suppose there is a basis $B$ for $N$ where $\delta$ can be described by structure coefficients satisfying $a_{ii} = 0$ and $a_{12} = u$ for a unit $u \in R$. Let $\overline{N} = \mathrm{span}_{k}\{\overline{x}_{3}, \ldots, \overline{x}_{n}\}$. Then 
$$
\overline{\delta}(\overline{x}_{i}) = \sum_{j \geq 3}(a_{ij} - a_{i2}\,u^{-1}\,a_{1j}) \otimes \overline{x}_{j}
$$
is a $D$-structure on $\overline{N}$. Furthermore, $(\overline{N}, \overline{\delta})$ is a deformation retraction of $(N,\delta)$.
\end{prop}

\noindent The homotopy equivalence $\iota$ is then $\iota(\overline{x}_{i}) = 1_{A} \otimes x_{i} - a_{i2}u^{-1} \otimes x_{1}$.\\
\ \\

\section{The type D-complex for an outside tangle}\label{sec:typeDcom}

\begin{figure}
\begin{center}
\includegraphics[scale=0.6]{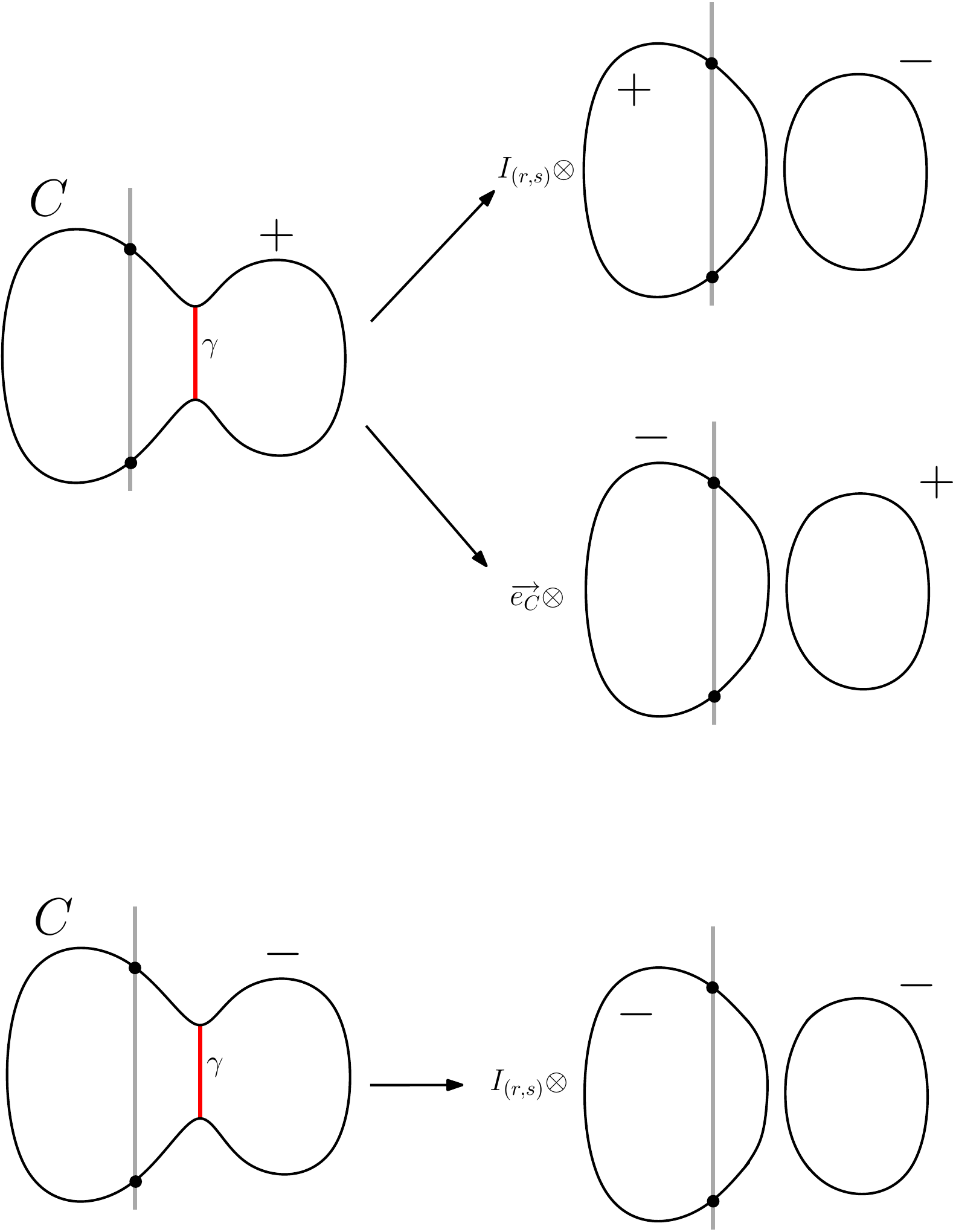}
\end{center}
\caption{}
\label{fig:rightdec}
\end{figure}

\begin{figure}
\begin{center}
\includegraphics[scale=0.6]{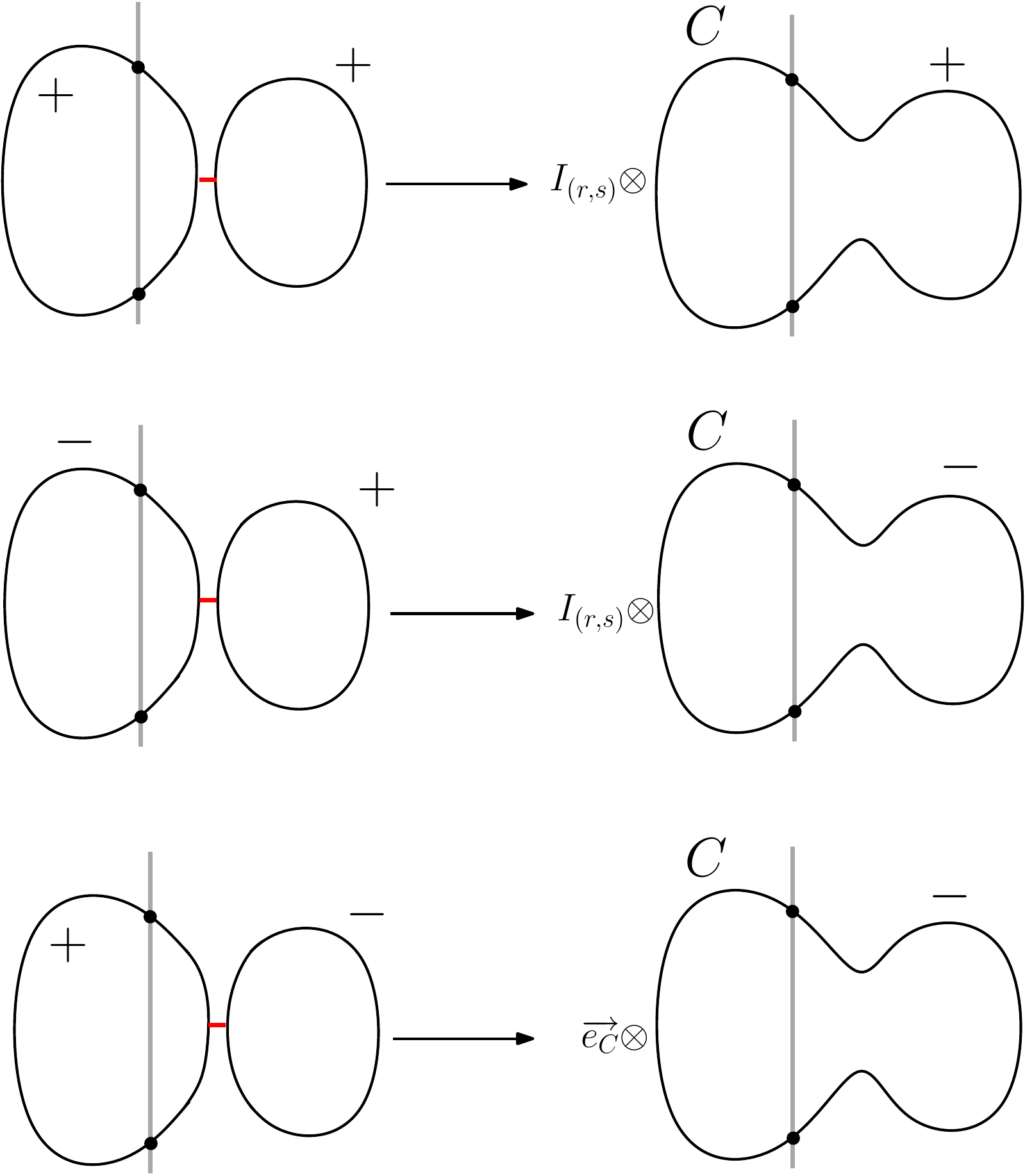}
\end{center}
\caption{}
\label{fig:rightdec2}
\end{figure}

\begin{figure}
\begin{center}
\includegraphics[scale=0.6]{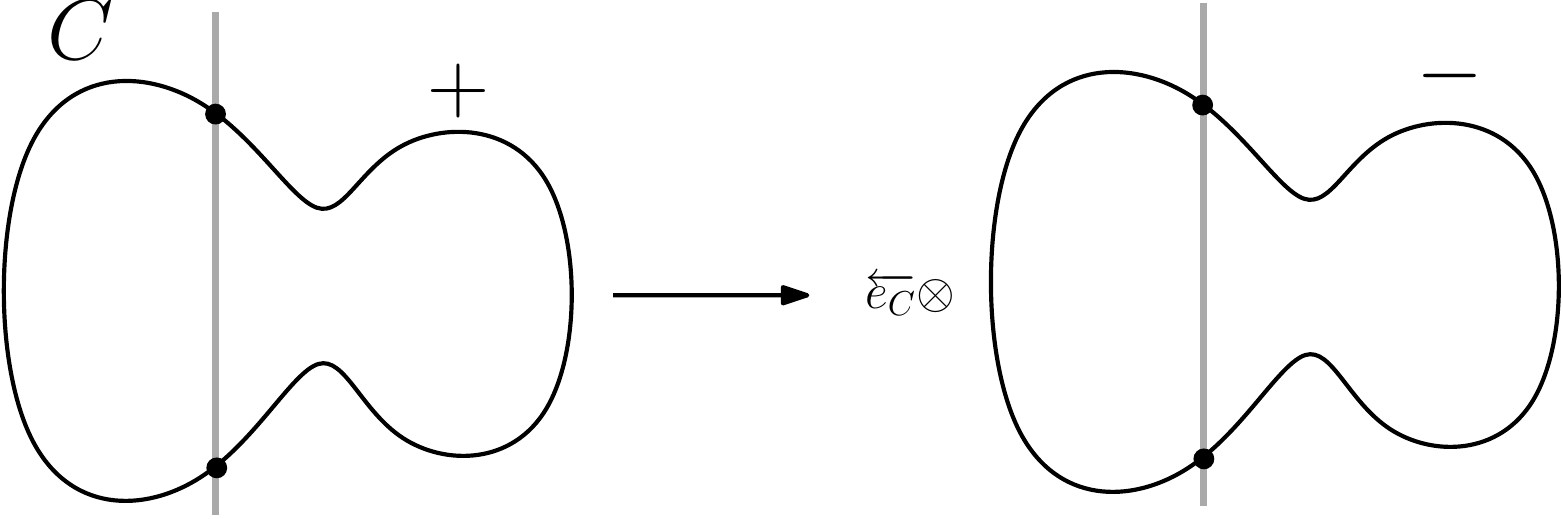}
\end{center}
\caption{}
\label{fig:leftdec}
\end{figure}

\noindent Given $\righty{T}$ our aim now is to describe a $\mathcal{I}_{n}$-{\em module} map
$$
\righty{\delta_{T}} : \rightComplex{\righty{T}} \longrightarrow \mathcal{B}\bridgeGraph{n} \otimes_{\mathcal{I}} \rightComplex{\righty{T}} [(-1,0)]
$$
which defines a right $D$-structure on $\rightComplex{\righty{T}}$.  \\
\ \\
\noindent We will define $\righty{\delta_{T}}$ by specifying its image on the generators $(r, s)$ of $\rightComplex{\righty{T}}$. Recall, that $(r, s)$ is  a pair with entries 1) a resolution of $T$ and 2) a sign assignment for each component of $r$. Then
\begin{equation}\label{eqn:delta}
\begin{array}{cl}
\righty{\delta_{T}}(r, s) =  & I_{\partial(r,s)} \otimes d(r,s)    + \disp{\sum_{\gamma \in \bridges{r}}  B(\gamma)} \\
\ & \ \\
\ & + \disp{\sum_{\gamma \in \dec{r,s}} (-1)^{I(r,\cross{\gamma})}\righty{e_{C(\gamma)}} \otimes (r_{\gamma}, s_{\gamma})} + \disp{\sum_{C \in \circles{\partial(r,s)}, s(C) = +} (-1)^{h(r)} \lefty{e_{C}} \otimes (r,s_{C}) }
\end{array}
\end{equation}
where
\begin{enumerate}

\item if $\gamma \in \rightBridges{r}$ then $\gamma$ corresponds to a crossing $c \in \cross{\righty{T}}$ and $r_{\gamma} = \gamma_{c}(r)$. If $\gamma$ merges two circles then $s_{\gamma}$ is computed from the Khovanov Frobenius algebra as for circles in $\free{r}$ and $$B(\gamma) =(-1)^{I(r,c)}e_{(\gamma, \partial s, (\partial s)_{\gamma})} \otimes (r_{\gamma}, s_{\gamma})$$ When surgery on $\gamma$ divides a circle $C$ in $\mathrm{cl}(r)$ with $s(C) = +$ we have two terms $$B(\gamma) =(-1)^{I(r,c)} \big(e_{(\gamma, s, s^{1}_{\gamma})}\otimes(\gamma_{c}(r), s^{1}_{\gamma})   +  e_{(\gamma, s, s^{2}_{\gamma})}\otimes(\gamma_{c}(r), s^{2}_{\gamma})\big) $$ When $s(C) = -$ there is just one term for $s^{-}_{\gamma}$ assigning $-$ to both the resulting circles and $$B(\gamma)=(-1)^{I(r,c)}e_{(\gamma, \partial s, (\partial s)_{\gamma})} \otimes (r_{\gamma}, s^{-}_{\gamma})$$ 

\item if $r = (\rho, \lefty{m})$ and $\gamma \in \leftBridges{\lefty{m}}$ then we use the same formulation to adjust
the decoration as in $\cleaved{n}$. In particular $\gamma$ corresponds to an edge(s) $\partial(r,s) =$
$(L,\sigma) \rightarrow (L_{\gamma}, \ast)$. To compute $r_{\gamma}$ we take $r_{\gamma} = \lefty{m}_{\gamma}\# \rho(\righty{T})$, then
$$
B(\gamma) = (-1)^{h(r)} e_{\gamma} \otimes (r_{\gamma}, s_{\gamma})
$$ 
is based on the morphism corresponding to the edge equipped with a sign from the homological grading. 

\item For $\gamma \in \dec{r,s}$ we have two cases both of which involve only one $C \in \mathrm{cl}(r)$. This circle is the $C(\gamma)$ in the formula. As a reminder the two cases are

\begin{enumerate}
\item both feet are on the same component of $C \cap \rightHalf$ for some $C \in \mathrm{cl}(r)$ with $s(C) = +$, 
 or
\item one foot of $\gamma_{c}$ is on $C \in \mathrm{cl}(r)$ with $s(C) = +$ and the other foot is on $C' \in \free{r}$ with $s(C') = -$
\end{enumerate} 

For these $r_{\gamma} = \gamma_{c}(r)$ as before. For $s_{\gamma}$ we take $s_{\gamma}(C(\gamma)) = -$ in both cases. $s_{\gamma}$ is $s$ on the uninvolved circles. In the case where both feet are on $C(\gamma)$ there is a new circle $D$ in $\gamma_{c}(r)$ and we require $s_{\gamma}(D) = +$. 
\end{enumerate}
\ \\
\noindent Before analyzing the properties of $\righty{\delta_{T}}$, we note that Figures \ref{fig:rightdec}, \ref{fig:rightdec2} illustrate the types of bridges used in the second to last summand, while \ref{fig:leftdec} illustrates the occurence of terms in the last summand.

\begin{prop}
$\righty{\delta_{T}}$ is grading preserving.
\end{prop}

\noindent{\bf Proof:} We show that each of the terms in the definition of $\righty{\delta_{T}}$ occurs in grading $(1,0)$ more than $(r,s)$ when considered in $\rightComplex{\righty{T}}$. First, in $(r,s) \rightarrow \I_{\partial(r,s)} \otimes d(r,s)$ we know that $d$ changes $(h,q)$ to $(h+1,q)$. For $(r,s) \rightarrow \lefty{e_{C}} \otimes (r,s_{C})$ when $s(C) = +$ we have $(h,\overline{q} + 1/2) \rightarrow (1,1) + (h, \overline{q}-1/2) = (h+1, \overline{q} + 1/2)$ since $C$ contribute $\pm 1/2$ to the quantum grading. For bridge edges, we distinguish merges and divides. For a right merge we have either $(h,\overline{q} + 1/2 + 1/2) \rightarrow (0,-1/2) + (h, \overline{q}+1/2) + (1,1)$ $= (h+1, \overline{q} + 1)$
or $(h,\overline{q} + 1/2 - 1/2) \rightarrow (0,-1/2) + (h, \overline{q}-1/2) + (1,1)$, where the first term in the sum comes from the bigrading of $\righty{e_{\gamma}}$, the second term comes from merging two $+$ circles or a $+$ and a $-$ circle (counted by $\pm 1/2$'s) and the last term comes from the homological/quantum shift occurring with a change in resolution at a crossing. For a right division we have either a division of a $+$ circle: $(h,\overline{q} + 1/2) \rightarrow (0,-1/2) + (h, \overline{q}+1/2-1/2) + (1,1) = (h+1,\overline{q} + 1/2)$, or a division of a $-$ circle: $(h,\overline{q} - 1/2) \rightarrow (0,-1/2) + (h, \overline{q}-1/2-1/2) + (1,1) = (h+1,\overline{q} - 1/2)$. For left merges there is no change in resolution at a crossing, and we ignore the decoration on the cleaved circles when calculating the bigrading, so we have either $(h,\overline{q} + 1/2 + 1/2) \rightarrow (1,1/2) + (h, \overline{q}+1/2)= (h+1,\overline{q} + 1)$ or $(h,\overline{q} + 1/2 - 1/2) \rightarrow (1,1/2) + (h, \overline{q}-1/2)= (h+1,\overline{q})$. For left divisions, we have
$(h,\overline{q} + 1/2) \rightarrow (1,1/2) + (h, \overline{q}+1/2-1/2)= (h+1,\overline{q}+1/2)$ or $(h,\overline{q} - 1/2) \rightarrow (1,1/2) + (h, \overline{q}-1/2-1/2)= (h+1,\overline{q}-1/2)$. This leaves terms of the form $(r,s) \rightarrow 
\righty{e_{C(\gamma)}} \otimes (r_{\gamma}, s_{\gamma})$ for $\gamma \in \dec{r,s}$. Such a term comes from changing the $+$ on a cleaved circle to a $-$ while consuming/creating a free circle, as a crossing resolution changes. For consumption we have use a $-$ free circle: $(h,\overline{q} + 1/2 - 1) \rightarrow (0,-1) + (h,\overline{q}-1/2) + (1,1) = (h+1,\overline{q}-1/2)$, while for creation we get a new $+$ free circle: $(h,\overline{q} + 1/2) \rightarrow (0,-1) + (h,\overline{q}-1/2+1) + (1,1) = (h+1,\overline{q}+1/2)$. In all cases, the bigrading changes by $(1,0)$. 
$\Diamond$


\begin{thm}
The map $\righty{\delta_{T}} : \complex{\righty{T}} \longrightarrow  \mathcal{B}\bridgeGraph{n}\otimes_{\mathcal{I}}\complex{\righty{T}}$ satisfies the relation
$$
(\mu_{\mathcal{B}\bridgeGraph{n}} \otimes \I)\,(\I \otimes \righty{\delta_{T}}) \, \righty{\delta_{T}}  +  (d_{\Gamma_{n}} \otimes |\I|)\,\righty{\delta_{T}} = 0
$$
where $\mu_{\mathcal{B}\bridgeGraph{n}}: \mathcal{B}\bridgeGraph{n} \otimes \mathcal{B}\bridgeGraph{n} \rightarrow \mathcal{B}\bridgeGraph{n}$ is the multiplication map on $\mathcal{B}\bridgeGraph{n}$.
\end{thm}

\noindent{\bf Proof:} Let $(r,s)$ be a generator for $\complex{\righty{T}}$ and
$$
\Xi = (\mu_{\mathcal{B}\bridgeGraph{n}} \otimes \I)\,(\I \otimes \righty{\delta_{T}}) \, \righty{\delta_{T}}(r,s)
$$
We start by computing the terms in $\Xi$.\\
\ \\
\noindent First $\Xi$ is a sum of terms of the form $A_{(r',s')} \otimes (r',s')$. We wish to show that $A_{(r',s')} = 0$ for every
pair $(r',s')$ once we include the term $(d_{\Gamma_{n}} \otimes \I)\,\righty{\delta_{T}}(\xi)$. Let $(L,\sigma) = \partial(r,s)$ and $(L',\sigma') = \partial(r',s')$. The coefficient $A_{(r',s')}$ is the sum
of signed products $e_{\alpha}e_{\beta}$ taken over all paths
$$ 
\begin{CD} 
(r,s) @>A>> (\overline{r},\overline{s}) @> B >> (r',s') \\
@VV\partial V @VV\partial V @VV\partial V\\
(L,\sigma) @>e_{\alpha}>> (\overline{L},\overline{\sigma}) @> e_{\beta} >> (L',\sigma') 
\end{CD}
$$ 
where $A$ and $B$ are either surgery on an active arc in their source resolution, or a bridge in $\leftBridges{r}$ or $\leftBridges{\overline{r}}$, or correspond to a cleaved circle decorated by a $+$. $e_{\alpha}$ and $e_{\beta}$ are the corresponding edges in $\bridgeGraph{n}$ coming from the effect on $\partial(r,s)$.\\
\ \\
\noindent As such we consider the possible paths between $(r,s)$ and $(r',s')$.  Consider the map $J: \xi = (r, s) \rightarrow (\#\free{r},\#\circles{\mathrm{cl}(r)}, \iota(\partial \xi))$. Terms contributing to $\I \otimes d$ change $J$ by $(\pm 1, 0, 0)$ depending on whether a merge or a division occurs. Terms coming from summing over bridges change the value of $J$ by $(0, \pm 1, -1)$. Terms coming from $\gamma \in \dec{r,s}$ change $J$ by $(\pm 1, 0, -2)$ while those coming from a circle with $s(C) = +$ change $J$ by $(0,0,-2)$.     \\
\ \\
\noindent We run through the possibilities for the two paths $A$ and $B$: 
\begin{center}
\begin{tabular}{|c||c|c|c|c|}
\hline
\ &                           $I \otimes d$ &       $\bridges{\xi}$ &         $\dec{\xi}$ &          $\circles{\mathrm{cl}(r)}$ \\
\hline
\hline
$I \otimes d$ &               $(\{\pm 2,0\}, 0,0)$ & $(\pm 1, \pm 1, -1)$ & $(\{\pm 2,0\}, 0, -2)$ & $(\pm 1, 0, -2)$ \\    
\hline
$\bridges{\xi}$ &              \ &                  $(0, \{\pm 2, 0\}, -2)$ &   $(\pm 1, \pm 1, -3)$ & $(0, \pm 1, -3)$ \\
\hline
$\dec{\xi}$ &                 \ &                   \                      &   $(\{\pm 2, 0\}, 0, -4)$ & $(\pm 1, 0, -4)$ \\
\hline
$\circles{\mathrm{cl}(r)}$&   \  &                   \ &                        \ &                       $(0,0,-4)$\\
\hline
\end{tabular}
\end{center}

\noindent Any path which goes from $(r,s)$ to $(r',s')$ will need to change $J$ in the same manner. For most of the changes above, the only way to achieve the change is for the types of $A$ and $B$ to be the same (or reversed). Thus if one is  a bridge, and the other is in $\dec{\xi}$ then all other paths will also need to have this form. There are two overlapping sets of entries for $(0,0,-4)$ and $(0,0,-2)$. We will return to the $(0,0,-2)$ case shortly.\\ 
\ \\
\noindent First, notice that two alterations to $r$ as in the upper left $3 \times 3$ subtable will change $r$ by surgering two active arcs. This changes the diagram to $\rightHalf$ in two distinct places. On the other hand, the entries in the far right column (and symmetrically in the last row) will surger at most one, and possibly no, active arcs. Consequently, the terms coming
from the right column will sum over $r'$'s that cannot come from the upper $3 \times 3$ subtable. In particular, the $(0,0,-4)$
repetition in the table will not result in canceling terms. We start therefore, by understanding the sums when either $A$ or $B$ comes from changing the decoration on a circle in $\mathrm{cl}(r)$ or $\mathrm{cl}(\overline{r})$ from $+$ to $-$ using a term from $\circles{\mathrm{cl}(r)}$. \\
\ \\
\noindent Suppose $(r',s')$ comes from two terms from $\circles{\xi}$ as in the lower right entry. Then $r'=r$ and $s'= s_{C_{1}, C_{2}}$ for distinct cleaved circles $C_{1}$ and $C_{2}$ with $s(C_{i}) = +$, $i = 1,2$. The table above shows that the only other pairs that can lead to $(r',s')$ from $(r,s)$ also come from two circles. Due to the form of $s'$ these must be the circles $C_{1}$ and $C_{2}$. For terms coming from $\leftHalf$, there is a $1-1$ identification from the changes in state to an edge in $\bridgeGraph{n}$. There are thus only two possible paths in $\bridgeGraph{n}$ corresponding to $(r,s) \rightarrow (r',s')$:

$$ 
\begin{CD} 
(L,\sigma) @>\lefty{e_{C_{1}}}>> (\overline{L},\sigma_{C_{1}}) @> \lefty{e_{C_{2}}} >> (L,\sigma_{C_{1},C_{2}}) 
\end{CD}
$$ 
and
$$ 
\begin{CD} 
(L,\sigma) @>\lefty{e_{C_{2}}}>> (\overline{L},\sigma_{C_{2}}) @> \lefty{e_{C_{1}}} >> (L,\sigma_{C_{2},C_{1}}) 
\end{CD}
$$ 
These yield $A_{(r',s')} = \big((-1)^{h(r)}\lefty{e_{C_{1}}}(-1)^{h(\overline{r})}\lefty{e_{C_{2}}} + (-1)^{h(r)}\lefty{e_{C_{2}}}(-1)^{h(\overline{r})}\lefty{e_{C_{1}}}\big)$ $= \big(\lefty{e_{C_{1}}}\lefty{e_{C_{2}}} + \lefty{e_{C_{2}}}\lefty{e_{C_{1}}}\big)$ since changing the decoration on a circles does not change $r$, i.e. $r=r'=\overline{r}$. The elements $\lefty{e_{C}}$ are odd in $\mathcal{B}\Gamma_{n}$, so they satisfy anti-commutative relations which imply $A_{(r',s')} = 0$.\\
\ \\ 
\noindent Essentially the same argument applies when $A$ and $B$ are a pair coming from a cleaved circle $C$ with $s(C) = +$ and an active arc $\gamma \in \dec{r,s}$. The active arc changes the decoration on a circle $D \in \mathrm{cl}(\overline{r})$, but $D \neq C$ since $C$ will have a $-$ decoration. Consequently, $D \in \mathrm{cl}(r)$ with $s(D) = +$. We can read $C$ and $D$ from $\partial(r',s')$ and
$\gamma$ from the resolution diagram of $r'$. Thus the only other way to go from $(r,s)$ to $(r',s')$ is to use $\gamma \in \dec{r,s}$ first, and then follow with changing the decoration on $C$. This yields two paths 
$$ 
\begin{CD} 
(L,\sigma) @>\lefty{e_{C}}>> (L,\sigma_{C}) @> \righty{e_{D}} >> (L,\sigma_{C,D}) 
\end{CD}
$$ 
and
$$ 
\begin{CD} 
(L,\sigma) @>\righty{e_{D}}>> (L,\sigma_{C_{2}}) @> \lefty{e_{C}} >> (L,\sigma_{D,C}) 
\end{CD}
$$ 
Then $A_{(r',s')} = (-1)^{I(r,c)}\righty{e_{D}}(-1)^{h(r_{\gamma})}\lefty{e_{C}} +  (-1)^{h(r)}\lefty{e_{C}}(-1)^{I(r,c)}\righty{e_{D}} = (-1)^{I(r,c)+h(r)}\big(-\righty{e_{D}}\lefty{e_{C}} + \lefty{e_{C}}\righty{e_{D}}\big) = 0 $ with the cancellation coming from  the (commuting) relations for disjoint supports. Note also that since $\gamma \in \actor{r}$ that $r_{\gamma}$ has homological grading $1$ higher than $r$.\\
\ \\
\noindent If $\gamma$ instead contributes to $d$, then there will be a term in $\I \otimes d$ from it. However, $\gamma$ does not change $\partial(r,s)$ in this case. Furthermore, the contribution is the same regardless of the decoration $\sigma = \partial(s)$. In other words, from $(r',s')$ we can read off $C$ since it is the only circle which will change decoration, and we can read off $\gamma$ due to the change in the resolution diagram. These commute due to the definition of $d$ not distinguishing the decoration on a cleaved circle. There are thus two possible paths

$$ 
\begin{CD} 
(L,\sigma) @>\I_{(L,\sigma)}>> (L,\sigma) @> \lefty{e_{C}} >> (L,\sigma_{C}) 
\end{CD}
$$ 
and
$$ 
\begin{CD} 
(L,\sigma) @>\lefty{e_{C}}>> (L,\sigma_{C}) @> \I_{(L,\sigma_{C})} >> (L,\sigma_{C}) 
\end{CD}
$$ 
which make $A_{(r',s')} = (-1)^{I(r,\cross{c})}(-1)^{h(r_{\gamma})} \lefty{e_{C}} + (-1)^{h(r)}(-1)^{I(r,\cross{c})}\lefty{e_{C}} = (-1)^{I(r,\cross{c})+h(r)} \big( -\lefty{e_{C}} + \lefty{e_{C}}\big) = 0$ in this case. \\
\ \\
\noindent Now suppose we have $(r,s)$ with $C \in \mathrm{cl}(r)$ and $s(C) = +$, and $\gamma \in \bridges{r}$ with
an arrow $(L,\sigma) \rightarrow (L_{\gamma}, \sigma_{\gamma})$. There are two cases to consider: 1) when $C$ is not in the support of $\gamma$, and 2) when $\gamma$ has at least one foot on $C$. Suppose $C$ is not in the support of $\gamma$. 
Looking at $(r',s')$ we will see a cleaved circle whose decoration has changed from $+$ to $-$ and surgery along an arc that is also a bridge. This identifies $\gamma$. Furthermore, each does not affect the decorations of the support of the other. Consequently, there are only two paths from $(r,s)$ to $(r',s')$ which give rise to coefficients from
$$ 
\begin{CD} 
(L,\sigma) @>\lefty{e_{C}}>> (L,\sigma_{C}) @> e_{\gamma} >> (L_{\gamma},\sigma_{C,\gamma}) 
\end{CD}
$$ 
and
$$ 
\begin{CD} 
(L,\sigma) @>e_{\gamma}>> (L,\sigma_{\gamma}) @> \lefty{e_{C}} >> (L,\sigma_{C, \gamma}) 
\end{CD}
$$ 
When $\gamma \in \leftBridges{r}$, these paths yield $A_{(r',s')} =(-1)^{h(r)} \lefty{e_{C}}(-1)^{h(r)}e_{\gamma} + (-1)^{h(r)} e_{\gamma}(-1)^{h(r)}\lefty{e_{C}} = 0$  due to the anti-commutative relations for disjoint supports. When $\gamma \in \rightBridges{r}$ these paths yield  $A_{(r',s')} = (-1)^{h(r)}\lefty{e_{C}}(-1)^{I(r,\cross{\gamma})}e_{\gamma} + (-1)^{I(r,\cross{\gamma})} e_{\gamma}(-1)^{h(r_{\gamma})}\lefty{e_{C}} =$ $(-1)^{I(r,\cross{\gamma})+h(r)}\big(\lefty{e_{C}}e_{\gamma} - e_{\gamma}\lefty{e_{C}}\big) = 0$.\\
\ \\\ 
\noindent We only describe the case of $C$ being in the support of $\gamma$ in the case that $\gamma \in \merge{L}$. In particular, there are two cleaved circles $D_{1}$ and $D_{2}$ joined by $\gamma$ to form a circle $D$. $C$ will be in the support of $\gamma$ only if $C \in \{D_{1}, D_{2}, D\}$. Let $C = D_{1}$, then there is a path
$$ 
\begin{CD} 
(L,\sigma) @>\lefty{e_{D_{1}}}>> (L,\sigma_{D_{1}}) @> e_{\gamma} >> (L_{\gamma},\sigma'') 
\end{CD}
$$ 
only when $\sigma(D_{1}) = \sigma(D_{2}) = +$. In that case $\sigma'' = \sigma_{\gamma, D}$. Looking at $(r',s')$ we can
determine the bridge $\gamma$ from its resolution diagram and see that the decoration changed on $D$, but this can occur either as above, or from
$$ 
\begin{CD} 
(L,\sigma) @>\lefty{e_{D_{2}}}>> (L,\sigma_{D_{2}}) @> e_{\gamma} >> (L_{\gamma},\sigma_{\gamma, D}) 
\end{CD}
$$ 
or from
$$ 
\begin{CD} 
(L,\sigma) @> e_{\gamma}>> (L_{\gamma},\sigma_{\gamma}) @> \lefty{e_{D}} >> (L_{\gamma},\sigma_{\gamma, D}) 
\end{CD}
$$ 
since this is a merge, the decorations $\sigma_{\gamma}$ is uniquely determined. There are, however, no other possibilities
and any other decorations for $D_{1}$ and $D_{2}$ cannot produce paths of this type. If $\gamma \in \leftBridges{r}$ then
$A_{(r',s')} = \lefty{e_{D_{1}}} e_{\gamma} + \lefty{e_{D_{2}}}  e_{\gamma} + e_{\gamma} \lefty{e_{D_{C}}} = 0$ (each element also comes with a sign $(-1)^{h(r)}$, but these are the same for each element). On the other hand if $\gamma \in \rightBridges{r}$ then  $A_{(r',s')} = (-1)^{I(r,\cross{\gamma})+h(r)}\lefty{e_{D_{1}}} e_{\gamma} + (-1)^{I(r,\cross{\gamma})+h(r)} \lefty{e_{D_{2}}}  e_{\gamma} + (-1)^{I(r,\cross{\gamma})+h(r_{\gamma}}e_{\gamma} \lefty{e_{D_{C}}}$ $= (-1)^{I(r,\cross{\gamma})+h(r)}\big(\lefty{e_{D_{1}}} e_{\gamma} + \lefty{e_{D_{2}}}  e_{\gamma} - e_{\gamma} \lefty{e_{D_{C}}}\big) = 0$ by relation \ref{rel:lefty1}.\\
\ \\
\noindent The case when $\gamma$ is a division is similar.\\
\ \\
\noindent Thus the coefficient is zero on any $(r',s')$ which comes from a path corresponding to the right column. We turn
our attention to those paths which correspond to a pair $(\gamma_{1}, \gamma_{2})$, in that order, where each $\gamma_{i}$ is
either an active arc for $r$ or a bridge in $\leftBridges{r}$. \\
\ \\
\noindent For a pair of active arcs $(\gamma_{1}, \gamma_{2})$ giving a path from the top left corner, the coefficient is
$$ 
\begin{CD} 
(L,\sigma) @> \I_{(L,\sigma)}>> (L,\sigma) @> \I_{(L,\sigma)} >> (L_,\sigma) 
\end{CD}
$$
or just $\I_{(L,\sigma)}$. We know that such a pair can be switched to give another term in $d^{2}$ corresponding to $(\gamma_{2}, \gamma_{1})$ and that these pairs cancel (since $d$ is a differential). The table above shows that none of the new terms can be used to produce the same $(r',s')$, so we have verified that the coefficient is still $0$ for any such $(r',s')$. Note that the signs in this case follow from $I(r,\cross{\gamma_{1}}) + I(r_{\gamma_{1}}, \cross{\gamma_{2}})$ being $1$ larger or smaller than $I(r,\cross{\gamma_{2}}) + I(r_{\gamma_{2}}, \cross{\gamma_{1}})$, as in \cite{Bar1}.\\
\ \\
\noindent Now suppose $\gamma_{1}$ corresponds to a term in $\I \otimes d$, and $\gamma_{2} \in \bridges{r}$. By looking at
$(L',\sigma') = \partial(r',s')$ we can determine the type of bridge, and by looking at $r'$ we can determine the active arcs. Thus, the only paths to get from $(r,s)$ to $(r',s')$ in this case must come from $(\gamma_{1}, \gamma_{2})$ and $(\gamma_{2}, \gamma_{1})$. The table above ensures that each consists of an arc contributing to $\I \otimes d$ and one to the term coming from $\bridges{r}$. Since surgery on $\gamma_{1}$ doesn't alter $\partial(r,s)$ while surgery on $\gamma_{2}$ merges or divides cleaved circles, in $(\gamma_{2}, \gamma_{1})$ $\gamma_{2}$ still merges or divides cleaved circles, and $\gamma_{1}$ still contributes to $I \otimes d$. In particular, the decoration changes of one arc does not affect the other. Thus there are two possible paths contributing to the coefficient of $(r',s')$
$$ 
\begin{CD} 
(L,\sigma) @> \I_{(L,\sigma)}>> (L,\sigma) @> e_{\gamma_{2}} >> (L_{\gamma},\sigma') 
\end{CD}
$$
and
$$ 
\begin{CD} 
(L,\sigma) @> e_{\gamma_{2}}>> (L_{\gamma},\sigma') @> \I_{(L,\sigma)} >> (L_{\gamma},\sigma') 
\end{CD}
$$
which yields a coefficient of $(-1)^{I(r,\gamma_{1}) + I(r_{\gamma_{1}}, \gamma_{2})}e_{\gamma_{2}} + (-1)^{I(r,\gamma_{2}) + I(r_{\gamma_{2}}, \gamma_{1})}e_{\gamma} = 0$ when $\gamma_{2} \in \rightBridges{r}$. When $\gamma_{2} \in \leftBridges{r}$
the coefficient is $(-1)^{I(r,\gamma_{1})+h(r_{\gamma_{1}}}I \cdot e_{\gamma_{2}} + (-1)^{I(r_{\gamma_{2}}, \gamma_{1})+h(r_{\gamma_{2}})}e_{\gamma_{2}} \cdot I$. However  $I(r_{\gamma_{2}}, \gamma_{1}) = I(r, \gamma_{1})$ and $h(r_{\gamma_{2}}) = h(r)$ since only the $1$-resolutions along active arcs contribute to these sums. Consequently, the coefficient is $(-1)^{I(r,\gamma_{1})+h(r)}\big( -e_{\gamma_{2}} + e_{\gamma_{2}} \big) = 0$. \\
\ \\
\noindent Now suppose we have a sequence $\gamma_{1}$ and $\gamma_{2}$ where $\gamma_{1}$ is in $\dec{r,s}$ and $\gamma_{2} \in \dec{r_{\gamma_{1}}, s_{\gamma_{1}}}$. There is a cleaved circle $C_{1}$ in $r$ with $s(C_{1}) = +$, while $s_{\gamma_{1}}(C'_{1}) = -$ where $C_{1}'$ runs through the same points in $P_{n}$.  Consequently, $\gamma_{2}$ must have support on a {\em different} cleaved circle $C_{2}$ (it is cleaved by assumption). Thus $\gamma_{2} \in \dec{r,s}$ and $\gamma_{1} \in \dec{r_{\gamma_{2}}, s_{\gamma_{2}}}$. The resolution $r'$ identifies the two active arcs used in this process, while the change in decoration identifies that both arcs must be used in $\dec{r,s}$ or the equivalent. There can be at most two paths from $(r,s)$ to $(r',s')$. and these contribute a coefficient $(-1)^{I(r,\gamma_{1}) + I(r_{\gamma_{1}}, \gamma_{2})}\righty{e_{C_{1}}}\righty{e_{C_{2}}} + (-1)^{I(r,\gamma_{2}) + I(r_{\gamma_{2}}, \gamma_{1})}\righty{e_{C_{2}}}\righty{e_{C_{1}}}$ which is $0$ by the relations for disjoint supports.  \\
\ \\
\noindent  Suppose $\gamma_{1}$ and $\gamma_{2}$ are active arcs in $r$ with $\gamma_{1}$ contributing to $\I \otimes d$ and $\gamma_{2}$ contributing to $\dec{r_{\gamma_{1}}, s_{\gamma_{1}}}$, and $\gamma_{2}$ abutting only one cleaved circle $C$.   Since $\gamma_{2}$ is in $\dec{r_{\gamma_{1}}, s_{\gamma_{1}}}$, surgery on it changes the decoration on  $C$ from $+$ to $-$, but otherwise leaves $\mathrm{cl}(r)$ unchanged. Since we surgered two active arcs to get $(r',s')$, we can determined that $(r',s')$ results from surgery on two distinct active arcs, changing the decoration (and shape) of one cleaved circle. Furthermore, neither arc is a bridge in $r$. This is obvious for $\gamma_{1}$, but for $\gamma_{2}$ it follows from only abutting one cleaved circle in $r$. Since it is in $\dec{r_{\gamma_{1}}}$ it still only abuts one cleaved circle there, and that circle must be $C$. When switching $\gamma_{1}$ and $\gamma_{2}$, we must again have one in $\dec{}$ which changes the decoration on $C$, and the other in $\I \otimes d$. Consequently, the coefficient on $(r',s')$ will be $(-1)^{I(r,\gamma_{1}) + I(r_{\gamma_{1}}, \gamma_{2})}\righty{e_{C}} + (-1)^{I(r,\gamma_{2}) + I(r_{\gamma_{2}}, \gamma_{1})}\righty{e_{C}} = 0$. Note that which is in $\dec{}$ is determined by the decorations on $(r',s')$ and $(r,s)$. In the case where both divide $C$, the roles can switch. \\  
\ \\
\noindent If $\gamma_{1}$ is a bridge for $r$ and $\gamma_{2} \in \dec{r_{\gamma_{1}}, s_{\gamma_{1}}}$ then $(r',s')$ identifies that $\gamma_{1}$ is a bridge and there is another active arc for any path from $(r,s)$ to $(r',s')$. By comparing the decorations in $\partial(r',s')$ to those in $\partial(r_{\gamma_{1}}, s_{\gamma_{1}})$ we can see that $\gamma_{2}$ is in $\dec{r,s}$ as well. If these have disjoint supports then the triviality of the coefficient on $(r',s')$ follows as before from the relations for disjoint supports. If not, we use the relations for $\righty{e_{C}}$. For example, suppose $\gamma_{1}$ merges $C_{1}$ and $C_{2}$ and $\gamma_{2}$ has active circle $C_{2}$. For the path to exist, $C_{1} \#_{\gamma_{1}} C_{2}$ must have $+$  decoration, so $C_{1}$ and $C_{2}$ will as well. Thus we can surger $\gamma_{2}$ first in $\dec{r,s}$ and then surger $\gamma_{1}$. This gives two paths
$$ 
\begin{CD} 
(L,\sigma) @> m_{\gamma_{1}}>> (L_{\gamma_{1}},\sigma_{\gamma_{1}}) @> \righty{e_{C}} >> (L_{\gamma_{1}},\sigma_{\gamma_{1},C}) 
\end{CD}
$$
and
$$ 
\begin{CD} 
(L,\sigma) @> \righty{e_{C_{2}}}>> (L,\sigma_{C_{2}}) @> m_{\gamma_{1}} >> (L_{\gamma_{1}},\sigma_{\gamma_{1}, C}) 
\end{CD}
$$
When $\gamma_{1} \in \rightBridges{r}$ the resulting coefficient is $(-1)^{I(r,\gamma_{1}) + I(r_{\gamma_{1}}, \gamma_{2})}m_{\gamma_{1}}\righty{e_{C}} + (-1)^{I(r,\gamma_{2}) + I(r_{\gamma_{2}}, \gamma_{1})}\righty{e_{C_{2}}}m_{\gamma_{1}} = (-1)^{I(r,\gamma_{1}) + I(r_{\gamma_{1}}, \gamma_{2})}\big(m_{\gamma_{1}}\righty{e_{C}} - \righty{e_{C_{2}}}m_{\gamma_{1}}\big)=0$. Note that where $\gamma_{2}$ occurs, namely on $C_{2}$, determines which relation we use. For $\gamma_{2}$ on $C_{1}$ we would use $\righty{e_{C_{1}}}$, but we would also have a different resolution diagram for $r'$.\\
\ \\
\noindent On the other hand if $\gamma_{1} \in \leftBridges{r}$ then the resulting coefficient is $(-1)^{I(r_{\gamma_{1}}, \gamma_{2})+h(r)} m_{\gamma_{1}}\righty{e_{C}} + (-1)^{I(r, \gamma_{2}) +h(r_{\gamma_{2}})} \righty{e_{C_{2}}}m_{\gamma_{1}}$. But $I(r_{\gamma_{1}}, \gamma_{2}) = I(r, \gamma_{2})$ so
$A_{(r',s')} = (-1)^{I(r, \gamma_{2})+h(r)} \big(m_{\gamma_{1}}\righty{e_{C}} - \righty{e_{C_{2}}}m_{\gamma_{1}}\big) = 0$.\\
\ \\
\noindent If $\gamma_{1}$ and $\gamma_{2}$ are both bridges for $(r,s)$, then $(r',s')$ will have a resolution which will differ from $(r,s)$ by two local surgeries. When $\gamma_{2} \in B_{d}(\partial(r,s), \gamma_{1})$, that is when they have feet on different arcs intersecting the $y$-axis, or when one is in $\actor{r}$, the only other path to reach $(r',s')$ will be to surger $\gamma_{2}$ first, and then surger $\gamma_{1}$. Each of these surgeries will occur along a bridge, since surgeries on bridges do not produce free circles. 
The two paths produce coefficients on $(r',s')$ from 
$$ 
\begin{CD} 
(L,\sigma) @> e_{\gamma_{1}}>> (L_{\gamma_{1}},\sigma_{\gamma_{1}}) @> e_{\gamma_{2}} >> (L_{\gamma_{1}, \gamma_{2}},\sigma_{\gamma_{1},\gamma_{2}}) 
\end{CD}
$$
and
$$ 
\begin{CD} 
(L,\sigma) @> e_{\gamma_{2}}>> (L,\sigma_{\gamma_{2}}) @> e_{\gamma_{1}} >> (L_{\gamma_{1}, \gamma_{2}},\sigma_{\gamma_{1}, \gamma{2}}) 
\end{CD}
$$
Most of the possibilities have a unique pairing between two such paths (see the discussion of decorations in the bridge relations for section), and thus a coefficient $A_{(r',s')}$ given by
\begin{enumerate}
\item $(-1)^{I(r,\gamma_{1}) + I(r_{\gamma_{1}}, \gamma_{2})} e_{\gamma_{1}}e_{\gamma_{2}} + (-1)^{I(r,\gamma_{2}) + I(r_{\gamma_{2}}, \gamma_{1})} e_{\gamma_{2}}e_{\gamma_{1}} =0$ if both $\gamma_{1}$ and $\gamma_{2}$ are in $\actor{r}$.
\item $(-1)^{I(r,\gamma_{1})+h(r_{\gamma_{1}})} e_{\gamma_{1}}e_{\gamma_{2}} + (-1)^{h(r) + I(r_{\gamma_{2}}, \gamma_{1})} e_{\gamma_{2}}e_{\gamma_{1}} =$ $(-1)^{I(r,\gamma_{1})+h(r)}\big(-e_{\gamma_{1}}e_{\gamma_{2}} + e_{\gamma_{2}}e_{\gamma_{1}} \big)= 0$ when $\gamma_{1}$ is in $\actor{r}$ and $\gamma_{2} \in \leftBridges{r}$.
\item $(-1)^{h(r) + h(r)} \big(e_{\gamma_{1}}e_{\gamma_{2}} + e_{\gamma_{2}}e_{\gamma_{1}}\big) = 0$ when both $\gamma_{1}$ and $\gamma_{2} \in \leftBridges{r}$.  
\end{enumerate}
If $\gamma_{1}$ and $\gamma_{2}$ both divide the same cleaved circle $C$ then there are four possible paths to $(r',\sigma_{C})$. However, the bridge relations assert that they cancel or commute in pairs, hence their sum will also be $0$.\\
\ \\
\noindent When $\gamma_{1}, \gamma_{2}$ are bridges for $(r,s)$, and the preceding paragraph does not apply, then $\gamma_{1}, \gamma_{2} \in \leftBridges{r}$ and have endpoints on a shared arc. If $\gamma_{2} \in B_{o}(\partial(r,s),\gamma_{1})$,  there are two bridges $\gamma_{2}^{+}, \gamma_{2}^{-}$ for $r_{\gamma_{1}}$ which map to $\gamma_{2}$ under surgery along $\gamma_{1}^{\dagger}$. If we orient the arc containing the shared endpoint of $\gamma_{1}$ and $\gamma_{2}$, these correspond to representatives of $\gamma_{2}$ whose endpoints come either before or after that of $\gamma_{2}$. Likewise, there are two bridges $\gamma_{1}^{+}$ and $\gamma_{1}^{-}$ in $r_{\gamma_{2}}$. If we follows the paths for surgery on the pair $(\gamma_{1},\gamma_{2}^{+})$ to get $(r',s')$, then the only other way to get to $(r',s')$ is to surger $(\gamma_{2}, \gamma_{1}^{-})$. As these are left bridges, the corresponding coefficient for $(\mu \otimes \I)(\I \otimes \delta)\delta$ will be $(-1)^{h(r) + h(r)}(\lefty{e_{\gamma_{1}}}\lefty{e_{\gamma_{2}^{+}}} + \lefty{e_{\gamma_{2}}}\lefty{e_{\gamma_{1}^{-}}}) = 0$, since these elements anti-commute. On the other hand, surgery on $(\gamma_{1},\gamma_{2}^{-})$ yields the same resolution as surgery on $(\gamma_{2}, \gamma_{1}^{+})$, and the same argument shows that the coefficient $(-1)^{h(r) + h(r)}(\lefty{e_{\gamma_{1}}}\lefty{e_{\gamma_{2}^{-}}} + \lefty{e_{\gamma_{2}}}\lefty{e_{\gamma_{1}^{+}}})$ is also $0$.\\
\ \\
\noindent When $\gamma_{1}, \gamma_{2} \in \leftBridges{r}$ and $\gamma_{2} \in B_{s}(\partial(r,s),\gamma_{1})$, let $\gamma_{3}$ be the bridge found by sliding one endpoint of $\gamma_{1}$ across $\gamma_{2}$. Then in $r_{\gamma_{1}}$, $\gamma_{2}=\gamma_{3}$ and we will call this class $\eta_{23}$. Likewise, let $\eta_{12}$ and $\eta_{13}$ be the bridges determined by $\gamma_{1}$ in $r_{\gamma_{2}}$ and $r_{\gamma_{3}}$, respectively. Surgery on $\gamma_{i}$ followed by $\eta_{jk}$ with $\{i,j,k\}=\{1,2,3\}$ will give the same $r'$ for each of the three possibilities. Consequently, the coefficient on $(r',s')$ will be
$(-1)^{h(r)+h(r_{\gamma_{1}})}\lefty{e_{\gamma_{1}}}\lefty{e_{\eta_{23}}} +   (-1)^{h(r)+h(r_{\gamma_{2}})}\lefty{e_{\gamma_{2}}}\lefty{e_{\eta_{13}}} + (-1)^{h(r)+h(r_{\gamma_{3}})}\lefty{e_{\gamma_{3}}}\lefty{e_{\eta_{12}}} $ $ = \lefty{e_{\gamma_{1}}}\lefty{e_{\eta_{23}}} +   \lefty{e_{\gamma_{2}}}\lefty{e_{\eta_{13}}} + \lefty{e_{\gamma_{3}}}\lefty{e_{\eta_{12}}} = 0
$
since surgery on $\gamma_{i}$ does not change the homological grading.\\
\ \\
\noindent However, there is another way for two bridges to arise: $\gamma_{1}$ can be a bridge, and $\gamma_{2}$ is a bridge in  $r_{\gamma_{1}}$. For this to be different from the previous case, $\gamma_{2}$ must either intersect  $\gamma_{1}^{\dagger}$ or be isotopic through bridges to $\gamma_{1}^{\dagger}$. When $\gamma_{1} \in \leftBridges{r}$ it is possible that such a $\gamma_{2}$ does intersect $\gamma_{1}^{\dagger}$. Then we have can have 
a path
$$
\begin{CD}
(L,s) @>e_{(\gamma_{1},\sigma,\sigma')}>> (L_{\gamma_{1}},\sigma') @>e_{(\gamma_{2},\sigma',\sigma'')}>> (L_{\gamma,\gamma_{2}}, \sigma'')\\
\end{CD}
$$
contributing to $A_{(r',s')}$. There is no other way to obtain $(r',s')$, so its sole coefficient is $\pm e_{(\gamma_{1},\sigma,\sigma')}e_{(\gamma_{2},\sigma',\sigma'')} = 0$ by the other bridge relations. 
 \\
\ \\
\noindent Now suppose $\gamma_{2}$ is isotopic through bridges to $\gamma_{1}^{\dagger}$. We will distinguish the cases where $\gamma_{1} \in \rightBridges{r}$ from $\gamma_{1} \in \leftBridges{r}$. When $\gamma_{1} \in \rightBridges{r}$, $\gamma_{2}$ must be a different active arc for $r$. If we let $\eta \in \rightBridges{\mathrm{cl}(r)}$ be the bridge corresponding to $\gamma_{1}$, then the coefficient from surgery on $\gamma_{1}$ followed by $\gamma_{2}$ will be
$$ 
\begin{CD} 
(L,\sigma) @> e_{\eta}>> (L,\sigma_{\eta}) @> e_{\eta^{\dagger}} >> (L,\sigma_{C}) 
\end{CD}
$$
which gives $(-1)^{I(r,\gamma_{1}) + I(r_{\gamma_{1}}, \gamma_{2})} e_{\eta}e_{\eta^{\dagger}}$ for some circle $C$ in $\mathrm{cl}(r)$ with $s(C) = +$. Reversing the order we see that either 1) $\gamma_{2}$ will need to be in $\dec{r,s}$ with surgery changing the decoration on $C$, while $\gamma_{1}$ contributes to $\I \otimes d$ on $r_{\gamma_{2}}$, or 2) $\gamma_{2}$ will need to contribute to $\I \otimes d$ with surgery , while $\gamma_{1}$ contributes to $\dec{r_{\gamma_{2}},s_{\gamma_{2}}}$ and changing the decoration on $C$. Which we obtain depends on whether $\gamma_{2}$ abuts the active circle for $\eta$ or not. The coefficient in either of these cases will be $(-1)^{I(r,\gamma_{2}) + I(r_{\gamma_{2}}, \gamma_{1})}\righty{e_{C}} \otimes (r',s')$. However, if $C$ is the active circle for $\eta$ we will have 
$e_{\eta}e_{\eta^{\dagger}} = \righty{e_{C}}$ from the relations, so $A_{(r',s')} = (-1)^{I(r,\gamma_{1}) + I(r_{\gamma_{1}}, \gamma_{2})}\big(e_{\eta}e_{\eta^{\dagger}} - \righty{e_{C}}\big) = 0$. (Note: if $\eta$ divides a $+$ cleaved circle, there are two choices of decoration, both of which merge to make the circle $-$ after surgery by $\eta^{\dagger}$. This contributes two terms $(-1)^{I(r,\gamma_{1}) + I(r_{\gamma_{1}},\gamma_{2})}\big(e_{(\eta,\sigma,\sigma^{1}_{\eta})}e_{(\eta^{\dagger}, \sigma^{1}_{\eta},\sigma_{C})} + e_{(\eta,\sigma,\sigma^{2}_{\eta})}e_{(\eta^{\dagger},\sigma^{2}_{\eta},\sigma_{C})}\big)$. On the other hand, surgery on $\eta^{\dagger}$ followed by $\eta$ will now produce two terms $(-1)^{I(r,\gamma_{2}) + I(r_{\gamma_{2}}, \gamma_{1})}\big(\I \otimes \righty{e_{C}}\big)$ and $(-1)^{I(r,\gamma_{2}) + I(r_{\gamma_{2}}, \gamma_{1})}\big(\righty{e_{C}} \otimes \I\big)$ which will cancel the preceding sum after applying $\mu \otimes \I$. \\
\ \\
\noindent For $\gamma_{1} \in \leftBridges{r}$, $\gamma_{2} \in \leftBridges{r_{\gamma_{1}}}$ there is a bridge $\eta$ with $\gamma_{1} = \eta$ and $\gamma_{2} = \eta^{\dagger}$. If $r = \lefty{m} \# \rho(\righty{T})$ then we can ignore everything not in $\mathrm{cl}(r)$ and concentrate solely on bridges and cleaved circles. We will thus obtain the sum 
$$
\sum_{\eta \in \leftBridges{\mathrm{cl}(r)}} e_{(\eta, \sigma, \sigma')} e_{(\eta^{\dagger}, \sigma', \sigma_{C(\eta)})} \otimes (r,s_{C(\eta)}) 
$$
where $C(\eta)$ is the active circle for $\eta$ and the signs from the homological grading cancel in each term. Breaking this up along circles $C \in \circles{\partial(r,s)}$ with $s(C) = +$ we find that the coefficient of $(r,s_{C})$ is the sum 
$$
\sum  e_{(\eta, \sigma, \sigma')} e_{(\eta^{\dagger}, \sigma', \sigma_{C})}
$$
taken over all paths
$$ 
\begin{CD} 
(L,\sigma) @>e_{(\eta,\sigma,\sigma')}>> (L_{\eta},\sigma') @>e_{(\eta^{\dagger},\sigma',\sigma_{C})}>> (L,\sigma_{C}) 
\end{CD}
$$ 
Thus the coefficient on $(r,s_{C})$ is exactly $-d_{\Gamma}(\lefty{e_{C}}) \otimes (r,s_{C})$. Thus adding
$$
\begin{array}{cl}
(d_{\mathcal{B}\bridgeGraph{n}} \otimes |\I|)\righty{\delta_{T}}(r,s) = & (d_{\mathcal{B}\bridgeGraph{n}} \otimes |\I|)\left(
\disp{\sum_{C \in \circles{\partial(r,s)}, s(C) = +} (-1)^{h(r)}\lefty{e_{C}} \otimes (r,s_{C}) } \right) \\
\ & \ \\
\ & = \sum_{C} \left( -\sum_{\gamma \in \actor{C}} e_{\gamma}e_{\gamma^{\dagger}} \right)\otimes (r,s_{C})
\end{array}
$$
Consequently, all the terms in $(\mu_{2} \otimes \I)(\I \otimes \delta)\delta$ cancel except those in the last case. However, we may cancel the remaining terms by adding $(d_{\mathcal{B}\bridgeGraph{n}} \otimes |\I|)\righty{\delta_{T}}$ which confirms the relation for a type $D$-structure. $\Diamond$

\section{Reidemeister Invariance of the $D$-structure}\label{sec:invariance}

\begin{thm}\label{thm:invariance}
Let $\righty{\tangle{T}}$ be an outside tangle with diagram $\righty{T}$. The homotopy class of the $D$-structure $\righty{\delta_{T}}: \rightComplex{\righty{T}} \rightarrow \mathcal{B}\bridgeGraph{n} \otimes_{\mathcal{I}} \rightComplex{\righty{T}}$ is an invariant of the tangle $\righty{\tangle{T}}$.
\end{thm}

\noindent {\bf Proof:}  The theorem follows immediately once we verify the invariance under application of Reidemeister moves and reordering the crossings. We do this in two separate lemmas.

\begin{lemma}\label{lem:order}
Reordering the elements of $\cross{\righty{T}}$ does not change the homotopy class of $\righty{\delta_{T}}$
\end{lemma}

\noindent{\bf Proof of lemma \ref{lem:order}:}  The formula for $\righty{\delta_{T}}$ in equation \ref{eqn:delta} decomposes as $\delta_{R} + \delta_{L}$ where $\delta_{L}$ includes those terms corresponding to edges in $\leftGraph{n}$, and $\delta_{R}$ contains all the other terms. The terms in $\delta_{L}$ have sign determined by the homological degree of the resolution ($h(r)$), while those in $\delta_{R}$ have sign determined by $I(r,\cross{\gamma})$, the number of $1$-resolutions used by $r$ for crossings occurring after the crossing $\cross{\gamma}$.  Thus changing the ordering does not change $\delta_{L}$.  \\
\ \\
\noindent Let $\mathfrak{o}_{i}$, $i = 1, 2$ be two orderings of $\cross{\righty{T}}$. Let $r$ be a resolution of $\righty{T}$. The permutation of the crossings determined by changing from $\mathfrak{o}_{1}$ to $\mathfrak{o}_{2}$ restricts to the crossings to which $r$ assigns a $1$ resolution as a permutation $\sigma_{r}$. We define a map from $\rightComplex{\righty{T}, \mathfrak{o}_{1}}$ to  $\rightComplex{\righty{T}, \mathfrak{o}_{2}}$ by taking $\Psi: (r,s) \rightarrow (-1)^{\mathrm{sgn}\sigma_{r}}(r,s)$. The sign on a term in $(\Psi \circ \delta_{L,1})(r,s)$ is $(-1)^{\mathrm{sgn}\sigma_{r}}(-1)^{h(r)}$ while that on $(\delta_{L,2} \circ \Psi)(r,s)$ is $(-1)^{h(r)} (-1)^{\mathrm{sgn}\sigma_{r}}$ as well, since no crossing is resolved differently when we change on the left of the $y$-axis. Consequently, $\Psi$ commutes with $\delta_{L}$. \\
\ \\
\noindent On the other hand, the terms in $(\Psi \circ \delta_{R,1})(r,s)$ and $(\delta_{R,2} \circ \Psi)(r,s)$ have signs $(-1)^{\mathrm{sgn}\sigma_{r_{\gamma}}}(-1)^{I_{1}(r,\cross{\gamma})}$ and
$(-1)^{I_{2}(r,\cross{\gamma})}(-1)^{\mathrm{sgn}\sigma_{r}}$, respectively. $\cross{\gamma}$ is an additional $1$-resolved crossing in $r_{\gamma}$, so $\mathrm{sgn}\sigma_{r_{\gamma}} - \mathrm{sgn}\sigma_{r}$ will equal the number of transpositions in $\sigma_{r_{\gamma}}$ involving $\cross{\gamma}$. As each transposition changes the location of $\cross{\gamma}$ by $1$ in the ordering, the number of such transpositions, modulo 2, equals the number of positions $\cross{\gamma}$ shifts, which is just $|I_{1}(r,\cross{\gamma}) - I_{2}(r,\cross{\gamma})|$. Thus the two products of signs will be the same. From this it follows that $\Psi$ also commutes with $\delta_{R}$. $\Diamond$\\
\ \\
\noindent We will use proposition \ref{prop:cancelD} to prove the second lemma.

\begin{lemma}\label{lem:inv}
The homotopy class of the $D$-structure $\righty{\delta_{T}}$  is invariant under the first, second, and third Redemeister moves.
\end{lemma}

\noindent{\bf Proof of lemma \ref{lem:inv}:} To prove that the homotopy class of $(\rightComplex{T}, \righty{\delta_{T}})$ (as a type $D$ structure) is invariant under the first, second, and third Reidemeister moves, we will repeatedly use the following observations. First, in the diagrams for each move there will be resolutions whose diagrams include a {\em free circle} which will either merge into another circle (free or cleaved) or be divided from another circle under the action of $\righty{\delta_{T}}$. When we merge the {\em free} circle into another circle, the free circle will be decorated with a $+$. In both $\mu$ (used for merging free circles) and $\lambda$ (used for merging the free circle into a cleaved circle) merging a $+$ decorated free circle acts as the identity and thus gives the idempotent $I_{\partial(r,s)}$ as the coefficient in $\mathcal{B}\bridgeGraph{n}$. For divisions we will gain a $-$ decorated free circle. Either $\Delta$ or $\delta$ will give a coefficient of $I_{\partial(r,s)}$ for this new configuration. Thus, the terms we will cancel occur in $\righty{\delta_{T}}$ as terms in the $x_{i} \longrightarrow I_{\partial(r,s)} \otimes d(x_{i})$ summand of $\righty{\delta_{T}}$.  Thus we can apply proposition \ref{prop:cancelD} to cancel these terms to obtain a homotopy equivalent $D$-structure. This $d$ structure will occur on same module as $\rightComplex{\righty{T}'}$ where $T'$ is the diagram which results after applying the Reidemeister move. To finish we will need to see that the new $D$ structure is identical with the $D$-structure $\righty{\delta_{T'}}$ coming directly from the diagram as in section \ref{sec:typeDcom}. 

\begin{center}
\begin{figure}
\includegraphics[scale=0.75]{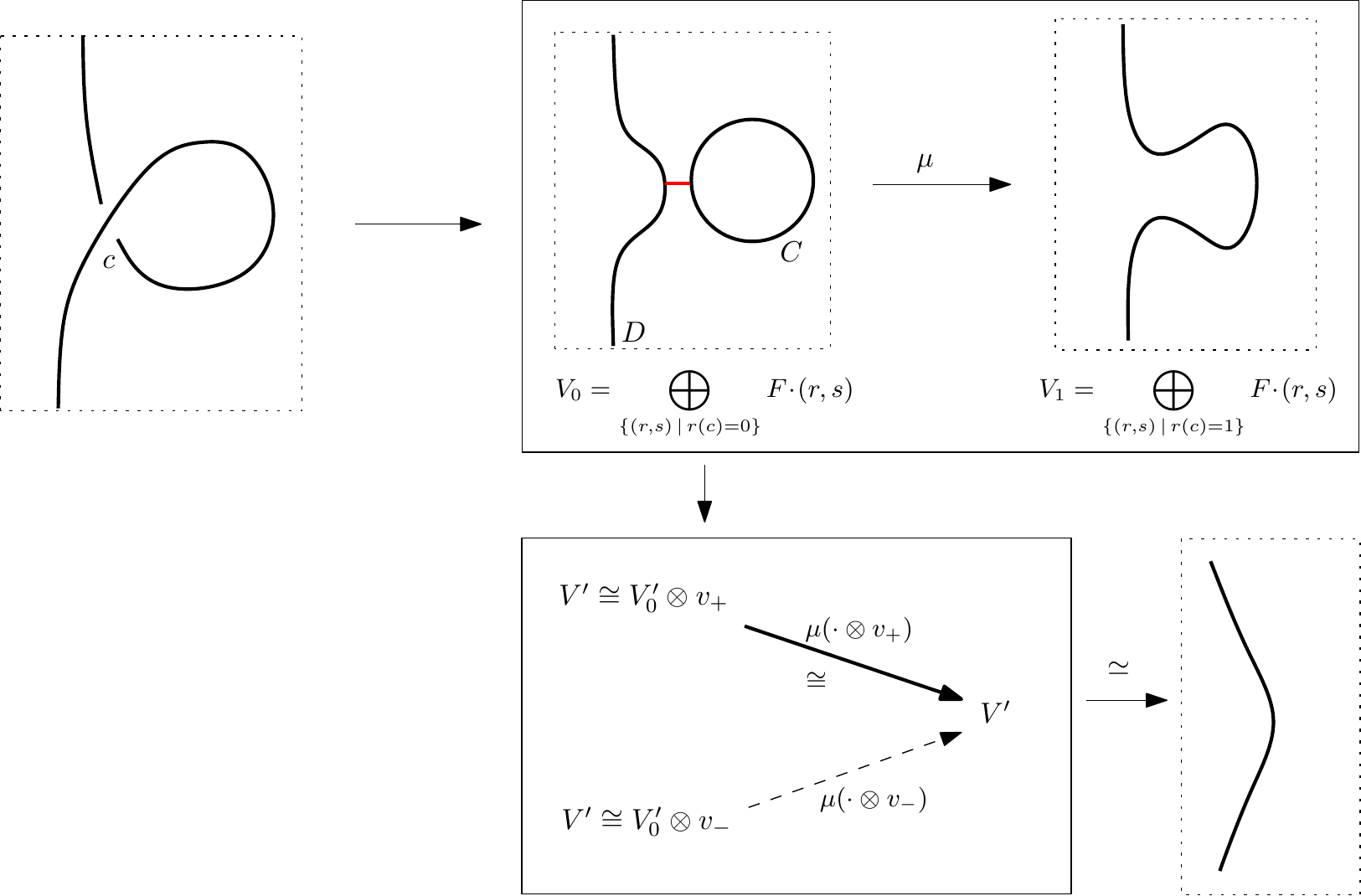}
\label{fig:RI}
\caption{}
\end{figure}
\end{center}

\noindent{\bf Reidemeister I:} Suppose that $\righty{T}$ has an RI reducible right-handed crossing $c \in \cross{\righty{T}}$, as (locally) depicted in Figure \ref{fig:RI}. Since we can reorder the crossings, we will take $c$ to be last (so the signs will be $+1$'s). $\rightComplex{\righty{T}}$ decomposes as $V_{0} \oplus V_{1}$ corresponding to the states where $\rho(c) = 0$ or $\rho(c) = 1$. Since the crossing is right-handed, the resolutions generating $V_{0}$ will have an additional {\em free} circle $C$. Changing $\rho$ from $0$ to $1$ at $c$ merges this circle into the other local component $D$, which may belong to either a free circle or a cleaved circle. Using the decoration $s(C)$ we will further divide $V_{0}$ as $V'_{0} \otimes v_{+}$ and $V'_{0} \otimes v_{-}$. The APS-differential $d_{APS}$ induces an isomorphism $I: V'_{0} \otimes v_{+} \rightarrow V_{1}$. We may identify $V_{1}$ with $V_{0}'$ so that when using the states as generators $I$ is a diagonal $\pm 1$ matrix. Each non-zero entry in $I$ gives rise to a term in $\righty{\delta_{T}}$ which can be canceled using proposition \ref{prop:cancelD}. More formally, if $(r,s_{+})$ has $r(c)=0$, $s_{+}(C) = +$, then $\righty{\delta_{T}}(r,s_{+})$ has one term from a generator in $V_{1}$, which is $ I_{\partial(r,s_{+})}\otimes (r_{c},s')$ where $s'(D)=s_{+}(D)$ and $\partial(r_{c},s')=\partial(r,s_{+})$. Since $I_{\partial(r,s_{+})}$ is an idempotent, we may cancel this term. \\
\ \\
\noindent If we cancel all the terms from $I$ we will be left solely with $V' = V'_{0} \otimes v_{-}$, with a new $D$-structure $\delta': V' \rightarrow \mathcal{B}\Gamma_{n} \otimes V'$. $\delta'$ is the image of $\delta|_{V'}$ in $V'$, plus a perturbation term. Proposition \ref{prop:cancelD} also describes how to compute this perturbation term. It arises from the image of $V'_{0} \otimes v_{-} \rightarrow V_{1}$, when the local arc $D$ in Figure \ref{fig:RI} belongs to a circle decorated with $+$. There are two cases to consider:
\begin{enumerate}
\item If $D$ belongs to free circle, then the image of $(r,s)$ with $r(c) = 0$, $s(C)=-$ and $s(D) = +$ in $V_{1}$ will be $I_{\partial(r,s)} \otimes (r_{c},s_{D})$,  where $s_{D}$ prescribes the same decorations with $s_{D}(D)=-$. This is the result of $d_{APS}$, and $\partial (r_{c},s_{D}) = \partial(r,s)$. However, this will be canceled by $(r,s_{+,D})$ where $s_{+,D}(D)=-$ and $s_{+,D}(C)=+$. The perturbation term is therefore the sum of terms of the form $- (a\cdot I_{(r,s_{+,D})}^{-1} \cdot I_{\partial(r,s)}) \otimes (r'',s'')$ where $a \otimes (r'',s'')$ is a different term in $\righty{\delta_{T}}(r,s_{D,+})$. The other terms, however, involve crossing or decoration changes outside the local picture (since $D$ belongs to a free circle), and thus $r''(c) = 0$ and $s''(C) = +$. So $(r'',s'')$ be a state generating $V_{0}\otimes v_{+}$ and thus be used to cancel another state. In the canceling process, eventually $(r'',s'')$ will be eliminated, so these perturbation terms will vanish. Thus, after all the cancellations are performed, the perturbation terms from this case will be trivial. 
\item When $D$ belongs to a cleaved circle then the image of $(r,s)$ with $r(c) = 0$, $s(C)=-$ and $s(D) = +$ in $V_{1}$ will be $\righty{e_{D}} \otimes (r_{c},s_{D})$. $I_{\partial(r_{c},s_{D})} \otimes (r_{c},s_{D})$ is the term in $\righty{\delta_{T}}$ which comes from $d_{APS}(r,s_{+,D})$. The perturbation term is therefore the sum of terms of the form $- (a\cdot I_{(r,s_{+,D})}^{-1} \cdot \righty{e_{D}}) \otimes (r'',s'')$ where $a \otimes (r'',s'')$ is a different term in $\righty{\delta_{T}}(r,s_{D,+})$. Once again $(r'',s'')$ must be a generator in $V_{0}\otimes v_{+}$ and thus be canceling. So again, the perturbation term is trivial. 
\end{enumerate}
\noindent Consequently, the perturbation term will be $0$ and $\delta'$ is just the projection of $\delta|_{V'}$ to the terms with states generating $V'$. Since we chose $c$ to be last in the crossing ordering, and $C$ can only interact with the rest of the diagram at $c$, $\delta'$ is the same $D$-structure we get from the diagram with $C$ removed, shifted by the quantum grading from $s(C) = -$. As usual, the shifts based on positive/negative crossings will realign the bigrading. \\
\ \\
\noindent We can obtain the invariance for removing a left-handed RI move from $\righty{T}$ by a combination of introducing a new crossing using a right-handed RI move and then simplifying using a RII move. We thus proceed to showing invariance under RII, thereby also showing the invariance under left-handed RI moves.\\
\ \\ 

\begin{figure}
\begin{center}
\includegraphics[scale=0.5]{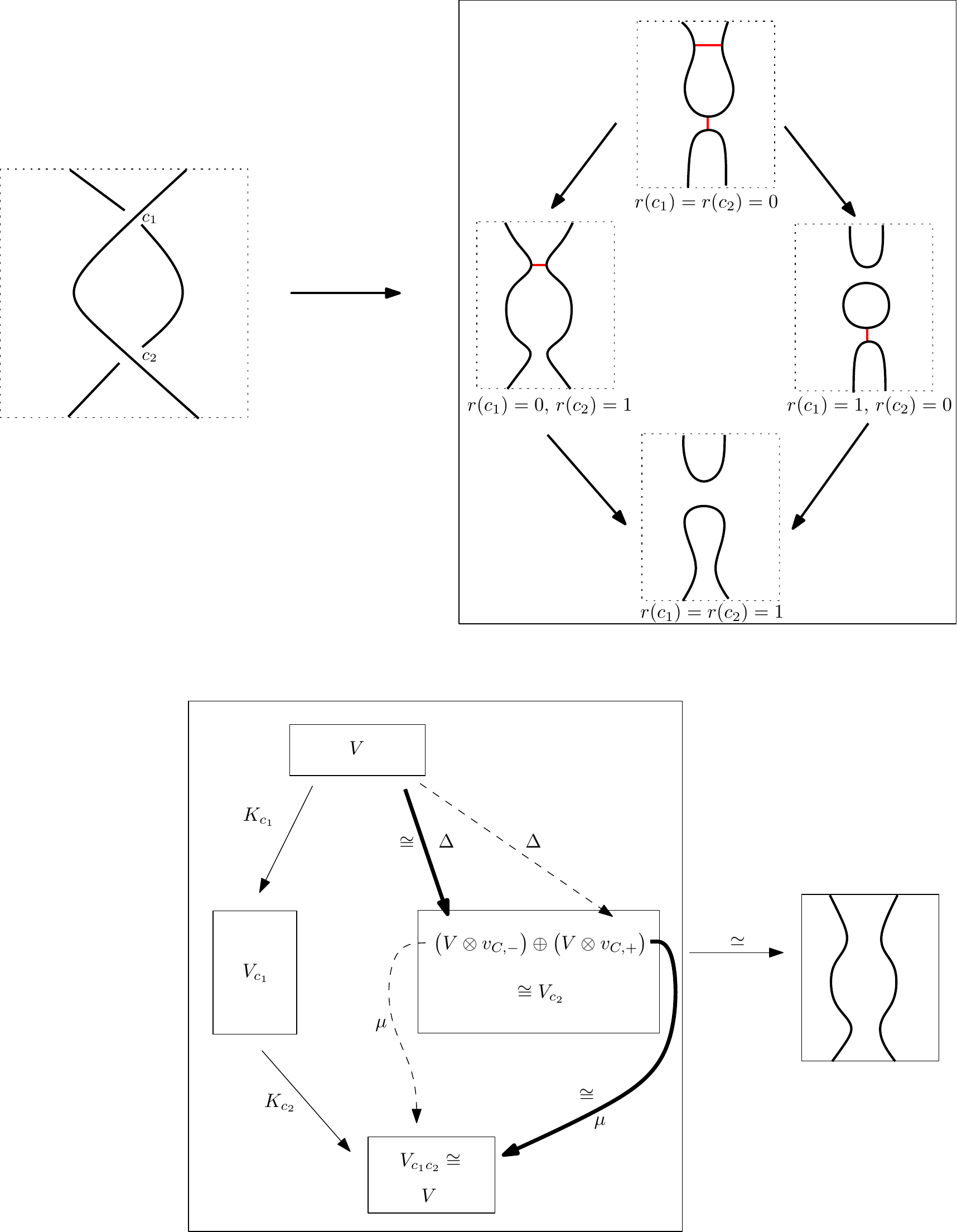}
\end{center}
\label{fig:RII}
\caption{}
\end{figure}

\noindent{\bf Reidemeister II:} Much the same argument can be made for why the standard proof of invariance under the second Reidemeister move will continue to work in our context. We will use the notation from Figure \ref{fig:RII}. Let $c_{0}$ and $c_{1}$ be the local crossings in a diagram for $\sigma_{1}\sigma_{1}^{-1}$ where resolving both with $\rho(c_{i}) = i$ yields the identity braid, which will be contained in $\righty{T}'$. Let $V_{ij}$ be the resolutions $(r,s)$ where $c_{0}$ is resolved with $i$ and $c_{1}$ is resolved with $j$. In $V_{10}$ there is a free circle $C$. We divide $V_{10}$ into two parts $V'_{10} \otimes v_{+}$ and  $V'_{10} \otimes v_{-}$ based on the decoration of $C$. Then $I_{\partial(r,s)} \otimes d_{APS}$ induces isomorphisms $I_{1}: V_{00} \rightarrow \I \otimes \left( V'_{10} \otimes v_{-}\right)$ (from $\Delta$) and $I_{2}: V'_{10} \otimes v_{+} \rightarrow \I \otimes V_{11}$ (from $\mu,\lambda$). \\
\ \\
\noindent When we cancel $I_{2}$ we can obtain a perturbation of $\righty{\delta_{T}}$ on elements of $V_{00} \oplus V_{01} \oplus \big(V'_{10} \otimes v_{-}\big)$ with image in this space (tensored with the algebra). The image of $\righty{\delta_{T}}$ on states generating 
$V'_{10} \otimes v_{+}$, however, must either resolve a crossing, or change the decoration on a cleaved circles. In particular, the free circle $C$ will remain, and have decoration $+$, except when we resolve $c$ to get a state in $V_{11}$. As those terms occur in $I_{2}$, the perturbation term will be supported in $V'_{10} \otimes v_{+}$. After canceling all of $I_{2}$, no such generators remain. Thus, after canceling $I_{2}$, the new type $D$ structure is just $\righty{\delta_{T}}$ restricted to $V_{00} \oplus V_{01} \oplus \big(V'_{10} \otimes v_{-}\big)$, with image projected to the same. However, no further perturbation of $\delta|_{V_{01}}$ can result from the cancellation of $I_{1}$ since there are no terms in the image of $\delta|_{V_{01}}(\xi)$ supported on elements of $V'_{10} \otimes v_{-}$ (due to the resolutions of the crossings).  \\
\ \\
\noindent Thus, after canceling both isomorphisms we obtain the homotopy equivalent $D$-structure $(V_{01}, \righty{\delta_{T}}|_{V_{01}})$. It remains to observe that this is identical with $(\rightComplex{\righty{T}'}, \righty{\delta_{T'}})$, as the states are identical in the two cases, and the calculation of $\righty{\delta_{T}}|_{V_{01}}$ only involves the crossings external to the local diagrams. \\
\ \\
\begin{figure}
\begin{center}
\begin{tabular}{ccc}
\includegraphics[scale=0.4]{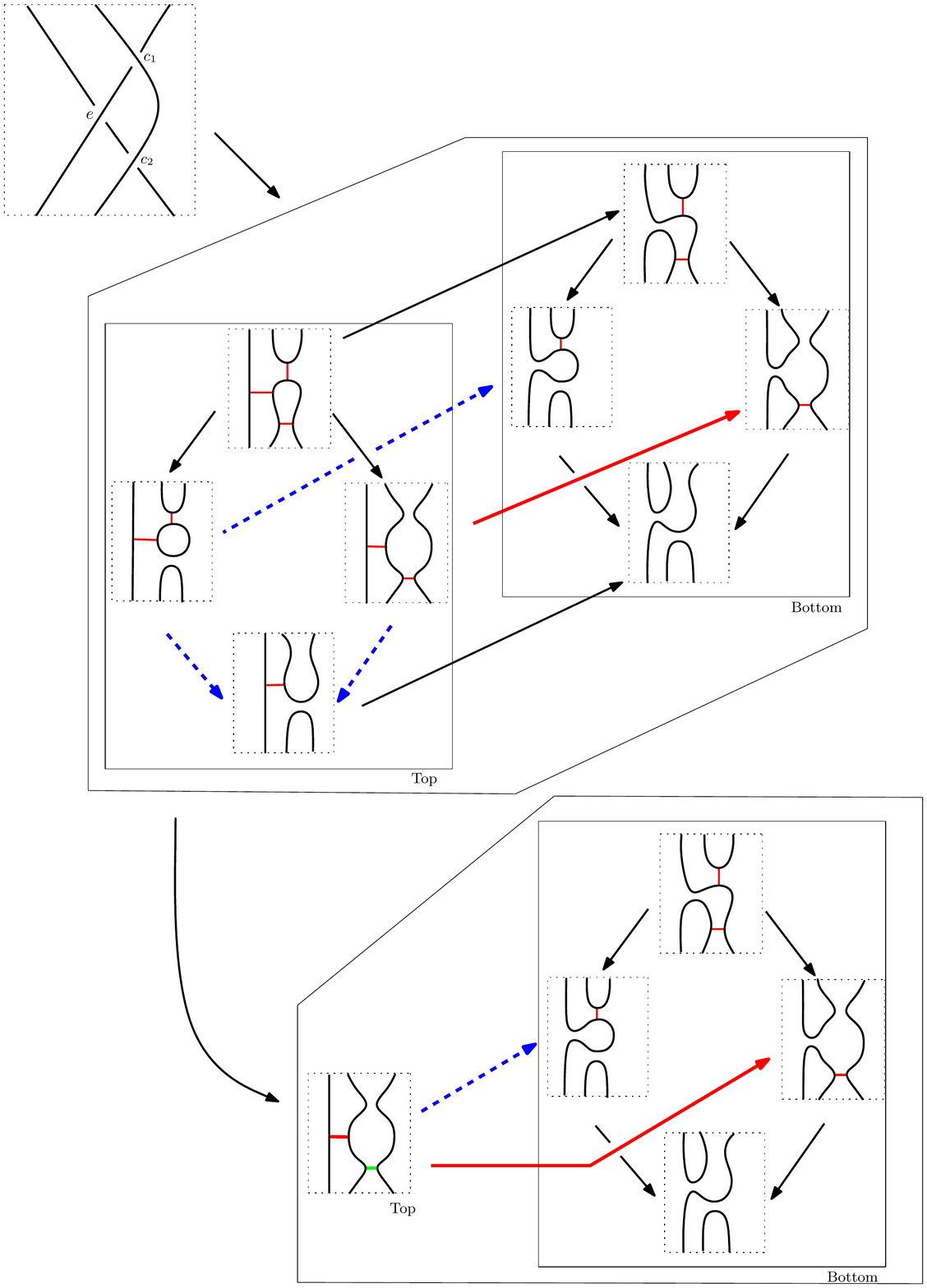}&
\hspace{0.33in}&
\includegraphics[scale=0.4]{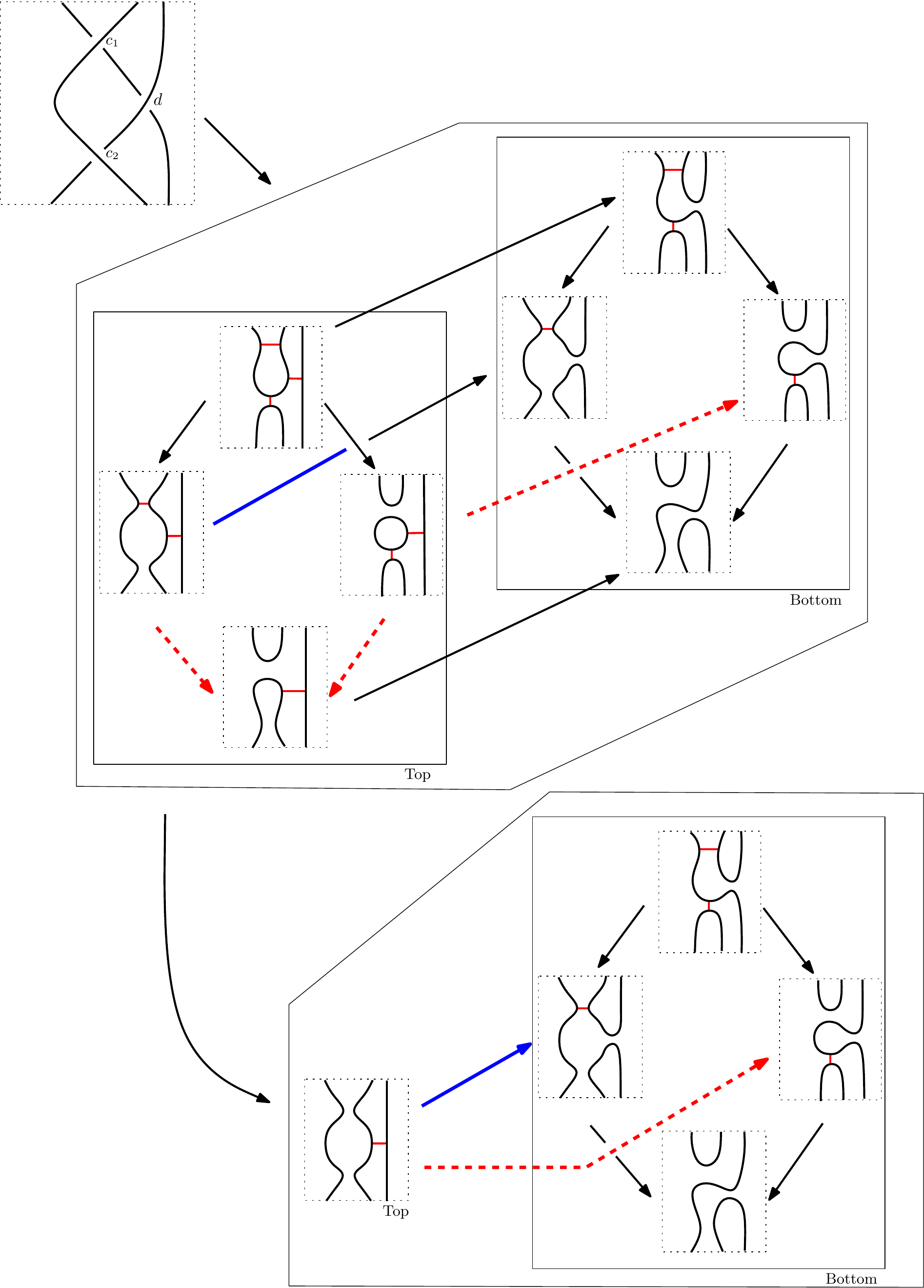}
\end{tabular}
\end{center}
\label{fig:RIII}
\caption{}
\end{figure}

\noindent{\bf Reidemeister III:} As with the RI and RII moves, we modify the usual proof of the invariance of Khovanov homology under an RIII move. The complexes before and after the RIII move are depicted schematically in \ref{fig:RIII}. We observe, as usual, that 1) the diagrams in the ``top'' rectangle, in both middle rows, replicate the local resolution diagrams from an RII move using $c_{1}$ and $c_{2}$, and 2) the diagrams in the ``bottom'' rectangles are isotopic in each of the four corners, including isotopies of the local bridges. Consequently,since they cannot change the resolutions at $e$ and $d$ from $1$ to $0$ $\righty{\delta_{T^{L}}}$ and  $\righty{\delta_{T^{R}}}$ restrict to the bottom rectangles to give type $D$ structures, and due to the isotopy will be identical on them. Furthermore, we can follow the cancellation sequence for the RII move to simplify the ``top'' rectangle. Again, the resulting maps will be identical {\em when considering only the top rectangles}. This leaves the maps between the top and bottom to consider. Our goal is to see, after the cancellations, that these maps are also identical.  However, to do that, we specify that $c_{2} < e < c_{1}$ on the left and $ c_{2} < d < c_{1}$ on the right, and that these are the terminal sequences in each ordering.\\
\ \\
\noindent If we look at the middle row, on the left, we will need to compute the contribution of the blue arrows when we cancel the states $(r,s)$ for the diagram with the free circle $C$, and $s(C) = +$. We want these to equal the term in the lower row, on the left, coming from a resolution change along the green bridge (the lowest one in ``top''). A similar argument for the diagrams on the right will then show that the two arrows on the left and right are the same, corresponding to the terms which come from changing the resolutions at the red and green bridges. The blue arrows form the following configuration, with the coefficients listed next to the arrow:

\begin{center}
\includegraphics[scale=0.6]{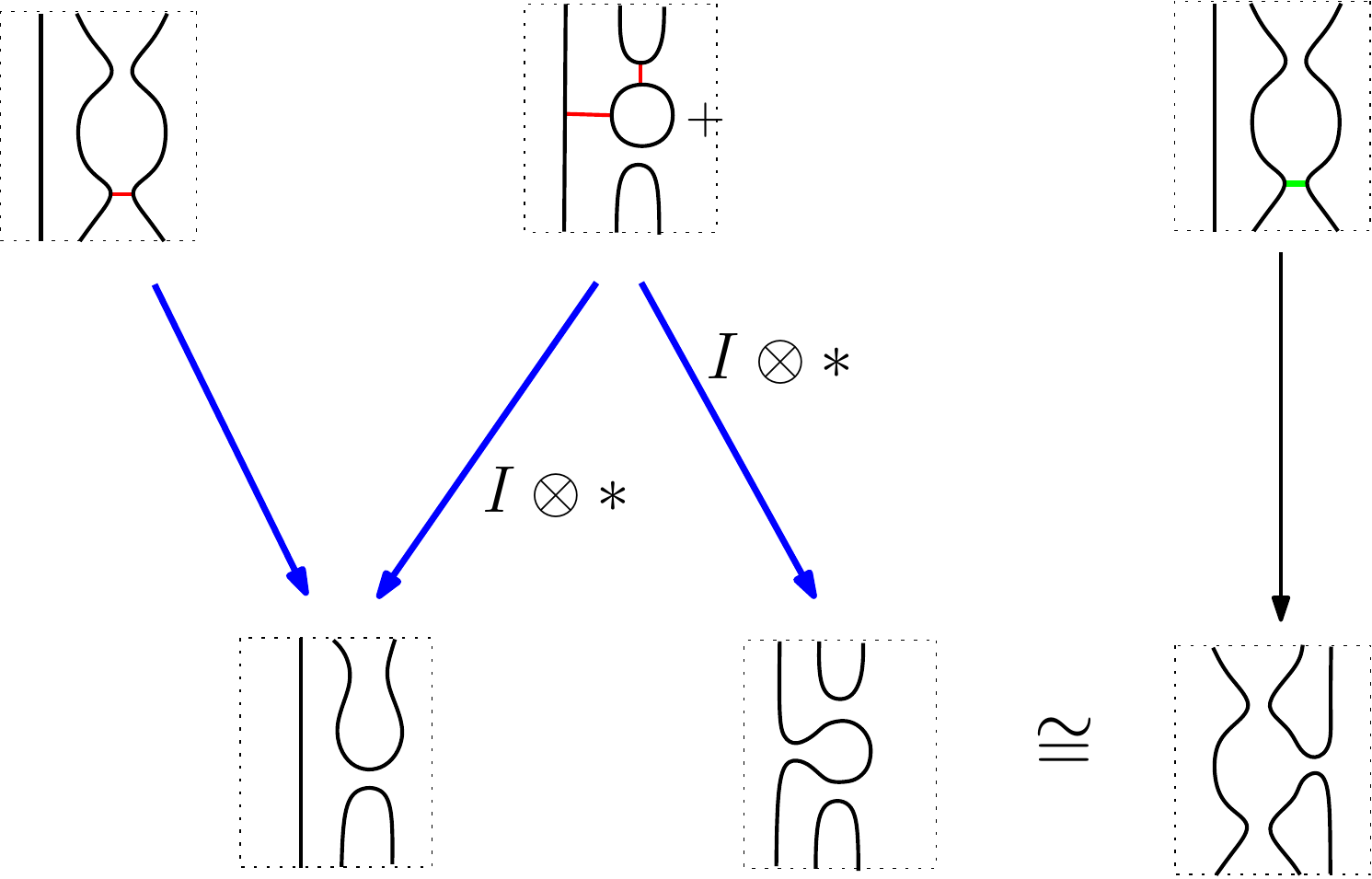}
\end{center}

\noindent Note that the arrow pointing down and to the left corresponds to the resolution code change $100 \rightarrow 101$ and thus has a positive coefficient. The arrow starting at the same point, and down and to the right, corresponds to $100 \rightarrow 110$, and thus also has a positive coefficient. The remaining arrow can be more complicated, but corresponds to $001 \rightarrow 101$, and thus is $- e \otimes \xi$ where $e \otimes \xi$ comes from the change in the red bridge. Canceling to get a term from $001 \rightarrow 110$ multiplies this by $-1$ (since the other arrow have idempotent coefficients). Thus, the new term is the same as we would get from the green bridge on the right. This latter diagram is in the complex for the $\righty{T}_{R}$. There it corresponds to the resolution change $100 \rightarrow 110$, and thus occurs with the same sign. Hence, these two arrows will be the same in the reduced complexes. \\
\ \\
\noindent We can repeat this argument for the arrow resulting from cancellation in $\righty{T}_{R}$. We wish to see that the bold arrow on the right in Figure \ref{fig:RIII} equals the term coming from $001 \rightarrow 011$ on the left, which is the term from changing the red bridge on the left, bottom times $-1$ due to the signs. On the right, the cancellations which give the map between isotopic diagrams are $100 \stackrel{+1}{\rightarrow} 101 \stackrel{-1}{\leftarrow} 001 \stackrel{-1}{\rightarrow} 011$, which, along with the $-1$ from the cancellation formula, means that the corresponding map for $\righty{T}_{R}$ is also $-1$ times the same bridge transition. Consequently, these two maps will also be equal.   
$\Diamond$

\section{Examples}\label{sec:examples}

\noindent We give three examples of calculations of the type $D$ structure. \\
\ \\
\noindent{\bf Example 1:} First, we consider the two planar tangles depicted below:
$$
\inlinediag[0.6]{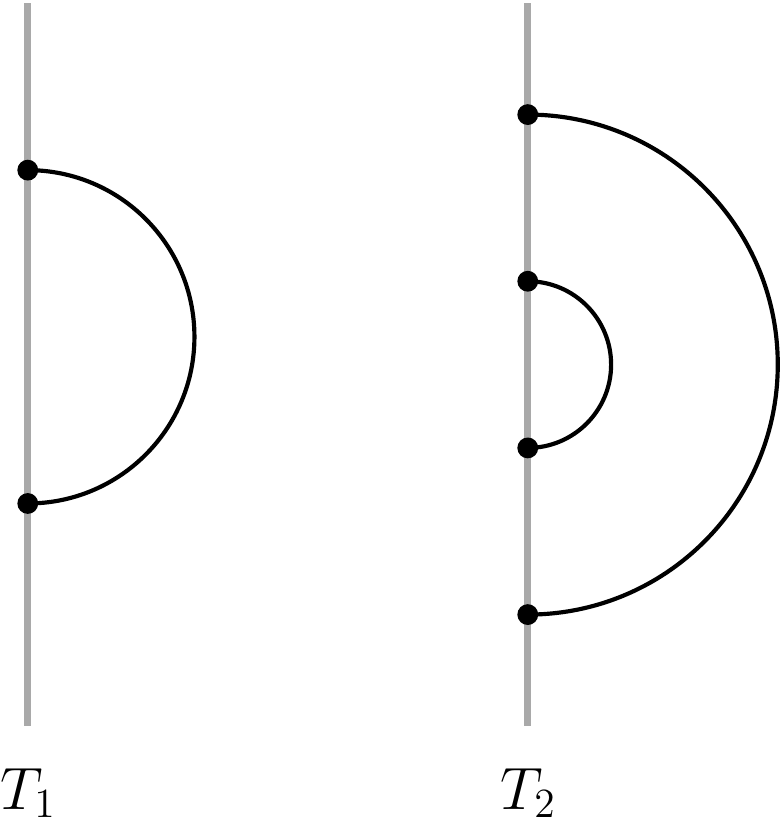}
$$
\noindent In neither are there any crossings to resolve. So the resolutions come solely from pairing with planar matchings in $\leftHalf$. For the tangle $T_{1}$, there is a unique matching $\lefty{m}$ on $P_{2}$. Thus, there is only one resolution. This resolution has one cleaved circle, which can be decorated with either a $+$ or a $-$. We will call the resulting states $s_{+}$ and $s_{-}$. $s_{\pm}$ will be in bigrading $(0,\pm 1/2)$. There are no bridges for $\lefty{m}$, so the entire sum defined $\righty{\delta}(s_{+})$ reduces to $\lefty{e_{C}}\otimes s_{-}$, with sign $(-1)^{0}=+1$ since $0$ is the homological grading. Furthermore, $\righty{\delta}(s_{-}) = 0$. Note that $\lefty{e_{C}}$ has bigrading $(1,1)$, so the bigrading of $\righty{\delta}(s_{+})$ is $(1,1-1/2)=(1,1/2)$, which is $(1,0)$ greater than that for $s_{+}$.\\
\ \\
\noindent For the right tangle $T_{2}$ with endpoints $P_{2}$, shown in the diagram above, there are two possible left matchings to pair with it: $\lefty{m}_{1}$ will pair the same points in $P_{2}$ as $T_{2}$, whereas $\lefty{m}_{2}$ will not. The states for this tangle are then
\begin{enumerate}
\item For $\lefty{m}_{1} \# T_{2}$ there are two circles $C$, the outer circle, and $D$ the inner circle. There are thus four states: $s^{1}_{++}$, $s^{1}_{+-}$, $s^{1}_{-+}$, $s^{1}_{--}$, where the first sign is the decoration on $C$ and the second is that on $D$. These occur in bigradings $(0,1)$, $(0,0)$, $(0,0)$, and $(0,-1)$, respectively. 
\item For $\lefty{m}_{2} \# T_{2}$ there is one circles $E$ and two states: $s^{2}_{+}$, $s^{2}_{-}$, corresponding to decoration on $E$. These occur in bigradings $(0,1/2)$ and $(0,-1/2)$, respectively. 
\end{enumerate}
Furthermore, there is one isotopy class of bridges in $\lefty{m}_{1}$, which we will denote by 
$\lefty{\gamma}$. Surgery on $\lefty{\gamma}$ in $\lefty{m}_{1}$ produces $\lefty{m}_{2}$. Thus we will denote the sole bridge in $\lefty{m}_{2}$ by $\lefty{\gamma}^{\dagger}$. \\
\ \\
\noindent We can now write down $\righty{\delta}$:
\begin{equation}
\begin{split}
\righty{\delta}(s^{1}_{++}) &= \lefty{\gamma}\otimes s^{2}_{+} + \lefty{e_{C}} \otimes s^{1}_{-+} + \lefty{e_{D}} \otimes s^{1}_{+-} \\
\righty{\delta}(s^{1}_{+-}) &= \lefty{\gamma}\otimes s^{2}_{-} + \lefty{e_{C}} \otimes s^{1}_{--} \\
\righty{\delta}(s^{1}_{-+}) &= \lefty{\gamma}\otimes s^{2}_{-} + \lefty{e_{D}} \otimes s^{1}_{--}\\
\righty{\delta}(s^{1}_{--}) &= 0\\
&\ \\
\righty{\delta}(s^{2}_{+})\ &= \lefty{\gamma_{1}}^{\dagger}\otimes s^{1}_{+-} + \lefty{\gamma_{2}}^{\dagger}\otimes s^{1}_{-+} + \lefty{e_{E}} \otimes s^{2}_{-}\\
\righty{\delta}(s^{2}_{-})\ &= \lefty{\gamma}^{\dagger}\otimes s^{1}_{--} \\
\end{split}
\end{equation}
We remind the reader that $\gamma^{\dagger}$ on $s^{2}_{+}$ gives two different elements of the algebra for $P_{2}$, since the corresponding edges will have two different targets. Note that  $\lefty{\gamma}$ has bigrading $(1,1/2)$, so $\lefty{\gamma}\otimes s^{2}_{+}$ has bigrading $(1,1)$, which is $(1,0)$ more than that for $s^{1}_{++}$. On the other hand $s^{1}_{--}$ is in grading $(0,-1)$, so $\lefty{\gamma}^{\dagger}\otimes s^{1}_{--}$ is in grading $(1,-1/2)$ which is $(1,0)$ more than the bigrading of $s^{2}_{1/2}$. \\
\ \\
\noindent We can use this example to illustrate the occurence of the relations in $\mathcal{B}\Gamma_{n}$. For instance,
\begin{align*}
(\mu \otimes \I)(\I \otimes \righty{\delta})\righty{\delta}(s^{1}_{++}) = (\lefty{\gamma}\lefty{\gamma}_{1}^{\dagger})&\otimes s^{1}_{+-} + 
(\lefty{\gamma}\lefty{\gamma}_{2}^{\dagger})\otimes s^{1}_{-+} + 
(\lefty{\gamma}\lefty{e_{E}} + \lefty{e_{C}}\lefty{\gamma} + \lefty{e_{D}}\lefty{\gamma})\otimes s^{2}_{-} \\
\ &+ (\lefty{e_{D}}\lefty{e_{C}} + \lefty{e_{D}}\lefty{e_{C}})\otimes s^{1}_{--}
\end{align*}
However, since $\lefty{e_{D}}$ and $\lefty{e_{C}}$ are odd elements, $\lefty{e_{D}}\lefty{e_{C}} + \lefty{e_{D}}\lefty{e_{C}}=0$. Furthemore, $\lefty{\gamma}\lefty{e_{E}} + \lefty{e_{C}}\lefty{\gamma} + \lefty{e_{D}}\lefty{\gamma} = 0$ by relation \ref{rel:lefty3}. The remaining two terms will cancel with $$(d_{\Gamma_{2}} \otimes |\I|)\righty{\delta}(s^{1}_{++}) = d\lefty{e_{C}} \otimes |s^{1}_{-+}| + d\lefty{e_{D}} \otimes |s^{1}_{-+}|$$
since $d\lefty{e_{C}} = - \lefty{\gamma}\lefty{\gamma}_{1}^{\dagger}$ and $s^{1}_{-+}$ is in bigrading $(0,0)$.  Likewise,  $d\lefty{e_{D}} = - \lefty{\gamma}\lefty{\gamma}_{2}^{\dagger}$.\\
\ \\
\noindent{\bf Example 2:} We consider the right tangle below, which is the left handed trefoil ($n_{-} = 3, n_{+} = 0$), with one arc removed. This has boundary $P_{2}$, so we will use $\mathcal{B}\Gamma_{2}$ for the algebra. \\
$$
\inlinediag[0.6]{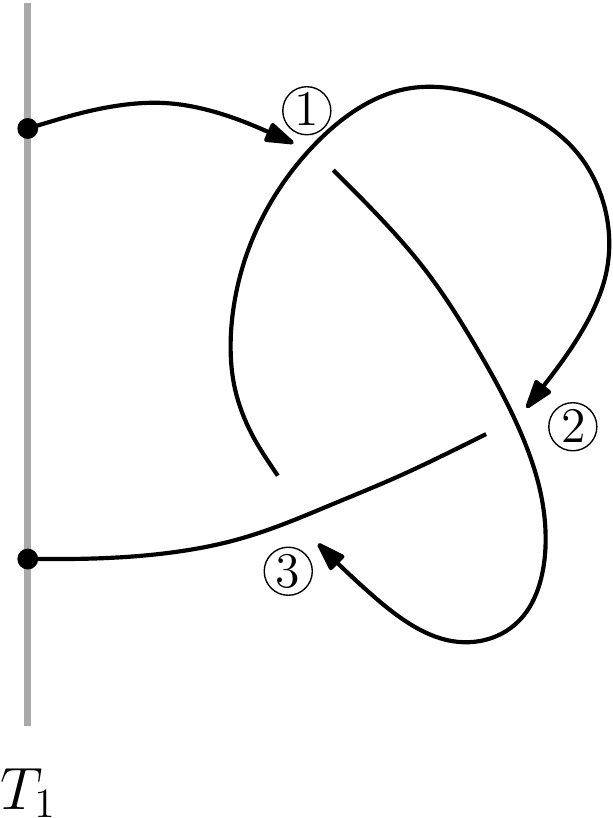}
$$
\ \\
\noindent If we group all the states based on whether the cleaved circle is decorated with a $+$ or a $-$ we obtain two copies of the APS-complex. In this case, the APS complex is isomorphic, up to grading shifts, with the reduced Khovanov homology of the left-handed trefoil. There are many cancellations used to calculate this reduced homology, and they can be implemented here to give a homotopy equivalent type $D$ structure. When we do this, we obtain a description of $\rightComplex{\righty{T}}$ up to homotopy,  which is the isotopy invariant of the tangle. \\
\ \\
\noindent We give a brief, tedious explanation of the computations: we will write our states as pairs $(r,d)$ where $r$ is an element of $\{0,1\}^{3}$ and specifies the resolution at each of the crossings (as labeled above), and $d$ specifies the decoration on each circle. We will write the decoration on the cleaved circle first, followed by the free circles (if there are any) from the top of the diagram to the bottom. \\
\ \\
\begin{enumerate}
\item There are eight states for $r=000$, but in $d_{APS}$ we have a $+1$ coefficient for $(000,+++) \rightarrow (010,++)$, and no other state in $000$ has non-zero coefficient in $\righty{\delta}$ on $(010,++)$. We can thus cancel these two states. A similar argument holds for $(000,++-) \rightarrow (100, +-)$. After this cancellation, the only state with non-zero coefficient on $(010,+-)$ will be $(00,+-+)$, so we may also cancel these. This leaves only $(000,+--)$ with decoration $+$ on the cleaved circle for $r=000$. Its image under $\righty{\delta}$, after the cancellations is $\lefty{e_{C}}\otimes (000,---) + \righty{e_{C}} \otimes (100,--) + \righty{e_{C}} \otimes (001,--)$. 
\item A similar argument applies to $(000,-++) \rightarrow (001,-+)$. However, when we cancel $(000,--+) \rightarrow (010,--)$ we encounter that $(000,-+-) \rightarrow (010,--)$ has $+1$ coefficient. $(000,--+)$ also has a non-zero term $(001,--)$, so after cancellation the image of $(000,-+-)$ becomes $(100,--) - (001,--)$. If we cancel along $(000,--+) \rightarrow (100,--)$ we need to adjust the image of $(000,+--)$ since it has non-zero coefficient on $\righty{e_{C}} \otimes (100,--)$. The adjustment is to remove  {\em subtract} $\righty{e_{C}} \otimes (100,--) + (-1)(1) \righty{e_{C}} \otimes (001,--)$ from $\righty{\delta} (000,+--)$. The image of $(000,+--)$ is now  $\lefty{e_{C}}\otimes (000,---) + 2 \righty{e_{C}} \otimes (001,--)$. No further cancellation is possible in homological degree $0$ since this image has algebra coefficients, and $(000,---)$ has trivial image under $\righty{\delta}$. \\
\item Remaining in homological degree $1$ are the six generators $(100,++)$, $(001,++)$, $(001,+-)$, $(100,-+)$, $(010,-+)$, and $(001,--)$. $\righty{\delta}(100,++) = -\lefty{e_{C}} \otimes (100,-+) + (110,+) + (101,+)$, so we will cancel along $(110,+)$. This state occurs in no other term in the image of $\righty{\delta}$. Likewise, $(001,++) \rightarrow (101,+)$ occurs with coefficient $-1$ (due to the crossing ordering), and once $(100,++)$ has been canceled, canceling along this term introduces no alteration to $\righty{\delta}$.
\item  Now $\righty{\delta}(001,+-) = -\lefty{e_{C}} \otimes (001,--) -\righty{e_{C}}\otimes\big[(101,-) + (011,-)\big]$ and there is a non-zero term $(011,+) \rightarrow \lefty{e_{C}}\otimes (011,-)$. When we cancel along the second term in $d_{APS}(010,-+) = -(110,-) + (011, -)$, we obtain  $\righty{\delta}(001,+-) = -\lefty{e_{C}} \otimes (001,--) -\righty{e_{C}}\otimes\big[(101,-) + (110,-)\big]$ for the new type $D$ structure, and $(011,+) \rightarrow \lefty{e_{C}}\otimes (110,-)$ as the replacement term. However, $\righty{\delta}(100,-+) = d_{APS}(100,-+) = (101,-) + (110,-)$. Canceling on the $(101,-)$ term results in $\righty{\delta}(001,+-) = -\lefty{e_{C}} \otimes (001,--)$, and $(011,+) \rightarrow  \lefty{e_{C}}\otimes (110,-)$. Only $(100,-+)$ and $(001,--)$ remaing in homological grading $1$. 
\item Remaining in homological grading $+2$ are the generators $(011,+)$ and $(110,-)$. The other four generators have been canceled already. $\righty{\delta}(011,+) = \lefty{e_{C}}\otimes (110,-) + \righty{e_{C}}(111,-+) + (111,+-)$, while $\righty{\delta}(110,-) = (111,--)$. We cancel along both the $+1$ coefficient terms. This annihilates the generators in homological grading $+2$. 
\item Remaining in homological grading $+3$ are $(111,++)$ and $(111,-+)$. $\righty{\delta}(111,++) = - \lefty{e_{C}} \otimes (111,-+)$.
\end{enumerate}
We are left with six generators. We shift by $(-n_{-}, n_{+} - 2n_{-}) = (-3,-6)$ and relabel  for ease of reference:
\begin{enumerate}
\item In homological grading $-3$ there are $s^{-3}_{+} = (000,+--)$ and $s^{-3}_{-} = (000,---)$ with bigradings $(-3,-8 \pm 1/2)$. Furthermore,
\item In homological grading $-2$ there are $s^{-2}_{+} = (001,+-)$ in bigrading $(-2,-11/2)$ and $s^{-2}_{-} = (001,--)$ in bigrading $(-2,-13/2)$. (Recall that there is an additional $+1$ shift in the quantum grading due to this corresponding to homological grading $+1$.)
\item In homological grading $0$ there are $s^{0}_{+} = (111,++)$ in bigrading $(0,-3/2)$ and $s^{0}_{-} = (111,-+)$ in bigrading $(0,-5/2)$.
\end{enumerate}
The map $\righty{\delta}$ is then non-zero on the $+$ generators only, where
\begin{equation}
\begin{split}
\righty{\delta}(s^{-3}_{+}) &= 2\,\righty{e_{C}} \otimes s^{-2}_{-} + \lefty{e_{C}} \otimes s^{-3}_{-} \\
\righty{\delta}(s^{-2}_{+}) &= - \lefty{e_{C}} \otimes s^{-2}_{-}\\
\righty{\delta}(s^{0}_{+})\ &= - \lefty{e_{C}} \otimes s^{0}_{-}\\
\end{split}
\end{equation}
We notice that the bigrading on $\righty{e_{C}} \otimes s^{-2}_{-}$ is $(0,-1) + (-2, -13/2) = (-2, -15/2)$ while that on $\lefty{e_{C}} \otimes s^{-3}_{-}$ is $(1,1) + (-3, -17/2) = (-2,-15/2)$. Both are $(1,0)$ bigger that $(-3,-15/2)$, which is the bigrading on $s^{-3}_{+}$.\\
\ \\
\noindent More interesting is the occurrence of $2$ as the coefficient. This will ultimately be responsible for the $2$-torsion in the full Khovanov homology of the left-handed trefoil, when we specify how to use $\rightComplex{\righty{T}}$ to compute this homology. \\
\ \\
 
\appendix
\section{Proof of the cancellation lemma for type $D$ structures}\label{sec:cancelDProof}

\noindent Here we give the algebraic calculation underlying proposition \ref{prop:cancelD}, which is repeated below. We employ the definition of type $D$ structures given in section \ref{sec:typeD}.

\begin{prop}[Cancellation]
Let $\delta$ be a $D$-structure on $N$. Suppose there is a basis $B$ for $N$ where $\delta$ can be described by structure coefficients satisfying $a_{ii} = 0$ and $a_{12} = 1_{A}$. Let $\overline{N} = \mathrm{span}_{R}\{\overline{x}_{3}, \ldots, \overline{x}_{n}\}$. Then 
$$
\overline{\delta}(\overline{x}_{i}) = \sum_{j \geq 3}(a_{ij} - a_{i2}\,a_{1j}) \otimes \overline{x}_{j}
$$
defines a $D$-structure on $\overline{N}$. Furthermore, the maps
$$
\begin{array}{lcl}
\iota : \overline{N} \rightarrow A \otimes N & \hspace{0.75in} & \iota(\overline{x}_{i}) = 1_{A} \otimes x_{i} -  a_{i2} \otimes x_{1} \\
\ &\ \\
\pi :  N \rightarrow A \otimes \overline{N} & \hspace{0.75in} & \pi(x_{i}) = \left\{
\begin{array}{ll}
0 & i = 1\\
-\sum_{j \geq 3}a_{1j} \otimes \overline{x}_{j} & i = 2\\ 
1_{A} \otimes \overline{x}_{i} & i \geq 3
\end{array}
\right.
\end{array}
$$
realize $\overline{N}$ as a deformation retraction of $N$ with $\iota \ast \pi \simeq_{H} \I_{N}$ using the {\em homotopy} $H: N \rightarrow A \otimes N[-1]$
$$
H(x_{i}) = \left\{\begin{array}{ll} -1_{A} \otimes x_{1} & i = 2\\ 0 & i \neq 2 \end{array}\right.
$$ 
\end{prop}
\ \\
\ \\
\noindent{\bf Proof:} We start by showing that $\overline{\delta}$ is a $D$-structure on $\overline{N}$. We will often use the following relation
\begin{equation}\label{eqn:dIdent}
\sum_{j} a_{1j}a_{j2} = \sum_{j \geq 3} a_{1j}a_{j2}  = 0
\end{equation}
That the first sum is zero follows from $d(a_{12}) = d(1_{A}) = 0$ and the structure identity \ref{eqn:structure}. Since $a_{11} = a_{22} = 0$, the sum can be truncated to $j \geq 3$. \\
\ \\
\noindent To verify that $\overline{\delta}$ is a $D$-structure we will verify \ref{eqn:structure} for $a_{ij}' = a_{ij} - a_{i2}a_{1j}$. First
\begin{equation}
\begin{split}
(-1)^{|x_{j'}|} d(a_{ij} - a_{i2}a_{1j}) &= (-1)^{|x_{j'}|}\big(d(a_{ij}) - (-1)^{|a_{1j}|} d(a_{i2})a_{1j} -  a_{i2}d(a_{1j})\big)\\
& = -\sum_{k=1}^{n} a_{ik}a_{kj}  - (-1)^{|x_{j}|+|a_{1j}|}d(a_{i2})a_{1j} - a_{i2}\big((-1)^{|x_{j}|}d(a_{1j})\big)\\
& = -\sum_{k=1}^{n} a_{ik}a_{kj}  - (-1)^{|x_{2}|}d(a_{i2})a_{1j} + a_{i2}\sum_{k=1}^{n} a_{1k}a_{kj}\\
\end{split}
\end{equation}
since $a_{1j}\otimes x_{j}$ is in the image of $\delta(x_{1})$ we have $|a_{1j}| + |x_{j}| = |x_{1}| + 1$. Since $a_{12} = 1_{A}$ we have that $|x_{2}| = |x_{1}| + 1$, by the same identity. Continuing the sequence of equations:
\begin{equation}
\begin{split}
(-1)^{|x_{j'}|} d(a_{ij} - a_{i2}a_{1j}) &= -\sum_{k=1}^{n} a_{ik}a_{kj}  + \sum_{k=1}^{n} a_{ik}a_{k2}a_{1j} + \sum_{k=1}^{n} a_{i2}a_{1k}a_{kj}\\
& = - a_{i1}a_{1j} - a_{i2}a_{2j} + a_{i1}a_{12}a_{1j} + a_{i2}a_{22}a_{1j} \\
&  \hspace{0.5in} +  a_{i2}a_{11}a_{1j} + a_{i2}a_{12}a_{2j}  \\
& \hspace{0.5in} + \sum_{k=3}^{n} \big[- a_{ik}a_{kj} + a_{ik}a_{k2}a_{1j} + a_{i2}a_{1k}a_{kj}\big]\\
& = - a_{i1}a_{1j} - a_{i2}a_{2j} + a_{i1}a_{1j} + 0 + 0  + a_{i2}a_{2j} \\
& \hspace{0.5in} + \sum_{k=3}^{n} \big[- a_{ik}a_{kj} + a_{ik}a_{k2}a_{1j} + a_{i2}a_{1k}a_{kj}\big]\\
& = \sum_{k=3}^{n} \big[- a_{ik}a_{kj} + a_{ik}a_{k2}a_{1j} + a_{i2}a_{1k}a_{kj}\big]
\end{split}
\end{equation}
where we use $a_{11} = a_{22} = 0$ and $a_{12} = 1_{A}$ to simplify the sum. On the other hand, we can use equation \ref{eqn:dIdent} to simplify
\begin{equation}
\begin{split}
\sum_{k=3}^{n} &\big(a_{ik} - a_{i2}a_{1k}\big)\big(a_{kj} - a_{k2}a_{1j}\big) \\
&= \sum_{k=3}^{n}a_{ik}a_{kj}  - \sum_{k=3}^{n}a_{ik}a_{k2}a_{1j} - \sum_{k=3}^{n} a_{i2}a_{1k}a_{kj} + \sum_{k=3}^{n}a_{i2}a_{1k}a_{k2}a_{1j} \\
&= \sum_{k=3}^{n}\big[a_{ik}a_{kj}  - a_{ik}a_{k2}a_{1j} - a_{i2}a_{1k}a_{kj}\big] + a_{i2}\left(\sum_{k=3}^{n}a_{1k}a_{k2}\right)a_{1j}\\
&= \sum_{k=3}^{n}\big[a_{ik}a_{kj}  - a_{ik}a_{k2}a_{1j} - a_{i2}a_{1k}a_{kj}\big]
\end{split}
\end{equation}
Since adding the results of these two computations will give $0$, the $a_{ij}'$ coefficients satisfy equation \ref{eqn:structure} and $\overline{\delta}$ is a type $D$ structure. \\
\ \\
\noindent We now show that $\pi$ and $\iota$ are $D$-structure morphisms. We start with $\iota$. We need to verify that
$$
(\mu_{A} \otimes \I)\,(\I \otimes \delta_{N})\,\iota - (\mu_{A} \otimes \I)\,(\I \otimes \iota)\,\delta_{\overline{N}}  +  (d \otimes |\I|)\,\iota \equiv 0 
$$
First,
\begin{equation}
\begin{split}
(\mu_{A} \otimes \I)\,(\I \otimes \iota)\,\delta_{\overline{N}}(\overline{x}_{i}) &= (\mu_{A} \otimes \I)\left(\sum_{j \geq 3} (a_{ij} - a_{i2}a_{1j}) \otimes (1_{A} \otimes x_{j} - a_{j2} \otimes x_{1})\right) \\
&=  \sum_{j \geq 3} (a_{ij} - a_{i2}a_{1j}) \otimes  x_{j} -  \left(\sum_{j \geq 3} (a_{ij}a_{j2} - a_{i2}a_{1j}a_{j2})\right) \otimes x_{1}\\
&=  \sum_{j \geq 3} (a_{ij} - a_{i2}a_{1j}) \otimes  x_{j} -  \left(\sum_{j \geq 3} a_{ij}a_{j2} \right) \otimes x_{1} \\
& \hspace{.5in} + a_{i2}\left( \sum_{j \geq 3} a_{1j}a_{j2}\right) \otimes x_{1}\\
&= \sum_{j \geq 3} (a_{ij} - a_{i2}a_{1j}) \otimes  x_{j} +  \big((-1)^{|x_{2}|}d(a_{i2}) + a_{i1}\big)\otimes x_{1}  + 0
\end{split}
\end{equation}
Furthermore, $(d\otimes |\I|)\,\iota(\overline{x}_{i}) = (d \otimes |\I|)(1_{A} \otimes x_{i} - a_{i2}\otimes x_{1}) = -(-1)^{|x_{1}|}d(a_{i2}) \otimes x_{1}$ $= (-1)^{|x_{2}|}d(a_{12})\otimes x_{1}$, so
$$
- (\mu_{A} \otimes \I)\,(\I \otimes \iota)\,\delta_{\overline{N}}  + (d \otimes |\I|)\,\iota(\overline{x}_{i}) = -\sum_{j \geq 3} (a_{ij} - a_{i2}a_{1j}) \otimes  x_{j} - a_{i1} \otimes x_{1} 
$$
On the other hand,
\begin{equation}
\begin{split}
(\mu_{A} \otimes \I)\,(\I \otimes \delta_{N})\,\iota(\overline{x}_{i}) &= (\mu_{A} \otimes \I)\,\left(1_{A} \otimes \big(\delta_{N}(x_{i})\big) -  a_{i2}\otimes\delta_{N}(x_{1})\right)\\
&= \sum_{j} a_{ij} \otimes x_{j} - \big(a_{i2}\sum_{j} a_{1j} \otimes x_{j} \big)\\
&= \sum_{j} \big(a_{ij} - a_{i2}a_{1j}\big) \otimes x_{j} \\
&= (a_{i1} - a_{i2}a_{11}) \otimes x_{1} + (a_{i2} -  a_{i2}a_{12}) \otimes x_{2} + \sum_{j \geq 3} \big(a_{ij} -  a_{i2}a_{1j}\big) \otimes x_{j} \\
&= a_{i1} \otimes x_{1} + (a_{i2} - a_{i2}) \otimes x_{2} + \sum_{j \geq 3} \big(a_{ij} -  a_{i2}a_{1j}\big) \otimes x_{j} \\
&= a_{i1} \otimes x_{1}  + \sum_{j \geq 3} \big(a_{ij} - a_{i2}a_{1j}\big) \otimes x_{j}
\end{split}
\end{equation}
Adding this to the result for $- (\mu_{A} \otimes \I)\,(\I \otimes \iota)\,\delta_{\overline{N}} + (d \otimes |\I|)\,\iota(\overline{x}_{i})$ will give $0$, verifying that $\iota$ is a $D$-structure morphism.\\
\ \\
\noindent For $\pi: N \longrightarrow A \otimes \overline{N}$ the verification is similar and will be shortened. First
\begin{equation}
(\mu_{A} \otimes \I)\,(\I \otimes \delta_{\overline{N}})\,\pi(x_{i}) = \left\{
\begin{array}{ll}
0 & i = 1\\
-\sum_{j\geq 3} \sum_{k \geq 3} \big(a_{1j}a_{jk} - a_{1j}a_{j2}a_{1k}\big) \otimes \overline{x}_{k} & i = 2\\
\sum_{k\geq 3} \big(a_{ik} - a_{i2}a_{1k}\big)\otimes \overline{x}_{k} & i \geq 3
\end{array}
\right.
\end{equation}
However,
\begin{equation}
\begin{split}
\sum_{j\geq 3} \sum_{k \geq 3} \big(a_{1j}a_{jk} - a_{1j}a_{j2}a_{1k}\big) \otimes \overline{x}_{k} &= 
\sum_{k \geq 3} \left(\sum_{j\geq 3} a_{1j}a_{jk} - \left[\sum_{j\geq 3}a_{1j}a_{j2}\right]a_{1k}\right) \otimes \overline{x}_{k}\\
&= \sum_{k\geq 3} \big(-(-1)^{|x_{k}|}(d(a_{1k})) - a_{11}a_{1k} - a_{12}a_{2k} \big) \otimes x_{k}\\
&= \sum_{k\geq 3} \big(-(-1)^{|x_{k}|}(d(a_{1k})) - a_{2k} \big) \otimes x_{k}
\end{split}
\end{equation}
Consequently,
\begin{equation}
(\mu_{A} \otimes \I)\,(\I \otimes \delta_{\overline{N}})\,\pi(x_{i}) = \left\{
\begin{array}{ll}
0 & i = 1\\
\sum_{k\geq 3} \big(a_{2k} +(-1)^{|x_{k}|}(d(a_{1k})) \big) \otimes x_{k} & i = 2\\
\sum_{k\geq 3} \big(a_{ik} - a_{i2}a_{1k}\big)\otimes \overline{x}_{k} & i \geq 3
\end{array}
\right.
\end{equation}
On the other hand,
\begin{equation}
(\mu_{A} \otimes \I)\,(\I \otimes \pi)\,\delta_{N}(x_{i}) = \sum_{k \geq 3} \big(a_{ik} - a_{i2}a_{1k}\big)\otimes \overline{x}_{k}
\end{equation}
Subtracting from the previous result calculates $(\mu_{A} \otimes \I)\,(\I \otimes \delta_{\overline{N}})\,\pi(x_{i}) - 
(\mu_{A} \otimes \I)\,(\I \otimes \pi)\,\delta_{N}(x_{i})$ to be
\begin{equation}
\left\{
\begin{array}{ll}
-\sum_{k \geq 3} \big(a_{1k} - a_{12}a_{1k}\big) & i = 1\\
\sum_{k\geq 3} \big((-1)^{|x_{k}|}(d(a_{1k})) + a_{2k} - a_{2k} + a_{22}a_{1k}\big) \otimes x_{k} & i = 2\\
0 & i \geq 3 \\
\end{array}
\right\}
= 
\left\{
\begin{array}{ll}
0 & i = 1\\
\sum_{j\geq 3} (-1)^{|x_{k}|}(d(a_{1k})) \otimes x_{k} & i = 2\\
0 & i \geq 3
\end{array}
\right.
\end{equation}
However, 
$$
(d \otimes |\I|)\pi(x_{i}) = \left\{
\begin{array}{ll}
0 & i = 1\\
-\sum_{j\geq 3} (d(a_{1k})) \otimes (-1)^{|x_{k}|}x_{k} & i = 2\\
0 & i \geq 3
\end{array}
\right.
$$ 
since $(d\otimes |\I|)\big(\I \otimes \overline{x}_{j}\big) = 0$. Thus, $\pi$ is a $D$-structure morphism:
$$
(\mu_{A} \otimes \I)\,(\I \otimes \delta_{\overline{N}})\,\pi(x_{i}) - 
(\mu_{A} \otimes \I)\,(\I \otimes \pi)\,\delta_{N}(x_{i}) +(d \otimes |\I|)\pi(x_{i}) = 0
$$

\noindent It remains to see that $\pi \ast \iota = I_{\overline{N}}$ and $\iota \ast \pi \simeq_{H} I_{N}$. The first is a simple calculation. For $i \geq 3$
$$
\overline{x}_{i} \stackrel{\iota}{\longrightarrow} 1_{A} \otimes x_{i} - a_{i2}\otimes x_{1} \stackrel{\I \otimes \pi}{\longrightarrow} 1_{A} \otimes 1_{A} \otimes \overline{x}_{i} + 0 \stackrel{\mu_{A} \otimes \I}{\longrightarrow} 1_{A} \otimes \overline{x}_{i}
$$
which is just the $D$-structure morphism $I_{\overline{N}}$. The argument that $\iota \ast \pi \simeq_{H} I_{N}$ is somewhat longer, so we provide the high points:
$$
\iota \ast \pi - I_{N} = (\mu \otimes \I)\,(\I \otimes \iota)\,\pi - \I_{N} =
\left\{
\begin{array}{l} 
x_{1} \longrightarrow -1_{A}\otimes x_{1}\\
x_{2} \longrightarrow -1_{A} \otimes x_{2} - \sum_{j \geq 3} a_{1j} \otimes x_{j}\\
x_{i} \longrightarrow - a_{i2}\otimes x_{1}, \hspace{0.5in} i \in \{3,4,\ldots, n\}
\end{array}
\right.
$$
On the other hand
$$
\begin{aligned}
(\mu\otimes I)(\I \otimes H) \delta(x_{i}) &= - a_{i2} \otimes x_{1}\\
(\mu\otimes I)(\I \otimes \delta) H(x_{i}) &= {\left\{\begin{array}{ll} 0 & i \neq 2\\ -\delta(x_{1}) = -1_{A} \otimes x_{2} -\sum_{j\geq 3}a_{1j} \otimes x_{j} & i = 2\\  \end{array} \right.}\\
(d\otimes |I|) H(x_{i}) &= 0\\
\end{aligned}
$$
When we add these we get $\iota \ast \pi - I_{N}$. For $i=1$ only the first line contributes, but since $a_{12} = 1_{A}$ we obtain equality. For $i=2$, only the second line contributes to the sum, since $a_{22} = 0$. For $i > 2$ only the first line contributes, and consequently the sum is equal to  $\iota \ast \pi - I_{N}$. Therefore, $\iota \ast \pi \simeq_{H} I_{N}$. $\Diamond$\\
\ \\

\end{document}